\documentclass[11pt,reqno]{amsart}

\usepackage{texdraw,multicol,rotating}
\usepackage{amsmath}
\usepackage{amssymb}
\usepackage{epsfig}
\usepackage{amssymb,amsthm,amsfonts,amscd}
\input{txdtools.tex}
\theoremstyle{plain}
\newtheorem{thm}{Theorem}[section]
\newtheorem{prop}[thm]{Proposition}

\theoremstyle{definition}
\newtheorem{defi}[thm]{Definition}

\newtheorem{example}[thm]{Example}
\theoremstyle{remark}
\newtheorem{remark}[thm]{Remark}
\numberwithin{equation}{section}
\newbox{\tmpa}
\newbox{\tmpb}

\newcommand{\nc}{\newcommand}
\nc{\Uq}{U_q} \nc{\Z}{\mathbf{Z}} \nc{\C}{\mathbf{C}}
\nc{\Q}{\mathbf{Q}}
\nc{\op}{\oplus} \nc{\ot}{\otimes} \nc{\pv}{P^{\vee}}
\nc{\ali}{\alpha_i} \nc{\B}{\mathbf{B}} \nc{\F}{\mathbf{F}}
\nc{\bP}{\mathbf{P}} \nc{\V}{\mathbf{V}} \nc{\La}{\Lambda}
\nc{\la}{\lambda}
\nc{\nbinom}[2]{\genfrac{}{}{0pt}{1}{{#1}}{{#2}}}
\nc{\qbinom}[2]{\left[\genfrac{}{}{0pt}{1}{{#1}}{{#2}}\right]}
\nc{\path}{\mathcal{P}} \nc{\fit}{\tilde{f}_i}
\nc{\eit}{\tilde{e}_i} \nc{\Y}{\mathbf{Y}} \nc{\A}{\mathbf{A}}
\nc{\ra}{\rightarrow} \nc{\vep}{\varepsilon} \nc{\vphi}{\varphi}
\nc{\g}{\mathfrak{g}} \nc{\h}{\mathfrak{h}} \nc{\oP}{\overline{P}}
\nc{\pathp}{\mathbf{p}}
\nc{\tris}{ \bsegment \move(0 0)\lvec(10 0)\lvec(10 10)\lvec(0
0)\ifill f:0.7 \esegment } \nc{\recs}{ \bsegment \move(0
0)\lvec(10 0)\lvec(10 5)\lvec(0 5)\lvec(0 0)\ifill f:0.7 \esegment
}
\nc{\hcvec}[5]{%
\getpos(#1 #3)\spx\spy \getpos(#2 #3)\epx\epy \getpos(#4
#5)\xoff\yoff \realadd \spx \xoff \twox \realadd \epx {-\xoff}
\thrx \realadd \spy \yoff \posy \move({\spx} {\spy}) \clvec
({\twox} {\posy})({\thrx} {\posy})({\epx} {\epy}) \rmove(0 0) }
\nc{\ahead}[2]{%
\cossin (0 0)({#1} {#2})\cosa\sina \bsegment
  \drawdim in \setunitscale 0.065
  \realmult {-0.5} \cosa \hcosa
  \realmult {-0.5} \sina \hsina
  \move({\hcosa} {\hsina}) \ravec({\cosa} {\sina})
\esegment }
\nc{\boxi}{%
{%
\savebox{\tmppic}{\begin{texdraw} \small \drawdim em \textref h:C
v:C \setunitscale 0.55 \htext(0 0){$i$} \move(-1 -1)\lvec(-1
1)\lvec(1 1)\lvec(1 -1)\lvec(-1 -1)
\end{texdraw}}%
\raisebox{-0.19\height}{\usebox{\tmppic}}%
}%
}
\nc{\boxj}{%
{%
\savebox{\tmppic}{\begin{texdraw} \small \drawdim em \textref h:C
v:C \setunitscale 0.55 \htext(0 0.1){$j$} \move(-1 -1)\lvec(-1
1)\lvec(1 1)\lvec(1 -1)\lvec(-1 -1)
\end{texdraw}}%
\raisebox{-0.19\height}{\usebox{\tmppic}}%
}%
}
\nc{\boxipo}{%
{%
\savebox{\tmppic}{\begin{texdraw} \small \drawdim em \textref h:C
v:C \setunitscale 0.55 \htext(0.15 0){$i\!\!+\!\!1$} \move(-1.4
-1)\lvec(-1.4 1)\lvec(1.4 1)\lvec(1.4 -1)\lvec(-1.4 -1)
\end{texdraw}}%
\raisebox{-0.19\height}{\usebox{\tmppic}}%
}%
} \everytexdraw{ \drawdim in \arrowheadsize l:0.065 w:0.03
\arrowheadtype t:F \textref h:C v:C }
\newsavebox{\tmppic}
\newsavebox{\tmpfig}
\newsavebox{\tmpdraw}
\newsavebox{\tmpfiga}
\newsavebox{\tmpfigb}
\newsavebox{\tmpfigc}
\newsavebox{\tmpfigd}
\newsavebox{\tmpfige}
\newsavebox{\tmpfigf}
\newsavebox{\tmpfigg}
\newsavebox{\tmpfigh}
\newsavebox{\tmpfigi}
\newsavebox{\tmpfigj}
\newsavebox{\tmpfigk}
\newsavebox{\tmpfigl}
\newsavebox{\tmpfigm}
\newsavebox{\tmpfign}
\newsavebox{\tmpfigo}
\newsavebox{\tmpfigp}
\newsavebox{\tmpfigq}
\newsavebox{\tmpfigr}
\newsavebox{\tmpfigs}
\newsavebox{\tmpfigt}
\newsavebox{\tmpfigu}
\newsavebox{\tmpfigv}
\newsavebox{\tmpfigw}
\newsavebox{\tmpfigx}
\newsavebox{\tmpfigy}
\newsavebox{\tmpfigz}
\newsavebox{\tmpfigaa}
\newsavebox{\tmpfigab}
\newsavebox{\tmpfigac}
\newsavebox{\tmpfigad}
\newsavebox{\tmpfigae}
\newsavebox{\tmpfigaf}
\newsavebox{\tmpfigag}
\newsavebox{\tmpfigah}
\newsavebox{\tmpfigai}
\newsavebox{\tmpfigaj}
\newsavebox{\tmpfigak}
\newsavebox{\tmpfigal}
\newsavebox{\tmpfigam}
\newsavebox{\tmpfigan}
\newsavebox{\tmpfigao}
\newsavebox{\tmpfigap}
\newsavebox{\tmpfigaq}
\newsavebox{\tmpfigar}
\newsavebox{\tmpfigas}
\newsavebox{\tmpfigat}
\newsavebox{\tmpfigau}
\newsavebox{\tmpfigav}
\newsavebox{\tmpfigaw}
\newsavebox{\tmpfigax}
\newsavebox{\tmpfigay}
\newsavebox{\tmpfigaz}
\newsavebox{\tmpfigba}
\newsavebox{\tmpfigbb}
\newsavebox{\tmpfigbc}
\newsavebox{\tmpfigbd}
\newsavebox{\tmpfigbe}
\newsavebox{\tmpfigbf}
\newsavebox{\tmpfigbg}
\newsavebox{\tmpfigbh}

\newsavebox{\spinaa}
\newsavebox{\spinab}
\newsavebox{\spinac}
\newsavebox{\spinad}
\newsavebox{\spinae}
\newsavebox{\spinaf}
\newsavebox{\spinag}
\newsavebox{\spinah}
\newsavebox{\spinai}
\newsavebox{\spinaj}
\newsavebox{\spinak}
\newsavebox{\spinal}
\newsavebox{\spinam}
\newsavebox{\spinba}
\newsavebox{\spinbb}
\newsavebox{\spinbc}
\newsavebox{\spinbd}
\newsavebox{\spinbe}
\newsavebox{\spinbf}
\newsavebox{\spinbg}
\newsavebox{\spinbh}
\nc{\node}{\lcir r:1 }
\nc{\sline}{\bsegment\savepos(10 0)(*ex *ey)
            \move(1 0)\rlvec(8 0)
            \esegment\move(*ex *ey)}
\nc{\dline}{\bsegment\savepos(10 0)(*ex *ey)
            \move(0.93 0.4)\rlvec(8.14 0)\rmove(0 -0.8)\rlvec(-8.14 0)
            \esegment\move(*ex *ey)}
\nc{\uline}{\bsegment\savepos(0 10)(*ex *ey)
            \move(0 1)\rlvec(0 8)
            \esegment\move(*ex *ey)}
\nc{\lpoint}{\savecurrpos(*ex *ey)
             \rmove(2.5 1.5)\rlvec(-1.5 -1.5)\rlvec(1.5 -1.5)
             \move(*ex *ey)}
\nc{\rpoint}{\savecurrpos(*ex *ey)
             \rmove(-2.5 -1.5)\rlvec(1.5 1.5)\rlvec(-1.5 1.5)
             \move(*ex *ey)}
\nc{\bline}{\bsegment\savepos(10 0)(*ex *ey)
            \linewd 0.6 \move(1.1 0)\rlvec(7.8 0)
            \esegment\move(*ex *ey)}
\nc{\araise}[1]{\raisebox{4.5pt}{#1}}
\nc{\braise}[1]{\raisebox{12.1pt}{#1}}
\nc{\craise}[1]{\raisebox{8pt}{#1}}
\nc{\draise}[1]{\raisebox{12pt}{#1}}
\nc{\eraise}[1]{\raisebox{14.5pt}{#1}}

\begin{document}

\title[]
{Correspondence between Young walls and Young tableaux
realizations of crystal bases for the classical Lie algebras}

\author[Jeong-Ah Kim and Dong-Uy Shin]
{Jeong-Ah Kim$^\diamond$ and Dong-Uy Shin$^{*}$}

\address{$^\diamond$Department of Mathematics\\
         Seoul National University\\
         Seoul 151-747, Korea}
\email{jakim@math.snu.ac.kr}
\address{$^*$School of Mathematics\\
         Korea Institute for Advanced Study\\
         Seoul 130-012, Korea}
\email{shindong@kias.re.kr}

\begin{abstract}
We give a 1-1 correspondence with the Young wall realization and
the Young tableau realization of the crystal bases for the
classical Lie algebras.
\end{abstract}
\maketitle
\vskip 1cm

\section*{Introduction}
{\it Young tableaux} and {\it Young walls} play important roles in
the interplay, which can be explained in a beautiful manner using
the {\it crystal base theory} for quantum groups, between the
fields of representation theory and combinatorics. Indeed, since
representation theory is known to be a vital tool in the solution
of certain kinds of two-dimensional solvable lattice models in
statistical mechanics, Young tableaux and Young walls are central
ingredients in mathematical physics.

\vskip 2mm The {\it quantum groups} introduced by Drinfeld and
Jimbo, independently are deformations of the universal enveloping
algebras of Kac-Moody algebras \cite{Drin,Jim}. More precisely,
let $\frak g$ be a Kac-Moody algebra and $U(\frak g)$ be its
universal enveloping algebra. Then, for each generic parameter
$q$, we associate a Hopf algebra $U_q(\frak g)$, called the
quantum group, whose structure tends to that of $U(\frak g)$ as
$q$ approaches $1$.

The important feature of quantum groups is that the representation
theory of $U(\frak g)$ is the same as that of $U_q(\frak g)$. More
precisely, let $M$ be a $U(\frak g)$-module in the category
${\mathcal O}_{int}$, which has a weight space decomposition
$M=\oplus_{\la\in P}M_{\la}$. Then, for each generic $q$, there
exists a $U_q(\frak g)$-module $M^q$ in  the category ${\mathcal
O}^q_{int}$ with a weight space decomposition $M^q=\oplus_{\la\in
P}M^q_{\la} $ such that $\text{dim}_{{\bf
C}(q)}M^q_{\la}=\text{dim}_{{\bf C}}M_{\la}$ for all $\la\in P$
and the structure of $M^q$ tends to that of $M$ as $q$ approaches
$1$. Therefore, to understand the structure of representations
over general quantum groups $U_q(\frak g)$, it is enough to
understand that of representations over $U_q(\frak g)$  for some
special parameter $q$ which is easy to treat.

\vskip 2mm  The {\it crystal bases}, introduced by Kashiwara
\cite{Kas90,Kas91}, can be viewed as bases at $q=0$ for the
integrable modules over quantum groups and they are given a
structure of colored oriented graph, called the {\it crystal
graphs}, which reflect the combinatorial structure of integrable
modules.

For classical Lie algebras, Kashiwara and Nakashima gave an
explicit realization of crystal bases for finite dimensional
irreducible modules \cite{KasNak}. In their work, crystal bases
were characterized as the sets of semistandard Young tableaux with
given shapes satisfying certain additional conditions. In
\cite{Li1}, Littelmann gave another description of crystal bases
for finite dimensional simple Lie algebras using the
Lakshmibai-Seshadri monomial theory. His approach was generalized
to the {\it path model theory} for all symmetrizable Kac-Moody
algebras \cite{Li2}. Littelmann's theory also gives rise to
colored oriented graphs, which turned out to be isomorphic to the
crystal graphs \cite{Kas96}. Moreover, in \cite{Nak}, Nakashima
gave a generalization of the Littlewood-Richardson rule for
$U_q(\frak g)$ associated with the classical Lie algebras using
the crystal base theory.

\vskip 2mm In \cite{HK} and \cite{Kang}, Young wall, an affine
combinatorial object, was introduced. The crystal bases for basic
representations (i.e., highest weight representation of level $1$)
for quantum affine algebras are realized as the sets of {\it
reduced proper Young walls}. Motivated by the fact that a
classical Lie algebra $\frak g$ lies inside an affine Lie algebra
$\hat{\frak g}$ and any crystal graph $B(\lambda)$ for a finite
dimensional irreducible $\mathfrak g$-module $V(\lambda)$ appears
as a connected component in the crystal graph $B(\Lambda)$ of a
basic representation $V(\Lambda)$ over $\hat{\frak g}$ without
$0$-arrows, Kang, Kim, Lee and Shin gave a new realization of
crystal bases for finite dimensional modules over classical Lie
algebras in terms of Young walls corresponding to the connected
component having the least number of blocks \cite{KKLS}.

\vskip 2mm In this paper, we give a new realization of crystal
bases for finite dimensional irreducible modules over classical
Lie algebras derived from the Young wall realization. The basis
vector is parameterized by certain tableaux with given shape which
is different from generalized Young diagram given by Kashiwara and
Nakashima. Moreover, motivated by the fact that the tableau
realization of $B(\la)$ given by Kashiwara and Nakashima can be
derived from the Young wall realization of the crystal graph
$B(\la)$ corresponding to the connected component having the
largest number of blocks, we give a crystal isomorphism from our
new realization to the realization of Kashiwara and Nakashima,
which will be a part of insertion scheme for the crystal of the
classical Lie algebras given in \cite{KS}.

\vskip 2mm The contents of this paper is organized as follows: In
Section 1, we review the crystal base theory and the tableaux
realization of crystal bases for the classical Lie algebras given
by Kashiwara and Nakashima. In Section 2, we recall the notions
and properties of Young walls and the Young wall realization of
crystal bases for classical Lie algebras given in \cite{KKLS}. In
Section 3, we give a new realization of crystal bases in terms of
new tableaux, which are obtained from Young wall realization. In
the end of Section 3, we give a 1-1 correspondence between our new
tableau realization and the tableau realization given by Kashiwara
and Nakashima, which turns out to be a crystal isomorphism.

\vskip 1cm
\section{Quantum groups and crystal bases}
Let $I$ be an index set and let $(A, \Pi, \Pi^{\vee}, P,
P^{\vee})$ be a Cartan datum, i.e.,

\begin{itemize}
\item [(a)] $A=(a_{ij})_{i,j \in I}$ be a {\it generalized Cartan
matrix},
\item [(b)] $\pv=\left(\bigoplus_{i \in I} \Z h_i\right) \bigoplus
\left(\bigoplus_{j=1}^{|I|-\text{rank}A}\Z d_j\right)$ is the {\it
dual weight lattice},
\item [(c)] ${\mathfrak h}= \C \otimes_{\Z} \pv$ is the {\it Cartan
subalgebra},
\item [(d)] $P=\{\lambda\in \mathfrak h^* | \lambda(P^{\vee})\subset \Z
\}$ is the {\it weight lattice},
\item [(e)] $\Pi^{\vee} = \{ h_i | \, i\in I\}$ is the set of {\it
simple coroots},
\item [(f)] $\Pi=\{\alpha_i | \, i\in I\}\subset \mathfrak
h^* $ is the set of {\it simple roots} defined by $\alpha_i (h_j)
= a_{ji}$, $\alpha_i (d_s) = 0$ or $1$ for $i,j \in I,\,
s=1,\cdots,|I|-\text{rank}A$.
\end{itemize}

To each Cartan datum, we can associate two algebras $\frak{g}$ and
$U_q(\frak g)$, called the {\it Kac-Moody algebra}
 and {\it quantized universal enveloping algebra} \cite{Kac90,HK}.
The generalized Cartan matrices are classified into three
categories: {\it finite} type, {\it affine} type and {\it
indefinite} type. (See \cite{Kac90} for more details.) The
indecomposable generalized Cartan matrices of finite type give
rise to {\it finite dimensional simple Lie algebras} and those of
affine type yield infinite dimensional Kac-Moody algebras called
the {\it affine Lie algebras}. The corresponding quantum groups
will be called the {\it quantum classical algebra} and the {\it
quantum affine algebra}, respectively.

\vskip 2mm Let $P^+=\{\,\la\in P\mid\la(h_i)\in {\bf Z}_{\ge
0}\,\,\,\text{for any $i\in I$}\}$ be the set of {\it dominant
integral weights} and let $\La_i\in P^+$ be the {\it fundamental
weight} defined by $\La_i(h_j)=\delta_{ij}$. From now on, we
denote by $\la_i$ ($i=1,\cdots,n$) and $\La_i$ ($i=0,1,\cdots,n$)
the fundamental weights for the quantum classical algebra and the
quantum affine algebra, respectively.

\vskip 2mm The {\it crystal base theory} is developed for
$U_q(\g)$-modules $M$ in the category ${\mathcal O}^q_{int}$
consisting of $U_q(\g)$-modules $M$ such that
\begin{itemize}
\item [(a)] $M$ has a weight space decomposition,
\item [(b)] there exist a finite number of elements
$\la_1$,$\cdots$,$\la_s\in P$ such that
\begin{center}
$\text{wt}(M)\subset D(\la_1)\cup\cdots\cup D(\la_s),$
\end{center}
where $\text{wt}(M)$ is the set of all weights of $M$ and
$D(\la)=\{\mu\in \frak h^*| \mu\le \la\}$,
\item [(c)] all $e_i$ and $f_i$ ($i\in I$) are {\it locally nilpotent} on $M$.
\end{itemize}

\vskip 2mm By the representation theory of $\Uq(\mathfrak{sl}_2)$,
every element $v\in M_\la$ can  be written uniquely as $v =
\sum_{k\geq0} f_i^{(k)} v_k$, where $k \ge -\la(h_i)$, $f_i^{(n)}
= f_i^n/[n]_i!$ and $v_k \in \ker e_i \cap M_{\la+k\ali}$. Then
the {\it Kashiwara operators} $\eit$ and $\fit$ on $M$ are defined
by
\begin{center}
$\eit v = \sum_{k\geq1} f_i^{(k-1)} v_k, \qquad \fit v =
\sum_{k\geq0} f_i^{(k+1)} v_k.$
\end{center}

\vskip 2mm Let $M$ be a $U_q(\g)$-module in the category
${\mathcal O}^q_{int}$ and let $ \A = \{ f/g \in \C(q)  \, | \, f,
g \in \C[q], \, g(0)\neq 0 \} $ be the subring of $\C(q)$
consisting of the rational functions in $q$ that are regular at
$q=0$.

\begin{defi}
A pair $(L, B)$ is called a {\it crystal base} if the following
conditions are satisfied:
\begin{itemize}
\item [(a)] $L$ is a free $\A$-submodule of $M$ such that $M \cong \C(q)
\otimes_{\A} L$,
\item [(b)] $B$ is a $\C$-basis of $L/qL$,
\item [(c)] $L=\bigoplus_{\lambda \in P} L_{\la}$, $B=\bigsqcup_{\la \in
P} B_{\la}$, where $L_{\la} = L \cap M_{\la}$ and $B_{\la}=B \cap
\left(L_{\la} / q L_{\la} \right)$,
\item [(d)] $\eit L \subset L$, $\fit L \subset L$ for all $i\in I$,
\item [(e)] $\eit B\subset B\cup \{0\}$, \ $\fit B\subset B\cup \{0\}$,
\item [(f)] For $b, b'\in B$, $\fit b = b'$ if and only if $b = \eit b'$.
\end{itemize}
\end{defi}

\vskip 3mm  For each $b\in B$ and $i\in I$, we define
\begin{center}
$\varepsilon_i(b)=\max \{k\ge 0 | \, \eit^k b \in B \}, \quad
\varphi_i(b) = \max \{ k\ge 0 | \, \fit^k b \in B \}.$
\end{center}
Let $M_j$ be a $U_q(\g)$-module in the category ${\mathcal
O}^q_{int}$ with a crystal basis $(L_j, B_j)$ \ $(j=1, \cdots,
N)$. Fix an index $i\in I$ and consider a vector $b=b_1 \ot \cdots
\ot b_N \in B_1 \ot \cdots \ot B_N$. To each $b_j \in B_j$ \,
$(j=1, \cdots , N)$, we assign a sequence of $-$'s and $+$'s with
as many $-$'s as $\vep_i(b_j)$ followed by as many $+$'s as
$\vphi_i(b_j)$:
\begin{equation*}
\begin{aligned}
b \, &  =  \, b_1 \otimes b_2 \otimes \cdots \otimes b_N \\
& \longmapsto
 (\underbrace{-, \cdots, -}_{\vep_i(b_1)},
\underbrace{+, \cdots, +}_{\vphi_i(b_1)}, \cdots \cdots,
\underbrace{-, \cdots, -}_{\vep_i(b_N)}, \underbrace{+, \cdots,
+}_{\vphi_i(b_N)}).
\end{aligned}
\end{equation*}
In this sequence, we cancel out all the $(+,-)$-pairs to obtain a
sequence of $-$'s followed by $+$'s, called the {\it
$i$-signature} of $b$:
\begin{center}
$\text{$i$-sgn}(b) = (-, -, \cdots, -, +, +, \cdots, +).$
\end{center}

\vskip 2mm Then the crystal bases have a nice behavior with
respect to the tensor product.

\vskip 3mm
\begin{prop} {\rm \cite{Kas90, Kas91}}\,
Let $M_j$ be a $U_q(\g)$-module in the category ${\mathcal
O}^q_{int}$ with a crystal basis $(L_j, B_j)$ \ $(j=1, \cdots, N)$
and let $b=b_1 \ot \cdots \ot b_N$ be a vector in $B_1 \ot \cdots
\ot B_N$. Then $\eit$ acts on $b_j$ corresponding to the
right-most $-$ in $\text{$i$-sgn}(b)$ and $\fit$ acts on $b_k$
corresponding to the left-most $+$ in $\text{$i$-sgn}(b)$:
\begin{equation}\label{eq:signature}
\begin{aligned}
& \eit b = b_1 \otimes \cdots \ot \eit b_j \ot \cdots \ot b_N, \\
& \fit b = b_1 \otimes \cdots \ot \fit b_k \ot \cdots \ot b_N.
\end{aligned}
\end{equation}
Moreover, $\eit b=0$ $($resp. $\fit b=0)$ if there is no $-$
$($resp. $+)$ in the $i$-signature of $b$.
\end{prop}

\vskip 2mm

\savebox{\tmpfiga}{\begin{texdraw}\fontsize{7}{7}\selectfont
\drawdim em \setunitscale 0.8

\rlvec(0 2) \rlvec(2 0)\rlvec(0 -2)\rlvec(-2 0)

\htext(1 1){$i$}
\end{texdraw}}
\savebox{\tmpfigb}{\begin{texdraw}\fontsize{7}{7}\selectfont
\drawdim em \setunitscale 0.8

\rlvec(0 2) \rlvec(2 0)\rlvec(0 -2)\rlvec(-2 0)

\htext(1 1){$0$}
\end{texdraw}}
\savebox{\tmpfigc}{\begin{texdraw}\fontsize{7}{7}\selectfont
\drawdim em \setunitscale 0.8

\rlvec(0 2) \rlvec(2 0)\rlvec(0 -2)\rlvec(-2 0)

\htext(1 1){$\overline{i}$}
\end{texdraw}}
\savebox{\tmpfigd}{\begin{texdraw}\fontsize{7}{7}\selectfont
\drawdim em \setunitscale 0.8

\rlvec(0 2) \rlvec(2 0)\rlvec(0 -2)\rlvec(-2 0)

\htext(1 1){$j$}
\end{texdraw}}
\savebox{\tmpfige}{\begin{texdraw}\fontsize{7}{7}\selectfont
\drawdim em \setunitscale 0.8

\rlvec(0 2) \rlvec(2 0)\rlvec(0 -2)\rlvec(-2 0)

\htext(1 1){$\overline{j}$}
\end{texdraw}}
\savebox{\tmpfigf}{\begin{texdraw}\fontsize{7}{7}\selectfont
\drawdim em \setunitscale 0.8

\rlvec(0 2) \rlvec(2 0)\rlvec(0 -2)\rlvec(-2 0)

\htext(1 1){$j\!\!-\!\!1$}
\end{texdraw}}
\savebox{\tmpfigg}{\begin{texdraw}\fontsize{7}{7}\selectfont
\drawdim em \setunitscale 0.8

\rlvec(0 2) \rlvec(2 0)\rlvec(0 -2)\rlvec(-2 0)

\htext(1 1){$j\!\!+\!\!1$}
\end{texdraw}}
\savebox{\tmpfigh}{\begin{texdraw}\fontsize{7}{7}\selectfont
\drawdim em \setunitscale 0.8

\rlvec(0 2) \rlvec(2 0)\rlvec(0 -2)\rlvec(-2 0)

\htext(1 1){$\overline{j\!\!-\!\!1}$}
\end{texdraw}}
\savebox{\tmpfigi}{\begin{texdraw}\fontsize{7}{7}\selectfont
\drawdim em \setunitscale 0.8

\rlvec(0 2) \rlvec(2 0)\rlvec(0 -2)\rlvec(-2 0)

\htext(1 1){$\overline{j\!\!+\!\!1}$}
\end{texdraw}}
\savebox{\tmpfigj}{\begin{texdraw}\fontsize{7}{7}\selectfont
\drawdim em \setunitscale 0.8

\rlvec(0 2) \rlvec(2 0)\rlvec(0 -2)\rlvec(-2 0)

\htext(1 1){$n$}
\end{texdraw}}
\savebox{\tmpfigk}{\begin{texdraw}\fontsize{7}{7}\selectfont
\drawdim em \setunitscale 0.8

\rlvec(0 2) \rlvec(2 0)\rlvec(0 -2)\rlvec(-2 0)

\htext(1 1){$\overline{n}$}
\end{texdraw}}

\begin{example}
The representation $V(\la_1)$ for $U_q(B_n)$, called the {\it
vector representation} and denoted by ${\bf V}$, has a basis
$\{\,\raisebox{-0.25\height}{\usebox{\tmpfiga}}\,,\,\,
\raisebox{-0.25\height}{\usebox{\tmpfigc}}\,; i\in {\mathcal A}\}
\cup\{\,\raisebox{-0.25\height}{\usebox{\tmpfigb}}\,\}$
(${\mathcal A}=\{1,\cdots,n\}$) and the action of generators of
$U_q(B_n)$ is given as follows:
$$
q^h\,
\raisebox{-0.25\height}{\usebox{\tmpfigd}}=q^{\epsilon_j(h)}\,
\raisebox{-0.25\height}{\usebox{\tmpfigd}}\,,\quad
q^h\,\raisebox{-0.25\height}{\usebox{\tmpfige}}=q^{-\epsilon_j(h)}\,
\raisebox{-0.25\height}{\usebox{\tmpfige}}\,,\quad
q^h\,\raisebox{-0.25\height}{\usebox{\tmpfigb}}=
\raisebox{-0.25\height}{\usebox{\tmpfigb}}\,,
$$
$$
\aligned &e_i\,
\raisebox{-0.25\height}{\usebox{\tmpfigd}}=\delta_{i+1,j}\,
\raisebox{-0.25\height}{\usebox{\tmpfigf}}\,,\quad e_i\,
\raisebox{-0.25\height}{\usebox{\tmpfige}}=\delta_{i,j}\,
\raisebox{-0.25\height}{\usebox{\tmpfigi}}\,,\quad e_i\,
\raisebox{-0.25\height}{\usebox{\tmpfigb}}=0\\
&f_i\, \raisebox{-0.25\height}{\usebox{\tmpfigd}}=\delta_{i,j}\,
\raisebox{-0.25\height}{\usebox{\tmpfigg}}\,,\quad f_i\,
\raisebox{-0.25\height}{\usebox{\tmpfige}}=\delta_{i+1,j}\,
\raisebox{-0.25\height}{\usebox{\tmpfigh}}\,,\quad f_i\,
\raisebox{-0.25\height}{\usebox{\tmpfigb}}=0
\endaligned
\quad (1\le i<n, 1\le j\le n),
$$
$$
\aligned &e_n\, \raisebox{-0.25\height}{\usebox{\tmpfigd}}=0,\quad
e_n\, \raisebox{-0.25\height}{\usebox{\tmpfige}}=\delta_{n,j}\,
\raisebox{-0.25\height}{\usebox{\tmpfigb}}\,,\quad e_n\,
\raisebox{-0.25\height}{\usebox{\tmpfigb}}=
[2]_n\,\raisebox{-0.25\height}{\usebox{\tmpfigj}}\\
&f_n\, \raisebox{-0.25\height}{\usebox{\tmpfigd}}=\delta_{n,j}\,
\raisebox{-0.25\height}{\usebox{\tmpfigb}}\,,\quad f_n\,
\raisebox{-0.25\height}{\usebox{\tmpfige}}=0,\quad f_n\,
\raisebox{-0.25\height}{\usebox{\tmpfigb}}=
[2]_n\,\raisebox{-0.25\height}{\usebox{\tmpfigk}}
\endaligned
\quad (1\le j\le n).
$$
Then the crystal base $({\bf L}, {\bf B})$ of ${\bf V}$ is given
by
$$
\aligned &{\bf L}=\bigoplus_{i=1}^n({\bf
A}\raisebox{-0.25\height}{\usebox{\tmpfiga}} \oplus{\bf A}
\raisebox{-0.25\height}{\usebox{\tmpfigc}}\,)\oplus
{\bf A}\raisebox{-0.25\height}{\usebox{\tmpfigb}}\,,\\
&{\bf B}=\{\,\raisebox{-0.25\height}{\usebox{\tmpfiga}}\,,\,\,
\raisebox{-0.25\height}{\usebox{\tmpfigc}}\,; i\in {\mathcal
A}\}\cup \{\,\raisebox{-0.25\height}{\usebox{\tmpfigb}}\,\}.
\endaligned
$$
Moreover, the crystal graph of ${\bf B}=B(\la_1)$ is given by

\vskip 3mm
\begin{center}
\begin{texdraw}\fontsize{7}{7}\selectfont
\drawdim em \setunitscale 0.8

\lvec(2 0)\rlvec(0 2)\rlvec(-2 0)\rlvec(0 -2) \move(2.5 1)
\arrowheadtype t:V  \ravec(2 0)

\move(5 0)\rlvec(2 0)\rlvec(0 2)\rlvec(-2 0)\rlvec(0 -2) \move(7.5
1) \arrowheadtype t:V  \ravec(2 0)

\move(12.5 1) \arrowheadtype t:V  \ravec(2 0)

\move(15 0)\rlvec(2 0)\rlvec(0 2)\rlvec(-2 0)\rlvec(0 -2)
\move(17.5 1) \arrowheadtype t:V  \ravec(2 0)

\move(20 0)\rlvec(2 0)\rlvec(0 2)\rlvec(-2 0)\rlvec(0 -2)
\move(22.5 1) \arrowheadtype t:V  \ravec(2 0)

\move(25 0)\rlvec(2 0)\rlvec(0 2)\rlvec(-2 0)\rlvec(0 -2)
\move(27.5 1) \arrowheadtype t:V  \ravec(2 0)

\move(32.5 1) \arrowheadtype t:V  \ravec(2 0)

\move(35 0)\rlvec(2 0)\rlvec(0 2)\rlvec(-2 0)\rlvec(0 -2)
\move(37.5 1) \arrowheadtype t:V  \ravec(2 0)

\move(40 0)\rlvec(2 0)\rlvec(0 2)\rlvec(-2 0)\rlvec(0 -2)

\htext(1 1){$1$}\htext(6 1){$2$} \htext(11 1){$\cdots$}\htext(16
1){$n$} \htext(21 1){$0$}\htext(26 1){$\overline{n}$}\htext(31
1){$\cdots$}\htext(36 1){$\overline{2}$}\htext(41
1){$\overline{1}$}\htext(42.5 0.2){$\cdot$}

\htext(3.5 1.5){$1$}\htext(8.5 1.5){$2$}\htext(13.5
1.5){$n\!\!-\!\!1$} \htext(18.5 1.5){$n$}\htext(23.5
1.5){$n$}\htext(28.5 1.5){$n\!\!-\!\!1$}\htext(33.5
1.5){$2$}\htext(38.5 1.5){$1$}
\end{texdraw}
\end{center}
\end{example}

\vskip 3mm We close this section with the realization of crystal
bases of irreducible highest weight modules for $U_q(\frak g)$
using Young tableaux, where $\frak g=A_n, C_n, B_n, D_n$. In
particular, we will focus on the classical Lie algebra $\frak
g=B_n$. The description of other types is referred to
\cite{KasNak}.

\vskip 3mm For a sequence of half integers $l_j\in \frac{1}{2}{\bf
Z}$ ($j=1,\cdots,n$), such that $l_j-l_{j+1}\in {\bf Z}_{\ge 0}$,
let $Y=(l_1,\cdots,l_n)$ be a diagram which has $n$ rows with the
$j$-th row of length $|l_j|$. Moreover, if $l_j\in \frac{1}{2}{\bf
Z}_{\ge 0}$ for all $j$, $Y$ is called the {\it generalized Young
diagram} of type $B_n$. Then there is a bijection $GY$ between the
set of generalized Young diagrams and the set of dominant integral
weights of $B_n$ defined by $GY: (l_1,\cdots,l_n)\mapsto
\sum_{j=1}^n l_j\epsilon_j$, where
($\epsilon_1,\cdots,\epsilon_n$) is the orthonormal base of $\frak
h^*$ \cite{Nak}.

\vskip 2mm The crystal graph $B(\la)$ is realized as the set of
tableaux of shape $GY^{-1}(\la)$ on $\{i,\overline{i}\,|\, i\in
{\mathcal A}\}\cup \{0\}$ with the linear order
\begin{center}
$1\prec 2\prec \cdots\prec n\prec 0\prec \overline{n}\prec \cdots
\prec \overline{2}\prec \overline{1}.$
\end{center}
To describe $B(\la)$, we need some definitions and conditions.
\vskip 3mm {\bf (1CC)} Given a column $C$ of length $N$ containing
the entries (reading from top to bottom)
$i_1$,$i_2$,$\cdots$,$i_N$, we say that $C$ satisfies the {\it
one-column condition} (1CC) if for $i_p=a$ and $i_q=\overline{a}$
($1\le a\le n$), $p+(N-q+1)\le a$.
 \vskip 3mm Let $C$ and $C'$ be adjacent
columns of length $N$ and $M$ ($1\le M\le N\le n$) consisting of
the entries (reading from top to bottom)
$i_1$,$i_2$,$\cdots$,$i_N$ and $j_1$,$j_2$,$\cdots$,$j_M$,
respectively. (Note that $C$ can be a half-column.)

\begin{defi}
(a) For $1\le a\le b<n$, we say that $C$ and $C'$ is in the
($a,b$)-{\it configuration} if there exist $1\le p\le q<r\le s\le
M$ such that $(i_p,i_q,i_r,j_s)=(a,b,\overline{b},\overline{a})$
or $(i_p,j_q,j_r,j_s)=(a,b,\overline{b},\overline{a})$.

(b) For $1\le a<n$, we say that $C$ and $C'$ is in the
($a,n$)-{\it configuration} if there exist $1\le p\le q<r=q+1\le
s\le M$ such that $i_p=a$, $j_s=\overline{a}$ and one of the
following conditions is satisfied:
\begin{center}
(i) $i_q$ and $i_r$ are $n$, $0$ or $\overline{n}$, \hskip 5mm
(ii) $j_q$ and $j_r$ are $n$, $0$ or $\overline{n}$.
\end{center}

(c) We say that $C$ and $C'$ is in the ($n,n$)-{\it configuration}
if there exist $1\le p<q \le M$ such that $(i_p,j_q)=(n,0)$,
$(n,\overline{n})$ or $(0,\overline{n})$.
\end{defi}

\vskip 2mm Now, for $C$ and $C'$ in the ($a,b$)-configuration, we
define $p(a,b;C,C')=(q-p)+(s-r)$.  In particular, if $a=b=n$, we
set $p(a,b;C,C')=0$.

\vskip 3mm {\bf (2CC)} Given adjacent two columns $C$ and $C'$, we
say that $C$ and $C'$ satisfy the {\it two-column condition} (2CC)
if for every $(a,b)$-configuration, $p(a,b;C,C')<b-a$.

\vskip 3mm
\begin{prop} \cite{KasNak}
For a dominant integral weight $\la$, the crystal graph $B(\la)$
is realized as the set of tableaux $T$ of shape $GY^{-1}(\la)$
with entries $\{i,\overline{i}\,|\, i\in {\mathcal A}\}\cup\{0\}$
such that
\begin{itemize}
\item [(a)] the entries of $T$ weakly increase along the rows, but the
element $0$ cannot appear more than once,
\item [(b)] the entries of $T$ strictly increase down the columns, but
the element $0$ can appear more than once,
\item [(c)] for a half column $C$ of $T$, $i$ and $\overline{i}$ can not
appear at the same time,
\item [(d)] for each column $C$ of $T$, (1CC) holds,
\item [(e)] for each pair of adjacent columns $C$, $C'$ of $T$, (2CC)
holds.
\end{itemize}
\end{prop}

\vskip 1cm
\section{Crystal bases in terms of Young walls for the classical Lie algebras}
In this section, we will review the realization of crystal bases
of irreducible highest weight modules for the classical Lie
algebras $\frak g=A_n,C_n,B_n$ and $D_n$ in terms of Young walls
given in \cite{KKLS}.

\vskip 2mm Young wall is a combinatorial object for realizing the
crystal bases for quantum affine algebras. They are built of
colored blocks with three different shapes. Given an affine
dominant integral weight $\Lambda$ of level 1 (i.e. $\Lambda(c)=1$
for the canonical central element $c$), we fix a frame
$Y_{\Lambda}$ called the {\it ground state wall} of weight
$\Lambda$ and on this frame, we build a wall of thickness less
than or equal to one unit with the rules for building the walls.
The coloring of blocks, description of ground state walls and the
patterns for building the walls are given in \cite{Kang}.

\vskip 2mm A column in a Young wall is called a {\it full column}
if its height is a multiple of the unit length and its top is of
unit thickness. A Young wall is said to be {\it proper} if none of
the full columns have the same height for the quantum affine
algebras of type $A_{2n-1}^{(2)}$, $B_n^{(1)}$ and $D_n^{(1)}$ and
every Young wall is defined to be {\it proper} for the quantum
affine algebras of type $A_n^{(1)}$.

\vskip 3mm
\begin{defi} (a) A block of color $i$  in a proper Young wall is called
a {\it removable $i$-block} if the wall remains a proper Young
wall after removing the block. A column in a proper Young wall is
called {\it $i$-removable} if the top of that column is a
removable $i$-block.

(b) A place in a proper Young wall where one may add an $i$-block
to obtain another proper Young wall is called an {\it
$i$-admissible slot}. A column in a proper Young wall is called
{\it $i$-admissible} if the top of that column is an
$i$-admissible slot.
\end{defi}

\vskip 2mm Let $\F(\Lambda)$ be the set of all proper Young walls
on $Y_{\Lambda}$ with a affine dominant integral weight $\Lambda$
of level $1$. We now define the action of Kashiwara operators
$\eit$, $\fit$ on $\F(\Lambda)$. Fix $i\in I$ and let
$Y=(y_k)_{k=0}^{\infty}\in \F(\Lambda)$ be a proper Young wall.
\begin{enumerate}
\item To each column $y_k$ of $Y$, we assign

\vskip 2mm
\begin{center}
$\begin{cases}
-- & \text{if $y_k$ is twice $i$-removable,} \\
- & \text{if $y_k$ is once $i$-removable but not $i$-admissible,} \\
-+ & \text{if $y_k$ is once $i$-removable and once
$i$-admissible,} \\
+ & \text{if $y_k$ is once $i$-admissible but not $i$-removable,}
\\
++ & \text{if $y_k$ is twice $i$-admissible,} \\
\ \ \cdot & \text{otherwise}.
\end{cases}$
\end{center}

\item From the (infinite) sequence of $+$'s and $-$'s, cancel out every
$(+,-)$-pair to obtain a finite sequence of $-$'s followed by
$+$'s, reading from left to right. This sequence $(-\cdots -
+\cdots +)$ is called the {\it $i$-signature} of the proper Young
wall $Y$.
\item We define $\eit Y$ to be the proper Young wall obtained from $Y$
by removing the $i$-block corresponding to the right-most $-$ in
the $i$-signature of $Y$. We define $\eit Y = 0$ if there exists
no $-$ in the $i$-signature of $Y$.
\item We define $\fit Y$ to be the proper Young wall obtained from $Y$
by adding an $i$-block to the column corresponding to the
left-most $+$ in the $i$-signature of $Y$. We define $\fit Y = 0$
if there exists no $+$ in the $i$-signature of $Y$.
\end{enumerate}

We also define
\begin{equation*}
\begin{split}
{\rm wt}(Y)&=\Lambda-\sum_{i\in I}k_i\alpha_i \in P, \\
\varepsilon_i(Y)&=\text{the number of $-$'s in the $i$-signature of $Y$}, \\
\varphi_i(Y)&=\text{the number of $+$'s in the $i$-signature of
$Y$},
\end{split}
\end{equation*}
where $k_i$ denotes the number of $i$-blocks in $Y$ that have been
added to $Y_{\Lambda}$.

\vskip 3mm
\begin{prop}\label {Kang}{\rm \cite{HK, Kang}}
The set $\F(\Lambda)$ together with the maps ${\rm wt}$,
$\varepsilon_i$, $\varphi_i$, $\tilde{e}_i$ and $\tilde{f}_i$ {\rm
($i\in I$)} becomes an affine crystal.
\end{prop}

\vskip 1mm Let $\delta=d_0\alpha_0+\cdots+d_n\alpha_n$ be the null
root for the quantum affine algebra $U_q(\widehat{\g})$, and set
$a_i=d_i$ if $\widehat{\g}\neq D_{n+1}^{(2)}$, $a_i=2d_i$ if
$\widehat{\g}=D_{n+1}^{(2)}$. The part of a column consisting of
$a_0$-many 0-blocks, $a_1$-many 1-blocks, $\cdots$, $a_n$-many
$n$-blocks in some cyclic order is called a {\it $\delta$-column}.

\vskip 3mm
\begin{defi} (a) A column in a proper Young wall is said to {\it contain a
removable $\delta$} if we may remove a $\delta$-column from $Y$
and still obtain a proper Young wall.

(b) A proper Young wall is said to be {\it reduced} if none of its
columns contain a removable $\delta$.

\end{defi}

\vskip 1mm Let $\Y(\Lambda)\subset \F(\Lambda)$ denote the set of
all reduced proper Young walls on $Y_{\Lambda}$ with a affine
dominant integral weight $\Lambda$ of level $1$. Then we have:

\vskip 3mm
\begin{prop} \cite{HK, Kang} \label{prop:aff}
Let $B(\lambda)$ be the crystal basis of the basic representation
$V(\Lambda)$ over quantum affine algebras, then there exists a
crystal isomorphism
$$
\Y(\Lambda) \stackrel{\sim} \longrightarrow B(\Lambda) \quad
\text{given by} \ \ Y_{\Lambda} \longmapsto u_{\Lambda},$$ where
$u_{\Lambda}$ is the highest weight vector in $B(\Lambda)$.
\end{prop}

\vskip 3mm Let $\frak g$ be a classical Lie algebra of type $A_n,
C_n,B_n$ and $D_n$. Then these Lie algebras lie inside an affine
Lie algebra $\widehat{\frak g}=A_n^{(1)},A_{2n-1}^{(2)},B_n^{(1)}$
and $D_n^{(1)}$, respectively. (This fact can be expected because
the Dynkin diagram of $\frak g$ can be obtained by removing the
$0$-node from the Dynkin diagram of $\widehat{\frak g}$.)

Fix such a pair $\frak g \subset \hat{\frak g}$ and let $\Lambda$
be a dominant integral weight of level 1 for the affine Lie
algebra $\hat{\frak g}$. Then by Proposition \ref{prop:aff}, the
crystal graph $B(\Lambda)$ is realized as the set $Y(\Lambda)$ of
all reduced proper Young walls built on the ground-state wall
$Y_{\Lambda}$. If we remove all $0$-arrows in $Y(\Lambda)$, then
it is decomposed into a disjoint union of infinitely many
connected components. Moreover, we can show that each connected
component is isomorphic to the crystal graph $B(\lambda)$ for some
dominant integral weight $\lambda$ for $\frak g$.

Conversely, any crystal graph $B(\lambda)$ for $\frak g$ arises in
this way. That is, given a dominant integral weight $\lambda$ for
$\frak g$, there is a dominant integral weight $\Lambda$ of level
$1$ for $\widehat{\frak g}$ such that $B(\lambda)$ appears as a
connected component in $B(\Lambda)$ without $0$-arrows.

\vskip 2mm The first step is to identify the highest weight vector
$H_{\lambda}$ for $B(\lambda)$ with some reduced proper Young wall
in $\Y(\Lambda)$ which is annihilated by all $\tilde{e_i}$ for
$i=1, \cdots, n$. However, given a dominant integral weight $\la$
for $\frak g$, there are infinitely many such Young walls in
$\Y(\Lambda)$. Equivalently, given $\lambda$, there are infinitely
many connected components of $\Y(\Lambda)$ without $0$-arrows that
are isomorphic to $B(\lambda)$. Among these, they choose the
characterization of $B(\lambda)$ corresponding to the connected
components having the least number of blocks \cite{KKLS}.

\vskip 2mm  From now on, we focus on the classical Lie algebra
$\frak g=B_n$. We define the linear functionals $\omega_i$ by
\begin{center}
$
\omega_i = \begin{cases}
\la_i \quad & \text{for} \ \ i=1, \cdots, n-1, \\
2 \la_n \quad & \text{for} \ \ i=n.
\end{cases}
$
\end{center}

Let $F(\la)\subset Y(\La)$ be the set of all reduced proper Young
walls lying between highest weight vector $H_{\la}$ and lowest
weight vector $L_{\la}$. (An algorithm of constructing the highest
weight vector $H_{\la}$ and lowest weight vector $L_{\la}$ inside
$\Y(\Lambda)$ was given in \cite{KKLS}.) Set
$\la=\omega_{i_1}+\cdots+\omega_{i_t}+b\lambda_n$ for $b=0$ or
$1$. For each $Y \in F(\la)$, we denote by $\overset {\circ}
{Y}_{\omega_{i_k}}$ ($k=1,\cdots,t$) (resp. $\overset
{\circ}{Y}_{\la_n}$) the part of $Y$ consisting of the blocks
lying above ${\overline H}_{\omega_{i_k}}$ (resp. $H_{\la_n}$) and
we denote by ${\overline Y_{\omega_{i_k}}}$ (resp. ${\overline
Y_{\la_n}}$)  the intersection of $Y$ and ${\overline
L}_{\omega_{i_k}}$ (resp. $L_{\la_n}$), where ${\overline
H}_{\omega_{i_k}}$ (resp. ${\overline L}_{\omega_{i_k}}$) is the
Young wall consisting of $H_{\omega_{i_k}}$ (resp.
$L_{\omega_{i_k}}$), and $i_k\times (t-k)$ or $i_k\times
(t-k+\frac{1}{2})$-many $\delta$-columns. Moreover, we denote by
${\overline Y_{\omega_{i_k}+\omega_{i_{k+1}}}}$ (resp. ${\overline
Y_{\omega_{i_t}+\la_n}}$) the union of ${\overline
Y_{\omega_{i_k}}}$ (resp. ${\overline Y_{\omega_{i_t}}}$) and
${\overline Y_{\omega_{i_{k+1}}}}$ (resp. ${\overline
Y_{\la_n}}$).

\savebox{\tmpfiga}{\begin{texdraw} \fontsize{7}{7}\selectfont
\drawdim em \setunitscale 0.85

\nc{\dtri}{ \bsegment \lvec(-2 0) \lvec(-2 2)\lvec(0 2)\lvec(0
0)\ifill f:0.7 \esegment }

\nc{\dtrii}{ \bsegment \lvec(-2 0) \lvec(-2 1)\lvec(0 1)\lvec(0
0)\ifill f:0.7 \esegment }

\nc{\dtriii}{ \bsegment \lvec(-2 0) \lvec(0 2)\lvec(0 0)\ifill
f:0.7 \esegment }

\nc{\dtriv}{ \bsegment \lvec(-2 0) \lvec(-2 -2)\lvec(0 0)\ifill
f:0.7 \esegment }

\nc{\dtrv}{ \bsegment \lvec(-2 0) \lvec(-2 2)\lvec(0 2)\lvec(0
0)\ifill f:0.85 \esegment }

\nc{\dtrvi}{ \bsegment \lvec(-2 0) \lvec(-2 1)\lvec(0 1)\lvec(0
0)\ifill f:0.85 \esegment }

\move(52 24)\dtrvi \move(48 22)\dtrv \move(50 22)\dtrv \move(48
20)\dtrv \move(50 24)\dtrvi\move(48 24)\dtrvi

\move(46 13)\dtrvi\move(46 14)\dtrv\move(46 16)\dtrv \move(44
13)\dtrvi\move(44 14)\dtrv\move(42 13)\dtrvi

\move(46 12)\dtrii\move(44 12)\dtrii\move(42 12)\dtrii\move(40
12)\dtrii \move(44 10)\dtri\move(42 10)\dtri\move(40 10)\dtri
\move(42 8)\dtri\move(40 8)\dtri\move(40 6)\dtriii

\move(38 1)\dtrii\move(36 1)\dtrii\move(34 1)\dtrii\move(32
1)\dtrii \move(38 2)\dtri\move(36 2)\dtri\move(34 2)\dtri \move(38
4)\dtri\move(36 4)\dtri\move(38 8)\dtriv

\setgray 0.6

\move(52 22)\rlvec(0 5) \move(52 26)\rlvec(-2 0)\rlvec(0 -26)
\move(52 25)\rlvec(-2 0)\move(52 24)\rlvec(-4 0)\rlvec(0
-24)\move(52 22)\rlvec(-4 0)\move(52 20)\rlvec(-6 0)\rlvec(0
-20)\move(52 18)\rlvec(-8 0)\rlvec(0 -18)\move(52 16)\rlvec(-8
0)\move(52 14)\rlvec(-10 0)\rlvec(0 -14)\move(52 13)\rlvec(-10 0)
\move(52 12)\rlvec(-12 0)\rlvec(0 -12)\move(52 10)\rlvec(-14
0)\rlvec(0 -10)\move(52 8)\rlvec(-16 0)\rlvec(0 -8)\move(52
6)\rlvec(-18 0)\rlvec(0 -6)\move(52 4)\rlvec(-20 0)\rlvec(0 -4)
\move(52 2)\rlvec(-22 0)\rlvec(0 -2)\move(52 1)\rlvec(-22 0)

\move(52 20)\rlvec(-2 -2) \move(50 20)\rlvec(-2 -2)\move(48
20)\rlvec(-2 -2) \move(52 8)\rlvec(-2 -2) \move(50 8)\rlvec(-2
-2)\move(48 8)\rlvec(-2 -2) \move(46 8)\rlvec(-2 -2)\move(44
8)\rlvec(-2 -2) \move(42 8)\rlvec(-2 -2)\move(40 8)\rlvec(-2
-2)\move(38 8)\rlvec(-2 -2)

\move(38 -2)\setgray 0.4 \rlvec(0 20)\move(46 -2)\rlvec(0 28)

\move(52 22)\linewd 0.15 \setgray 0 \rlvec(-2 0)\rlvec(0 -2)
\rlvec(-4 0) \rlvec(0 -8)\rlvec(-2 0)\rlvec(0 -2)\rlvec(-2
0)\rlvec(0 -2)\rlvec(-4 0)\rlvec(0 -7)\rlvec(-8 0)\rlvec(0
-1)\rlvec(22 0)\rlvec(0 22)\move(48 20)\rlvec(-2 -2)\move(40
8)\rlvec(-2 -2)

\htext(51 25.5){$4$} \htext(51 24.5){$4$}\htext(51
23){$3$}\htext(49 23){$3$}\htext(51 21){$2$}\htext(49 21){$2$}
\htext(51.5 18.5){$1$}\htext(50.5 19.5){$0$}\htext(49.5
18.5){$0$}\htext(48.5 19.5){$1$}\htext(47.5 18.5){$1$}\htext(46.5
19.5){$0$}\htext(51 17){$2$}\htext(49 17){$2$}\htext(47
17){$2$}\htext(45 17){$2$}\htext(51 15){$3$}\htext(49
15){$3$}\htext(47 15){$3$}\htext(45 15){$3$}

\htext(51 13.5){$4$}\htext(49 13.5){$4$}\htext(47
13.5){$4$}\htext(45 13.5){$4$}\htext(43 13.5){$4$} \htext(51
12.5){$4$}\htext(49 12.5){$4$}\htext(47 12.5){$4$}\htext(45
12.5){$4$}\htext(43 12.5){$4$} \htext(51 11){$3$}\htext(49
11){$3$}\htext(47 11){$3$}\htext(45 11){$3$}\htext(43
11){$3$}\htext(41 11){$3$} \htext(51 9){$2$}\htext(49
9){$2$}\htext(47 9){$2$}\htext(45 9){$2$}\htext(43
9){$2$}\htext(41 9){$2$}\htext(39 9){$2$}

\htext(51.5 6.5){$1$} \htext(50.5 7.5){$0$} \htext(49.5 6.5){$0$}
\htext(48.5 7.5){$1$}\htext(47.5 6.5){$1$}  \htext(46.5 7.5){$0$}
\htext(45.5 6.5){$0$} \htext(44.5 7.5){$1$}\htext(43.5 6.5){$1$}
\htext(42.5 7.5){$0$}\htext(41.5 6.5){$0$} \htext(40.5
7.5){$1$}\htext(39.5 6.5){$1$}  \htext(38.5 7.5){$0$}\htext(36.5
7.5){$1$}

\htext(51 5){$2$}\htext(49 5){$2$}\htext(47 5){$2$}\htext(45
5){$2$}\htext(43 5){$2$} \htext(41 5){$2$} \htext(39 5){$2$}
\htext(37 5){$2$} \htext(35 5){$2$} \htext(51 3){$3$} \htext(49
3){$3$} \htext(47 3){$3$} \htext(45 3){$3$}\htext(43 3){$3$}
\htext(41 3){$3$} \htext(39 3){$3$} \htext(37 3){$3$} \htext(35
3){$3$}\htext(33 3){$3$} \htext(51 1.5){$4$} \htext(49 1.5){$4$}
\htext(47 1.5){$4$} \htext(45 1.5){$4$}\htext(43 1.5){$4$}
\htext(41 1.5){$4$} \htext(39 1.5){$4$} \htext(37 1.5){$4$}
\htext(35 1.5){$4$}\htext(33 1.5){$4$}\htext(31 1.5){$4$}
\htext(51 0.5){$4$} \htext(49 0.5){$4$} \htext(47 0.5){$4$}
\htext(45 0.5){$4$}\htext(43 0.5){$4$} \htext(41 0.5){$4$}
\htext(39 0.5){$4$} \htext(37 0.5){$4$} \htext(35
0.5){$4$}\htext(33 0.5){$4$}\htext(31 0.5){$4$}

\htext(34 -1){$\overline{Y}_{\la_4}$}\htext(42
-1){$\overline{Y}_{\omega_4}$}\htext(49
-1){$\overline{Y}_{\omega_3}$}
\end{texdraw}}
\savebox{\tmpfigg}{\begin{texdraw} \fontsize{7}{7}\selectfont
\drawdim em \setunitscale 0.85

\rlvec(2 0)\rlvec(0 6)\rlvec(4 0)\move(0 0)\rlvec(0 2)\rlvec(4
0)\rlvec(0 6)\rlvec(2 0)\rlvec(0 -4)\rlvec(-4 0)\move(0 0)\rlvec(2
2) \move(4 7)\rlvec(2 0)

\htext(1 0.5){$1$}\htext(3 3){$2$}\htext(3 5){$3$} \htext(5
5){$3$}\htext(5 6.5){$4$}\htext(5 7.5){$4$}
\end{texdraw}}
\savebox{\tmpfigh}{\begin{texdraw} \fontsize{7}{7}\selectfont
\drawdim em \setunitscale 0.85

\rlvec(2 0)\rlvec(0 6)\rlvec(6 0)\rlvec(0 6)\rlvec(-2 0)\rlvec(0
-8)\rlvec(-6 0)\rlvec(0 -4)\rlvec(2 2)\move(0 2)\rlvec(4
0)\rlvec(0 6)\rlvec(4 0)

\move(4 7)\rlvec(4 0)\move(6 10)\rlvec(2 0)

\htext(1 0.5){$1$}\htext(1 3){$2$}\htext(3 3){$2$}\htext(3 5){$3$}
\htext(5 5){$3$}\htext(5 6.5){$4$}\htext(5 7.5){$4$}\htext(7
6.5){$4$}\htext(7 7.5){$4$} \htext(7 9){$3$}\htext(7 11){$2$}
\end{texdraw}}
\savebox{\tmpfigi}{\begin{texdraw} \fontsize{7}{7}\selectfont
\drawdim em \setunitscale 0.85

\rlvec(8 0)\rlvec(0 7)\rlvec(-2 0)\rlvec(0 -7) \move(8 5)\rlvec(-4
0)\rlvec(0 -5) \move(8 3)\rlvec(-6 0)\rlvec(0 -3) \move(8
1)\rlvec(-8 0)\rlvec(0 -1)\move(6 5)\rlvec(2 2)

\htext(7 0.5){$4$}\htext(7 2){$3$}\htext(7 4){$2$} \htext(7
6.5){$1$}\htext(5 0.5){$4$}\htext(5 2){$3$}\htext(5 4){$2$}
\htext(3 0.5){$4$}\htext(3 2){$3$}\htext(1 0.5){$4$}
\end{texdraw}}

\begin{example} \label{ex-kkls}
If $\frak g=B_4$, $\la=\omega_3+\omega_4+\la_4$, and let
$$Y=\raisebox{-0.4\height}{\usebox{\tmpfiga}}\,,$$ we have
$$
\overset
{\circ}{Y}_{\omega_3}=\raisebox{-0.4\height}{\usebox{\tmpfigg}},
\overset
{\circ}{Y}_{\omega_4}=\raisebox{-0.4\height}{\usebox{\tmpfigh}}\,\,\,\,
\text{and}\,\,\,\overset
{\circ}{Y}_{\la_4}=\raisebox{-0.4\height}{\usebox{\tmpfigi}}\,.
$$
\end{example}

\vskip 3mm
\begin{defi}
(a) For ${\overline Y_{\omega_{i_k}}}$, we define
$Y_{\omega_{i_k}}^+$ (resp. $Y_{\omega_{i_k}}^-$) for
$k=1,\cdots,t$ by the part consisting of the blocks lying above
(resp. below) the $n$-th row.

(b) For ${\overline Y_{\omega_{i_k}+\omega_{i_{k+1}}}}$ and
${\overline Y_{\omega_{i_t}+\la_n}}$ of $Y$, we define
$Y^{\omega_{i_k}}$ (resp. $Y^{\omega_{i_{k+1}}}$) in ${\overline
Y_{\omega_{i_k}+\omega_{i_{k+1}}}}$ is the intersection between
the part consisting of blocks lying above the highest weight
vector ${\overline H}_{\omega_{i_k}}$ (resp. ${\overline
H}_{\omega_{i_{k+1}}}$) and the part consisting of blocks lying
below the lowest weight vector ${\overline L}_{\omega_{i_{k+1}}}$
(resp. ${\overline L}_{\omega_{i_k}}$) reading from top (resp.
right) to bottom (resp. left). Moreover, $Y^{\omega_{i_t}}$ and
$Y^{\lambda_n}$ are defined by a similar way.

(c) For ${\overline Y_{\omega_{i_k}+\omega_{i_{k+1}}}}$ of $Y$, we
define $L_{(\omega_{i_k},\omega_{i_{k+1}})}^-$ (resp.
$L_{(\omega_{i_k},\omega_{i_{k+1}})}^+$) is the part of
$L_{\omega_{i_k}}$ (resp. $L_{\omega_{i_{k+1}}}$) consisting of
the right isosceles triangular blocks below (resp. above) the
$n$-row. Moreover, $L_{(\omega_{i_t},\la_n)}^-$ and
$L_{(\omega_{i_t},\la_n)}^+$ are defined by a similar way. Then
$$\aligned
&Y_{(\omega_{i_k},\omega_{i_{k+1}})}^-=Y\cap
L_{(\omega_{i_k},\omega_{i_{k+1}})}^-\,\,\,\text{and}\,\,\,
Y_{(\omega_{i_k},\omega_{i_{k+1}})}^+=Y\cap
L_{(\omega_{i_k},\omega_{i_{k+1}})}^+,\\
&Y_{(\omega_{i_t},\la_n)}^-=Y\cap
L_{(\omega_{i_k},\la_n)}^-\,\,\,\text{and}\,\,\,
Y_{(\omega_{i_k},\la_n)}^+=Y\cap L_{(\omega_{i_k},\la_n)}^+.
\endaligned
$$

(d) We define $|Y_{\omega_{i_k}}^{-}|$ by  the wall obtained by
reflecting $Y_{\omega_{i_k}}^{-}$ along the $n$-row and shifting
the blocks to the right as much as possible. Moreover, we define
$|Y_{(\omega_{i_k},\omega_{i_{k+1}})}^{-}|$ (resp.
$|Y_{(\omega_{i_t},\la_n)}^{-}|$) by the wall obtained by
reflecting $Y_{(\omega_{i_k},\omega_{i_{k+1}})}^{-}$ (resp.
$Y_{(\omega_{i_t},\la_n)}^{-}$) with respect to the upper $n$-row
and shifting the blocks to the right as much as possible.
\end{defi}

\savebox{\tmpfigb}{\begin{texdraw} \fontsize{7}{7}\selectfont
\drawdim em \setunitscale 0.85

\rlvec(2 0)\rlvec(0 6)\rlvec(4 0)\rlvec(0 -2)\rlvec(-4 0)\move(0
0)\rlvec(0 2)\rlvec(4 0)\rlvec(0 4)\move(0 0)\rlvec(2 2)

\htext(1.5 0.5){$1$} \htext(3 3){$2$}\htext(3 5){$3$}\htext(5
5){$3$}
\end{texdraw}}
\savebox{\tmpfigc}{\begin{texdraw} \fontsize{7}{7}\selectfont
\drawdim em \setunitscale 0.85

\rlvec(2 0)\rlvec(0 4)\rlvec(-2 0)\rlvec(0 -4)\move(0 2)\rlvec(2
0)

\htext(1 1){$3$} \htext(1 3){$2$}
\end{texdraw}}
\savebox{\tmpfigd}{\begin{texdraw} \fontsize{7}{7}\selectfont
\drawdim em \setunitscale 0.85

\rlvec(2 0)\rlvec(0 6)\rlvec(4 0)\rlvec(0 -2)\rlvec(-4 0)\move(0
0)\rlvec(0 2)\rlvec(4 0)\rlvec(0 4)\move(0 0)\rlvec(2 2)\move(2
4)\rlvec(-2 0)\rlvec(0 -2)

\htext(1.5 0.5){$1$}\htext(1 3){$2$} \htext(3 3){$2$}\htext(3
5){$3$}\htext(5 5){$3$}
\end{texdraw}}
\savebox{\tmpfige}{\begin{texdraw} \fontsize{7}{7}\selectfont
\drawdim em \setunitscale 0.85

\rlvec(2 0)\rlvec(0 4)\rlvec(2 0)\rlvec(0 -2)\rlvec(-4 0)\rlvec(0
-2)\move(2 3)\rlvec(2 0)

\htext(1 1){$3$}\htext(3 2.5){$4$} \htext(3 3.5){$4$}
\end{texdraw}}
\savebox{\tmpfigf}{\begin{texdraw} \fontsize{7}{7}\selectfont
\drawdim em \setunitscale 0.85

\rlvec(2 0)\rlvec(0 6)\rlvec(4 0)\move(0 0)\rlvec(0 4)\rlvec(6
0)\rlvec(0 6)\rlvec(-2 0)\rlvec(0 -8)\rlvec(-4 0)\move(2
5)\rlvec(4 0)\move(4 8)\rlvec(2 0)

\htext(1 1){$2$}\htext(1 3){$3$}\htext(3 3){$3$} \htext(3
4.5){$4$}\htext(3 5.5){$4$}\htext(5 4.5){$4$}\htext(5
5.5){$4$}\htext(5 7){$3$}\htext(5 9){$2$}
\end{texdraw}}
\savebox{\tmpfigg}{\begin{texdraw} \fontsize{7}{7}\selectfont
\drawdim em \setunitscale 0.85

\rlvec(4 0)\rlvec(0 5)\rlvec(-2 0)\rlvec(0 -5)\move(0 0)\rlvec(0
1)\rlvec(4 0)\move(2 3)\rlvec(2 0)

\htext(1 0.5){$4$}\htext(3 0.5){$4$}\htext(3 2){$3$} \htext(3
4){$2$}
\end{texdraw}}
\savebox{\tmpfigh}{\begin{texdraw} \fontsize{7}{7}\selectfont
\drawdim em \setunitscale 0.85

\rlvec(6 0)\rlvec(0 7)\rlvec(-2 0)\rlvec(0 -7)\move(0 0)\rlvec(0
1)\rlvec(6 0)\move(2 0)\rlvec(0 5)\rlvec(4 0)\move(2 3)\rlvec(4
0)\move(4 5)\rlvec(2 2)

\htext(1 0.5){$4$}\htext(3 0.5){$4$}\htext(3 2){$3$} \htext(3
4){$2$}\htext(5 0.5){$4$}\htext(5 2){$3$} \htext(5
4){$2$}\htext(4.5 6.5){$1$}
\end{texdraw}}
\savebox{\tmpfigi}{\begin{texdraw} \fontsize{7}{7}\selectfont
\drawdim em \setunitscale 0.85

\rlvec(2 0)\rlvec(0 3)\rlvec(2 0)\rlvec(0 -1)\rlvec(-4 0)\rlvec(0
-2)

 \htext(1 1){$3$}\htext(3 2.5){$4$}
\end{texdraw}}
\savebox{\tmpfigj}{\begin{texdraw} \fontsize{7}{7}\selectfont
\drawdim em \setunitscale 0.85

\rlvec(4 0)\rlvec(0 5)\rlvec(-2 0)\rlvec(0 -5)\move(0 0)\rlvec(0
1)\rlvec(4 0)\move(2 3)\rlvec(2 0)

\htext(1 0.5){$4$}\htext(3 0.5){$4$}\htext(3 2){$3$} \htext(3
4){$2$}
\end{texdraw}}
\savebox{\tmpfigk}{\begin{texdraw} \fontsize{7}{7}\selectfont
\drawdim em \setunitscale 0.85

\rlvec(2 0)\rlvec(0 6)\rlvec(6 0)\rlvec(0 1)\rlvec(-4 0)\rlvec(0
-5)\rlvec(-4 0) \move(0 0)\rlvec(0 4)\rlvec(6 0)\rlvec(0 3)\move(0
0)\rlvec(2 2)

\htext(1.5 0.5){$1$}\htext(1 3){$2$}\htext(3 3){$2$} \htext(3
5){$3$}\htext(5 5){$3$}\htext(5 6.5){$4$}\htext(7 6.5){$4$}
\end{texdraw}}
\savebox{\tmpfigl}{\begin{texdraw} \fontsize{7}{7}\selectfont
\drawdim em \setunitscale 0.85

\rlvec(8 0)\rlvec(0 7)\rlvec(-2 0)\rlvec(0 -7) \move(8 5)\rlvec(-4
0)\rlvec(0 -5) \move(8 3)\rlvec(-6 0)\rlvec(0 -3) \move(8
1)\rlvec(-8 0)\rlvec(0 -1)\move(6 5)\rlvec(2 2)

\htext(7 0.5){$4$}\htext(7 2){$3$}\htext(7 4){$2$} \htext(7
6.5){$1$}\htext(5 0.5){$4$}\htext(5 2){$3$}\htext(5 4){$2$}
\htext(3 0.5){$4$}\htext(3 2){$3$}\htext(1 0.5){$4$}
\end{texdraw}}

\vskip 2mm
\begin{example}
Let $Y$ be a Young wall given in Example \ref{ex-kkls}, then we
have
$$Y_{\omega_3}^+=\varnothing,
Y_{\omega_3}^-=\raisebox{-0.4\height}{\usebox{\tmpfigb}},
Y_{\omega_4}^+=\raisebox{-0.4\height}{\usebox{\tmpfigc}}\,\,\,\,
\text{and}\,\,\,Y_{\omega_4}^-=\raisebox{-0.4\height}{\usebox{\tmpfigd}}\,$$
$$Y^{\omega_3}=\raisebox{-0.4\height}{\usebox{\tmpfige}},\,\,\,\,
Y^{\omega_4}=\raisebox{-0.4\height}{\usebox{\tmpfigf}}\,\,\,\,\text{in}\,\,\,
{\overline Y_{\omega_3+\omega_4}}\,,$$
$$Y^{\omega_4}=\raisebox{-0.4\height}{\usebox{\tmpfigg}}\,,\,\,\,\,
Y^{\la_4}=\raisebox{-0.4\height}{\usebox{\tmpfigl}}\,\,\,\,\text{in}\,\,\,
{\overline Y_{\omega_4+\la_4}}\,.$$ Moreover, the shaded parts
represent $L_{(\omega_3,\omega_4)}^-$,
$L_{(\omega_3,\omega_4)}^+$, $Y_{(\omega_4,\la_4)}^-$ and
$Y_{(\omega_4,\la_4)}^+$, and
$$Y_{(\omega_3,\omega_4)}^-=\raisebox{-0.4\height}{\usebox{\tmpfigi}},\,\,\,\,
Y_{(\omega_3,\omega_4)}^+=\raisebox{-0.4\height}{\usebox{\tmpfigj}}\,\,\,\,
Y_{(\omega_4,\la_4)}^-=\raisebox{-0.4\height}{\usebox{\tmpfigk}},\,\,\,\,
Y_{(\omega_4,\la_4)}^+=\raisebox{-0.4\height}{\usebox{\tmpfigl}}\,.$$
\end{example}

\savebox{\tmpfiga}{\begin{texdraw} \fontsize{7}{7}\selectfont
\drawdim em \setunitscale 0.85

\move(0 -0.5)\lvec(0 3)\rlvec(2 0)\rlvec(0 -3.5)\move(0
1.2)\rlvec(2 0)

\htext(1 0.5){$a\!\!-\!\!2$}\htext(1 2.1){$a\!\!-\!\!1$}
\end{texdraw}}
\savebox{\tmpfigb}{\begin{texdraw} \fontsize{7}{7}\selectfont
\drawdim em \setunitscale 0.85

\move(0 -0.5)\lvec(0 3)\rlvec(2 0)\rlvec(0 -3.5)\move(0
1.2)\rlvec(2 0)

\htext(1 0.5){$a\!\!+\!\!1$}\htext(1 2.1){$a$}
\end{texdraw}}
\savebox{\tmpfigc}{\begin{texdraw} \fontsize{7}{7}\selectfont
\drawdim em \setunitscale 0.85

\move(0 -0.5)\lvec(0 3)\rlvec(2 0)\rlvec(0 -3.5)\move(0
1.2)\rlvec(2 0)\move(0 1.2)\rlvec(2 1.8)

\htext(1 0.5){$2$}\htext(1.5 1.87){$1$}\htext(0.5 2.5){$0$}
\end{texdraw}}
\savebox{\tmpfigd}{\begin{texdraw} \fontsize{7}{7}\selectfont
\drawdim em \setunitscale 0.85

\move(0 -0.5)\lvec(0 3)\rlvec(2 0)\rlvec(0 -3.5)\move(0
1.2)\rlvec(2 0)

\htext(1 0.5){$3$}\htext(1 2.1){$2$}
\end{texdraw}}
\savebox{\tmpfigf}{\begin{texdraw} \fontsize{7}{7}\selectfont
\drawdim em \setunitscale 0.85

\rlvec(2 0)\rlvec(0 2)\rlvec(-2 0)\rlvec(0 -2)

\htext(1 1){$3$}
\end{texdraw}}
\savebox{\tmpfigg}{\begin{texdraw} \fontsize{7}{7}\selectfont
\drawdim em \setunitscale 0.85

\rlvec(2 0)\rlvec(0 4)\rlvec(-2 0)\rlvec(0 -4)\move(0 2)\rlvec(2
0)

\htext(1 1){$3$}\htext(1 3){$2$}
\end{texdraw}}
\savebox{\tmpfigh}{\begin{texdraw} \fontsize{7}{7}\selectfont
\drawdim em \setunitscale 0.85

\rlvec(2 0)\rlvec(0 4)\move(0 0)\rlvec(0 4)\rlvec(4 0)\rlvec(0
-2)\rlvec(-4 0)

\htext(1 1){$2$}\htext(1 3){$3$}\htext(3 3){$3$}
\end{texdraw}}
\savebox{\tmpfigi}{\begin{texdraw} \fontsize{7}{7}\selectfont
\drawdim em \setunitscale 0.85

\rlvec(4 0)\rlvec(0 4)\rlvec(-2 0)\rlvec(0 -4)\move(0 0)\rlvec(0
2) \rlvec(4 0)

\htext(1 1){$3$}\htext(3 1){$3$}\htext(3 3){$2$}
\end{texdraw}}

\vskip 3mm  Now, we assume that ${\overline
Y_{\omega_{i_k}+\omega_{i_{k+1}}}}$ and ${\overline
Y_{\omega_{i_t}+\la_{n}}}$ satisfy the following condition: If the
top of the $p$-th column of ${\overline Y_{\omega_{i_k}}}$ (resp.
${\overline Y_{\omega_{i_t}}}$) from the right is
$\raisebox{-0.4\height}{\usebox{\tmpfiga}}$ and the top of the
$q$-th column of ${\overline Y_{\omega_{i_{k+1}}}}$ (resp.
${\overline Y_{\la_n}}$) from the right is
$\raisebox{-0.4\height}{\usebox{\tmpfigb}}$ with $p>q$, then it is
called  ${\overline Y_{\omega_{i_k}+\omega_{i_{k+1}}}}$ (resp.
${\overline Y_{\omega_{i_t}+\la_{n}}}$) satisfies {\bf (C1)}.

\begin{defi}
(a) We define $L_{\omega_{i_k}}^{+}(a;p,q)$ (resp.
$L_{\omega_{i_{k+1}}}^{+}(a;p,q)$) to be the right isosceles
triangle formed by $a$-block in the $q$-th column,
$(a+p-q-1)$-block in the $q$-th column and $(a+p-q-1)$-block in
the $(p-1)$-th column in $Y_{\omega_{i_k}}$ (resp.
$Y_{\omega_{i_{k+1}}}$).

(b) We define $L_{\omega_{i_k}}^{-}(a;p,q)$ (resp.
$L_{\omega_{i_{k+1}}}^{-}(a;p,q)$) by the wall obtained by
reflecting $L_{\omega_{i_k}}^{+}(a;p,q)$ (resp.
$L_{\omega_{i_{k+1}}}^{+}(a;p,q)$) with respect to the $n$-row.

(c) We also define $Y_{\omega_{i_k}}^{\pm}(a;p,q)$,
$Y_{\omega_{i_k}}^{\pm}(a;p,q)$ and
$Y_{\omega_{i_k}}^{\pm}(a;p,q)$ by
$$
\begin{aligned}
& Y_{\omega_{i_k}}^{\pm}(a;p,q)  = L_{\omega_{i_k}}^{\pm}(a;p,q)
\cap Y, \quad Y_{\omega_{i_{k+1}}}^{\pm}(a;p,q) =
L_{\omega_{i_{k+1}}}^{\pm}(a;p,q) \cap Y,\\
&Y_{\omega_{i_t}}^{\pm}(a;p,q) = L_{\omega_{i_t}}^{\pm}(a;p,q)
\cap Y,
\end{aligned}
$$

(d) $|Y_{\omega_{i_k}}^{-}(a;p,q)|$ is defined by the wall
obtained by reflecting $Y_{\omega_{i_k}}^{-}(a;p,q)$ with respect
to the $n$-row and shifting the blocks to the right as much as
possible
\end{defi}

\begin{example}
Let $Y$ be a Young wall given in Example \ref{ex-kkls}, then we
can see that there exist
\begin{center}
$\raisebox{-0.4\height}{\usebox{\tmpfigc}}$ \quad and \quad
$\raisebox{-0.4\height}{\usebox{\tmpfigd}}$ \end{center} on top of
the third and first column in ${\overline Y_{\omega_3}}$ and
${\overline Y_{\omega_4}}$, respectively. In this case, $a=2$ and
$$Y_{\omega_3}^+(2;3,1)=\varnothing,\,\,\,
Y_{\omega_3}^-(2;3,1)=\raisebox{-0.4\height}{\usebox{\tmpfigf}}\,,\,\,\,
Y_{\omega_4}^+(2;3,1)=\raisebox{-0.4\height}{\usebox{\tmpfigg}}\,\,\,\text{and}\,\,\,
Y_{\omega_4}^-(2;3,1)=\raisebox{-0.4\height}{\usebox{\tmpfigh}}\,.$$
\end{example}

\vskip 3mm Now, we close this section with the realization theorem
using Young walls for crystal bases.

\vskip 3mm
\begin{prop} \cite{KKLS}
Let $\la \in P^{+}$ be a dominant integral weight for $\frak g =
B_n$, and write
$$\aligned
&\la = \omega_{i_1} + \cdots + \omega_{i_t} \qquad (1 \le i_1 \le
\cdots \le i_t \le n)\quad \text{or}\\
&\la = \omega_{i_1} + \cdots + \omega_{i_t} + \la_n \qquad (1 \le
i_1 \le \cdots \le i_t \le n). \endaligned$$ We define $Y(\la)$ to
be the set of all reduced proper Young walls satisfying the
following conditions:

{\bf (Y1)} For each $k = 1, \cdots, t$, we have
$Y_{\omega_{i_k}}^{+} \subset |Y_{\omega_{i_k}}^{-}|$;

{\bf (Y2)} For each $k=1, \cdots, t-1$, we have $Y^{\omega_{i_k}}
\subset Y^{\omega_{i_{k+1}}}$, \ \ $Y^{\omega_{i_t}} \subset
Y^{\la_n}$;

{\bf (Y3)} For each $k=1, \cdots, t-1$, we have
$$|Y_{(\omega_{i_k},\omega_{i_{k+1}})}^{-}| \subset
Y_{(\omega_{i_k},\omega_{i_{k+1}})}^{+}, \quad
|Y_{(\omega_{i_t},\la_n)}^{-}| \subset
Y_{(\omega_{i_t},\la_n)}^{+};$$

{\bf (Y4)} For each $k=1, \cdots, t-1$, suppose that ${\overline
Y_{\omega_{i_k}+\omega_{i_{k+1}}}}$ or ${\overline
Y_{\omega_{i_t}+\la_n}}$ satisfy $(C1)$, then we have
$$\aligned
&Y_{\omega_{i_k}}^{+}(p,q, a) \subset
|Y_{\omega_{i_k}}^{-}(p,q,a)|, \quad Y_{\omega_{i_{k+1}}}^{+}(p,q,
a) \subset
|Y_{\omega_{i_{k+1}}}^{-}(p,q,a)|,\\
&Y_{\omega_{i_t}}^{+}(p,q,a) \subset
|Y_{\omega_{i_t}}^{-}(p,q,a)|.
\endaligned$$

Then there is an isomorphism of crystal graphs for
$U_q(B_n)$-modules
$$
Y(\lambda) \stackrel{\sim} \longrightarrow B(\la) \quad
\text{given by} \ \ H_{\la} \longmapsto u_{\la},
$$
where $u_\la$ is the highest weight vector in $B(\la)$.
\end{prop}

\vskip 1cm
\section
{A new realization of crystal bases}
In this section, we give a new realization of crystal bases for
finite dimensional irreducible modules over classical Lie algebras
derived from Young wall realization $Y(\la)$. The basis vectors
are parameterized by certain tableaux with given shapes, which is
different from generalized Young diagram given by Kashiwara and
Nakashima.

\vskip 2mm We define the linear functionals $\omega_i$ by

\vskip 2mm \hskip 3mm 1) $\frak g = A_n, \ C_n$:
\begin{equation*}
\omega_i = \la_i \quad \text{for} \ \ i=1, \cdots, n,
\end{equation*}

\hskip 3mm 2) $\frak g = D_n$:
\begin{equation*}
\omega_i = \begin{cases}
\la_i \quad & \text{for} \ \ i=1, \cdots, n-2, \\
\la_{n-1} + \la_n \quad & \text{for} \ \ i=n-1, \\
2 \la_n \quad & \text{for} \ \ i=n, \\
2 \la_{n-1} \quad & \text{for} \ \ i=n+1.
\end{cases}
\end{equation*}
Then a dominant integral weight $\la$ can be expressed as
\begin{equation}
\la=\left\{\begin{array}{ll} \omega_{i_1}+\cdots+\omega_{i_t}
&\text{if\,\,\, $\frak g=A_n$, $C_n$,}\\
\omega_{i_1}+\cdots+\omega_{i_t}+b\lambda_n
&\text{if\,\,\, $\frak g=B_n$,}\\
\omega_{i_1}+\cdots+\omega_{i_t}+b_1\lambda_{n-1}+b_2\lambda_n
&\text{if\,\,\, $\frak g=D_n$},
\end{array}\right.
\end{equation}
where $1\le i_1\le\cdots\le i_t\le n$, $b=0$ or $1$,
$(b_1,b_2)=(1,0)$ or $(0,1)$.

\vskip 2mm
\begin{defi}\, For a sequence of half integers $l_j\in \frac{1}{2}{\bf
Z}$ ($j=1,\cdots,n$), such that $l_j-l_{j+1}\in {\bf Z}_{\ge 0}$,
let $Y=(l_1,\cdots,l_n)$ be a diagram which has $n$ rows with the
$j$-th row (from bottom to top) of length $|l_j|$. Then $Y$ is
called {\it generalized reverse Young diagram} of type $A_n$ and
$C_n$ (resp. $B_n$) if all $l_j$ are non-negative integers (resp.
half integers). Moreover, $Y$ is called {\it generalized reverse
Young diagram} of type $D_n$ if all $l_j$ are half integers and
$l_1\ge l_2\ge\cdots\ge l_{n-1}\ge |l_n|$.
\end{defi}

\begin{remark} (a) It is just a diagram obtained by reflecting
generalized Young diagram to the origin.

(b) There is a bijection $GRY$ between the set of generalized
reverse Young diagrams and the set of dominant integral weights of
$\frak g$ defined by
$$GRY: (l_1,\cdots,l_n)\mapsto \sum_{j=1}^n l_j\epsilon_j,$$
where $(\epsilon_1,\cdots,\epsilon_n)$ is the orthonormal base of
$\frak h^*$.
\end{remark}

\vskip 2mm From now on, we will focus on the classical Lie algebra
$\frak g=B_n$. The basic data for other types are presented in the
end of this section. Let $Y\in Y(\la)$ be a reduced proper Young
wall with a dominant integral weight $\la$. Then we can associate
a tableau $T_Y$ of shape $GRY^{-1}(\la)$ determined by the
following steps.

\vskip 2mm {\bf Step 1.}\, At first, consider ${\overline
Y}_{\omega_{i_k}}=(y^1_{\omega_{i_k}},\cdots,y^{i_k}_{\omega_{i_k}})$
in $Y$, where $y^a_{\omega_{i_k}}$ ($a=1,\cdots,i_k$) is the
$a$-th column of  ${\overline Y_{\omega_{i_k}}}$  from the right.
Then each column $y^a_{\omega_{i_k}}$ corresponds to a box with
entry, denoted by $T^a_{\omega_{i_k}}$, as follows:

%
\savebox{\tmpfiga}{\begin{texdraw} \fontsize{7}{7}\selectfont
\drawdim em \setunitscale 1.7 \move(-1 0)\lvec(-1 1)\lvec(0
1)\lvec(-1 0)\lvec(0 0)\lvec(0 1) \htext(-0.7 0.75){$0$}
\end{texdraw}}
\savebox{\tmpfigb}{\begin{texdraw} \fontsize{7}{7}\selectfont
\drawdim em \setunitscale 1.7 \move(-1 0)\lvec(-1 1)\lvec(0
1)\lvec(-1 0)\lvec(0 0)\lvec(0 1) \htext(-0.3 0.3){$0$}
\end{texdraw}}
\savebox{\tmpfigc}{\begin{texdraw} \fontsize{7}{7}\selectfont
\drawdim em \setunitscale 1.7 \move(0 0)\lvec(-1 0)\lvec(-1
1)\lvec(0 1)\lvec(0 0) \htext(-0.5 0.5){$1$}
\end{texdraw}}
\savebox{\tmpfigd}{\begin{texdraw} \fontsize{7}{7}\selectfont
\drawdim em \setunitscale 1.7 \move(-1 0)\lvec(-1 1)\lvec(0
1)\lvec(-1 0)\lvec(0 0)\lvec(0 1) \htext(-0.7 0.75){$1$}
\end{texdraw}}
\savebox{\tmpfige}{\begin{texdraw} \fontsize{7}{7}\selectfont
\drawdim em \setunitscale 1.7 \move(-1 0)\lvec(-1 1)\lvec(0
1)\lvec(-1 0)\lvec(0 0)\lvec(0 1) \htext(-0.3 0.3){$1$}
\end{texdraw}}
\savebox{\tmpfigf}{\begin{texdraw} \fontsize{7}{7}\selectfont
\drawdim em \setunitscale 1.7 \move(0 0)\lvec(-1 0)\lvec(-1
1)\lvec(0 1)\lvec(0 0) \htext(-0.5 0.5){$\bar{1}$}
\end{texdraw}}
\savebox{\tmpfigi}{\begin{texdraw} \fontsize{7}{7}\selectfont
\drawdim em \setunitscale 1.7 \move(0 0)\lvec(-1 0)\lvec(-1
1)\lvec(0 1)\lvec(0 0) \htext(-0.5 0.5){$2$}
\end{texdraw}}
\savebox{\tmpfigl}{\begin{texdraw}\fontsize{7}{7}\selectfont
\drawdim em \setunitscale 1.7 \move(0 0)\lvec(-1 0) \move(0
1)\lvec(-1 1)\move(0 0)\lvec(0 2)\lvec(-1 2)\lvec(-1 0)
\htext(-0.5 0.5){$j\!\!-\!\!2$} \htext(-0.5 1.5){$j\!\!-\!\!1$}
\end{texdraw}}
\savebox{\tmpfigm}{\begin{texdraw} \fontsize{7}{7}\selectfont
\drawdim em \setunitscale 1.7 \move(0 0)\lvec(-1 0)\lvec(-1
1)\lvec(0 1)\lvec(0 0) \htext(-0.5 0.5){$j$}
\end{texdraw}}
\savebox{\tmpfign}{\begin{texdraw}\fontsize{7}{7}\selectfont
\drawdim em \setunitscale 1.7 \move(0 0)\lvec(-1 0) \move(0
1)\lvec(-1 1)\move(0 0)\lvec(0 2)\lvec(-1 2)\lvec(-1 0)
\htext(-0.5 1.5){$j$} \htext(-0.5 0.5){$j\!\!+\!\!1$}
\end{texdraw}}
\savebox{\tmpfigo}{\begin{texdraw}\fontsize{7}{7}\selectfont
\drawdim em \setunitscale 1.7 \move(0 0)\lvec(-1 0)\lvec(-1
1)\lvec(0 1)\lvec(0 0) \htext(-0.5 0.5){$\bar{j}$}
\end{texdraw}}
\savebox{\tmpfigq}{\begin{texdraw} \fontsize{7}{7}\selectfont
\drawdim em \setunitscale 1.7 \move(0 0)\lvec(-1 0)\lvec(-1
1)\lvec(0 1)\lvec(0 0) \htext(-0.5 0.5){$0$}
\end{texdraw}}
\savebox{\tmpfigw}{\begin{texdraw} \fontsize{7}{7}\selectfont
\drawdim em \setunitscale 1.7 \move(0 0)\lvec(-1 0)\lvec(-1
1)\lvec(0 1)\lvec(0 0) \htext(-0.5 0.5){$\bar{n}$}
\end{texdraw}}
\savebox{\tmpfigx}{\begin{texdraw} \fontsize{7}{7}\selectfont
\drawdim em \setunitscale 1.7 \move(0 0)\lvec(0 1)\lvec(-1
1)\lvec(-1 0) \move(0 0)\lvec(-1 0)\lvec(0 1) \htext(-0.3
0.3){$1$} \htext(-0.7 0.75){$0$}
\end{texdraw}}
\savebox{\tmpfigy}{\begin{texdraw} \fontsize{7}{7}\selectfont
\drawdim em \setunitscale 1.7 \move(0 0)\lvec(0 1)\lvec(-1
1)\lvec(-1 0) \move(0 0)\lvec(-1 0)\lvec(0 1) \htext(-0.3
0.3){$0$} \htext(-0.7 0.75){$1$}
\end{texdraw}}
\savebox{\tmpfigad}{\begin{texdraw} \fontsize{7}{7}\selectfont
\drawdim em \setunitscale 1.7 \move(0 0)\lvec(-1 0)\lvec(-1
1)\lvec(0 1)\lvec(0 0) \htext(-0.5 0.5){$\overline{n\!\!-\!\!1}$}
\end{texdraw}}
\savebox{\tmpfigal}{\begin{texdraw} \fontsize{7}{7}\selectfont
\drawdim em \setunitscale 1.7 \move(0 0)\lvec(0 0.5)\lvec(-1
0.5)\lvec(-1 0)\lvec(0 0) \htext(-0.5 0.26){$n$}
\end{texdraw}}
\savebox{\tmpfigam}{\begin{texdraw}\fontsize{7}{7}\selectfont
\drawdim em \setunitscale 1.7 \move(0 0)\lvec(0 1)\lvec(-1
1)\lvec(-1 0)\lvec(0 0) \move(0 0.5)\lvec(-1 0.5) \htext(-0.5
0.26){$n$} \htext(-0.5 0.76){$n$}
\end{texdraw}}
\savebox{\tmpfigap}{\begin{texdraw} \fontsize{7}{7}\selectfont
\drawdim em \setunitscale 1.7 \move(0 0)\lvec(0 1)\lvec(-1
1)\lvec(-1 0) \move(0 0)\lvec(-1 0)\lvec(0 1) \move(-1 1)\lvec(-1
2)\lvec(0 2)\lvec(0 1) \htext(-0.3 0.3){$1$} \htext(-0.7
0.75){$0$} \htext(-0.5 1.5){$2$}
\end{texdraw}}
\savebox{\tmpfigaq}{\begin{texdraw} \fontsize{7}{7}\selectfont
\drawdim em \setunitscale 1.7 \move(0 0)\lvec(0 1)\lvec(-1
1)\lvec(-1 0) \move(0 0)\lvec(-1 0)\lvec(0 1) \move(-1 1)\lvec(-1
2)\lvec(0 2)\lvec(0 1) \htext(-0.3 0.3){$0$} \htext(-0.7
0.75){$1$} \htext(-0.5 1.5){$2$}
\end{texdraw}}
\savebox{\tmpfigar}{\begin{texdraw}\fontsize{7}{7}\selectfont
\drawdim em \setunitscale 1.7 \move(0 0)\lvec(-1 0)\lvec(-1
1)\lvec(0 1)\lvec(0 0) \htext(-0.5 0.5){$3$}
\end{texdraw}}
\savebox{\tmpfigaw}{\begin{texdraw} \fontsize{7}{7}\selectfont
\drawdim em \setunitscale 1.7 \move(0 0)\lvec(-1 0)\lvec(-1
2)\lvec(0 2)\lvec(0 0) \move(-1 0.5)\lvec(0 0.5) \move(-1
1)\lvec(0 1) \htext(-0.5 0.26){$n$} \htext(-0.5 0.76){$n$}
\htext(-0.5 1.5){$n\!\!-\!\!1$}
\end{texdraw}}

\vskip 3mm
\begin{multicols}{2}
\hskip 5mm
\usebox{\tmpfiga} \araise{ or } \usebox{\tmpfigb}
\araise{$\longmapsto$} \usebox{\tmpfigc}\\[4pt]

\hskip 5mm
\usebox{\tmpfigx} \araise{ or } \usebox{\tmpfigy}
\araise{$\longmapsto$} \usebox{\tmpfigi}\\[4pt]

\hskip 5mm
\usebox{\tmpfigap} \braise{ or } \usebox{\tmpfigaq}
\braise{$\longmapsto$} \craise{\usebox{\tmpfigar}}\\[4pt]

\hskip 5mm
\usebox{\tmpfigl}
\braise{$\longmapsto$} \craise{\usebox{\tmpfigm}}
\quad \draise{($4\leq j \leq n$)}\\[4pt]

\hskip 5mm
\usebox{\tmpfigal}
\araise{$\longmapsto$} \usebox{\tmpfigq}\\[4pt]

\usebox{\tmpfigam}
\araise{$\longmapsto$} \usebox{\tmpfigw}\\[4pt]

\usebox{\tmpfigaw}
\braise{$\longmapsto$} \craise{\usebox{\tmpfigad}}\\[4pt]

\usebox{\tmpfign}
\braise{$\longmapsto$} \craise{\usebox{\tmpfigo}}
\quad \draise{($2\leq j \leq n-2$)}\\[4pt]

\usebox{\tmpfigd} \araise{ or } \usebox{\tmpfige}
\araise{$\longmapsto$} \usebox{\tmpfigf}
\end{multicols}

\noindent Here, the blocks of left-hand side represent the top
blocks in each column $y^a_{\omega_{i_k}}$ of ${\overline
Y}_{\omega_{i_k}}$ and the boxes with entry of right-hand side
represent $T^a_{\omega_{i_k}}$. Now, a tableau $T_{\omega_{i_k}}$
corresponded to ${\overline Y_{\omega_{i_k}}}$ is obtained by
stacking from $T^1_{\omega_{i_k}}$ to $T^{i_k}_{\omega_{i_k}}$.

\savebox{\tmpfiga}{\begin{texdraw} \fontsize{7}{7}\selectfont
\drawdim em \setunitscale 0.85

\move(-2 0)\rlvec(8 0)\rlvec(0 8)\move(2 0)\rlvec(0 4)\rlvec(4
0)\move(4 0) \rlvec(0 6)\rlvec(2 0)\move(0 0)\rlvec(2 2)\move(2
0)\rlvec(2 2)\move(4 0)\rlvec(2 2)\move(0 0)\rlvec(0 2)\rlvec(6 0)
\move(6 7)\rlvec(-2 0)\rlvec(0 -1)

\htext(0.5 1.5){$1$} \htext(1.5 0.5){$0$}\htext(2.5 1.5){$0$}
\htext(3.5 0.5){$1$}\htext(4.5 1.5){$1$}\htext(5.5 0.5){$0$}
\htext(3 3){$2$}\htext(5 3){$2$}\htext(5 5){$3$}\htext(5 6.5){$4$}
\end{texdraw}}
\savebox{\tmpfigb}{\begin{texdraw} \fontsize{7}{7}\selectfont
\drawdim em \setunitscale 0.85

\rlvec(2 0)\rlvec(0 6)\rlvec(-2 0)\rlvec(0 -6)\move(0 2)\rlvec(2
0)\move(0 4)\rlvec(2 0)

\htext(1 1){$0$} \htext(1 3){$3$}\htext(1 5){$2$}
\end{texdraw}}

\vskip 3mm
\begin{example}
If $\frak g=B_4$, $\la=\omega_3$ and
$$Y=\raisebox{-0.4\height}{\usebox{\tmpfiga}}\in Y(\omega_3), \quad
\text{then}\,\,\,
T_Y=\raisebox{-0.4\height}{\usebox{\tmpfigb}}\,.$$
\end{example}

\savebox{\spinaa}{\begin{texdraw}\fontsize{7}{7}\selectfont
\drawdim em \setunitscale 0.85

\move (-9 0)\lvec(0 0)\rlvec(0 5)\move(0 1)\rlvec(-8 0)\move(-2
1)\rlvec(0 -1)\move(-6 1)\rlvec(0 -1) \move(-8 1)\rlvec(0 -1)

\htext(-1 0.5){$n$}\htext(-4 0.5){$\cdots$}\htext(-7 0.5){$n$}
\end{texdraw}}
\savebox{\spinab}{\begin{texdraw}\fontsize{7}{7}\selectfont
\drawdim em \setunitscale 0.85

\rlvec(0 8)\rlvec(1 0)\rlvec(0 -8)\rlvec(-1 0)\move(0 2)\rlvec(1
0)\move(0 6)\rlvec(1 0)

\htext(0.5 1){$n$}\htext(0.5 4){$\vdots$}\htext(0.5 7){$1$}
\end{texdraw}}
\savebox{\spinac}{\begin{texdraw}\fontsize{7}{7}\selectfont
\drawdim em \setunitscale 0.85

\rlvec(0 2)\rlvec(2 0)\rlvec(0 -2)\rlvec(-2 0)\move(0 1)\rlvec(2
0)

\htext(1 0.5){$n$}\htext(1 1.5){$n$}
\end{texdraw}}
\savebox{\spinad}{\begin{texdraw}\fontsize{7}{7}\selectfont
\drawdim em \setunitscale 0.85

\rlvec(0 2)\rlvec(1 0)\rlvec(0 -2)\rlvec(-1 0)

\htext(0.5 1){$\overline{n}$}
\end{texdraw}}
\savebox{\spinae}{\begin{texdraw}\fontsize{7}{7}\selectfont
\drawdim em \setunitscale 0.85

\rlvec(0 4)\rlvec(2 0)\rlvec(0 -4)\rlvec(-2 0)\move(0 2)\rlvec(2
0)

\htext(1 1){$j\!\!+\!\!1$}\htext(1 3){$j$}
\end{texdraw}}
\savebox{\spinaf}{\begin{texdraw}\fontsize{7}{7}\selectfont
\drawdim em \setunitscale 0.85

\rlvec(0 2)\rlvec(1 0)\rlvec(0 -2)\rlvec(-1 0)

\htext(0.5 1){$\overline{j}$}
\end{texdraw}}
\savebox{\spinag}{\begin{texdraw}\fontsize{7}{7}\selectfont
\drawdim em \setunitscale 0.85

\rlvec(0 4)\rlvec(2 0)\rlvec(0 -4)\rlvec(-2 0)\move(0 2)\rlvec(2
0)\move(0 1)\rlvec(2 0)

\htext(1 0.5){$n$}\htext(1 1.5){$n$} \htext(1 3){$n\!\!-\!\!1$}
\end{texdraw}}
\savebox{\spinah}{\begin{texdraw}\fontsize{7}{7}\selectfont
\drawdim em \setunitscale 0.85

\rlvec(0 2)\rlvec(1 0)\rlvec(0 -2)\rlvec(-1 0)

\htext(0.5 1){$\overline{n\!\!-\!\!1}$}
\end{texdraw}}
\savebox{\spinai}{\begin{texdraw}\fontsize{7}{7}\selectfont
\drawdim em \setunitscale 0.85

\rlvec(0 2)\rlvec(2 0)\rlvec(0 -2)\rlvec(-2 0)\rlvec(2 2)

\htext(0.75 1.25){$1$}
\end{texdraw}}
\savebox{\spinaj}{\begin{texdraw}\fontsize{7}{7}\selectfont
\drawdim em \setunitscale 0.85

\rlvec(0 2)\rlvec(2 0)\rlvec(0 -2)\rlvec(-2 0)\rlvec(2 2)

\htext(1.25 0.75){$1$}
\end{texdraw}}
\savebox{\spinak}{\begin{texdraw}\fontsize{7}{7}\selectfont
\drawdim em \setunitscale 0.85

\rlvec(0 2)\rlvec(1 0)\rlvec(0 -2)\rlvec(-1 0)

\htext(0.5 1){$\overline{1}$}
\end{texdraw}}

\vskip 3mm {\bf Step 2.}\, Consider
$\overline{Y}_{\la_n}=(y^1_{\la_n},\cdots,y^n_{\la_n})$, where
$y^a_{\la_n}$ ($a=1,,\cdots,n$) is the $a$-th column of
$\overline{Y}_{\la_n}$ from the right. On the one hand, if
${\overline Y_{\la_n}}=H_{\la_n}$, then it corresponds to a
tableau $\raisebox{-0.4\height}{\usebox{\spinab}}\,.$

Suppose that there are some added blocks above $H_{\la_n}$ in
${\overline Y_{\lambda_n}}$. Then the columns having added blocks
above $H_{\la_n}$ correspond to boxes with entries as follows.
$$\begin{array}{ll}
\raisebox{-0.25\height}{\usebox{\spinac}}\longmapsto
\raisebox{-0.25\height}{\usebox{\spinad}}&\qquad\quad
\raisebox{-0.4\height}{\usebox{\spinae}}\longmapsto
\raisebox{-0.25\height}{\usebox{\spinaf}}\,\,\, (2\le j\le n-2)\\
\raisebox{-0.4\height}{\usebox{\spinag}}\longmapsto
\raisebox{-0.4\height}{\usebox{\spinah}}&\qquad\quad
\raisebox{-0.25\height}{\usebox{\spinai}}\,\,\,\text{or}\,\,\,
\raisebox{-0.25\height}{\usebox{\spinaj}}\longmapsto
\raisebox{-0.25\height}{\usebox{\spinak}}
\end{array}$$
Moreover, if $k$-many columns with added blocks above $H_{\la_n}$
in ${\overline Y_{\lambda_n}}$ correspond to $k$-many boxes with
entries $i_1,\cdots,i_k\in\{\overline{n},\cdots,\overline{1}\}$,
then ($n-k$)-many columns with no added block above $H_{\la_n}$
correspond to ($n-k$)-many boxes with entry on
$\{1,\cdots,n\}\setminus\{\overline{i_1},\cdots,\overline{i_k}\}$.
Finally, ${\overline Y}_{\la_n}$ corresponds to a tableau of half
column of length $n$ obtained by stacking $T^a_{\la_n}$
($a=1,\cdots,n$) for the entries of $T^a_{\la_n}$ to strictly
increase down the columns.

\begin{remark}
$T(\la_n)$ is just the tableau realization of $B_{sp}$ given by
Kashiwara and Nakashima.
\end{remark}

%
\savebox{\spinaf}{
\begin{texdraw}\fontsize{7}{7}\selectfont
\drawdim em \setunitscale 0.85 \move (-8 0)\lvec(0 0)\rlvec(0
7)\move(0 1)\rlvec(-7 0)\move(-2 1)\rlvec(0 -1)\move(-4 1)\rlvec(0
-1)\move(-6 1)\rlvec(0 -1)\move(0 2)\rlvec(-2 0)\rlvec(0
-1)\move(0 4)\rlvec(-2 0)\rlvec(0 -2)\rlvec(-2 0)\rlvec(0
-1)\move(0 6)\rlvec(-2 -2)\rlvec(0 -2)\move(0 6)\rlvec(-2
0)\rlvec(0 -2)

\htext(-1 0.5){$3$}\htext(-3 0.5){$3$}\htext(-5 0.5){$3$}\htext(-1
1.5){$3$}\htext(-1 3){$2$}\htext(-0.7 4.3){$1$}\htext(-3 1.5){$3$}
\end{texdraw}
}
\savebox{\spinbf}{\begin{texdraw}\fontsize{7}{7}\selectfont
\drawdim em \setunitscale 0.8 \lvec(0 6)\rlvec(-1 0)\rlvec(0
-6)\rlvec(1 0)\move(0 2)\rlvec(-1 0)\move(0 4)\rlvec(-1 0)
\htext(-0.5 1){$\overline{1}$}\htext(-0.5 3){$\overline{3}$}
\htext(-0.5 5){$2$}
\end{texdraw}}
\savebox{\spinag}{\begin{texdraw}\fontsize{7}{7}\selectfont
\drawdim em \setunitscale 0.8

\lvec(0 2)\rlvec(-1 0)\rlvec(0 -2)\rlvec(1 0)

\htext(-0.5 1){$\overline{1}$}
\end{texdraw}}
\savebox{\spinah}{\begin{texdraw}\fontsize{7}{7}\selectfont
\drawdim em \setunitscale 0.8

\lvec(0 2)\rlvec(-1 0)\rlvec(0 -2)\rlvec(1 0)

\htext(-0.5 1){$\overline{3}$}
\end{texdraw}}

\vskip 3mm \begin{example}\, Let $\frak g=B_3$, $\la=\la_3$ and
$$Y=(y^1_{\la_3},y^2_{\la_3},y^3_{\la_3})=\raisebox{-0.4\height}{\usebox{\spinaf}}\,.$$
Since $y^1_{\la_3}$ and $y^2_{\la_3}$ correspond to
$\raisebox{-0.15\height}{\usebox{\spinag}}$ and
$\raisebox{-0.15\height}{\usebox{\spinah}}$\,, respectively, we
have
$$T_Y=\raisebox{-0.4\height}{\usebox{\spinbf}}\,.$$

\end{example}

\vskip 3mm {\bf Step 3.}\, Finally, a tableau $T_Y$ corresponded
to $Y\in Y(\la)$ of shape $GRY^{-1}(\la)$ is obtained by attaching
$T_{\omega_{i_k}}$ ($k=1,\cdots,t-1$) and $T_{\omega_{i_t}}$ to
the left-hand side of $T_{\omega_{i_{k+1}}}$ and $T_{\la_n}$ to be
bottom-justified.

%
\savebox{\tmpfiga}{\begin{texdraw} \fontsize{7}{7}\selectfont
\drawdim em \setunitscale 0.85 \nc{\dtri}{ \bsegment \lvec(-2 0)
\lvec(-2 2)\lvec(0 2)\lvec(0 0)\ifill f:0.7 \esegment }
\nc{\dtrii}{ \bsegment \lvec(-2 0) \lvec(-2 1)\lvec(0 1)\lvec(0
0)\ifill f:0.7 \esegment } \bsegment \setgray 0.6 \move(-3 0)
\lvec(18 0)\lvec(18 25.5) \move(18 24)\rlvec(-2 0)\rlvec(0
-24)\move(18 22)\rlvec(-4 0)\rlvec(0 -22) \move(18 20)\rlvec(-4
0)\move(18 18)\rlvec(-4 0)\move(18 16)\rlvec(-4 0)\move(18
14)\rlvec(-4 0)\move(18 13)\rlvec(-8 0)\rlvec(0 -13)\move(18
12)\rlvec(-8 0) \move(18 10)\rlvec(-8 0)\move(18 8)\rlvec(-10
0)\rlvec(0 -8) \move(18 6)\rlvec(-12 0) \rlvec(0 -6)\move(18 4)
\rlvec(-14 0)\rlvec(0 -4)\move(18 2)\rlvec(-16 0)\rlvec(0
-2)\move(18 1)\rlvec(-20 0)\move(12 13)\rlvec(0 -13) \move(18
20)\rlvec(-2 -2)\move(16 20)\rlvec(-2 -2) \move(18 8)\rlvec(-2
-2)\move(16 8)\rlvec(-2 -2)\move(14 8)\rlvec(-2 -2)\move(12
8)\rlvec(-2 -2)\move(10 8)\rlvec(-2 -2)

\move(18 20)\linewd 0.15 \setgray 0.4 \rlvec(-2 0)\rlvec(-2 -2)
\rlvec(0 -8)\rlvec(-2 0)\rlvec(0 -2)\rlvec(-4 0) \rlvec(0
-7)\rlvec(-8 0)\rlvec(0 -1)
\htext(17 23){$3$}\htext(17 21){$2$}\htext(15 21){$2$} \htext(17.5
18.5){$1$}\htext(16.5 19.5){$0$}\htext(15.5 18.5){$0$}\htext(14.5
19.5){$1$} \htext(17 17){$2$}\htext(15 17){$2$} \htext(17 15){$3$}
\htext(15 15){$3$} \htext(17 13.5){$4$} \htext(15
13.5){$4$}\htext(17 12.5){$4$} \htext(15 12.5){$4$}\htext(13
12.5){$4$} \htext(11 12.5){$4$} \htext(17 11){$3$} \htext(15
11){$3$}\htext(13 11){$3$} \htext(11 11){$3$} \htext(17
9){$2$}\htext(15 9){$2$}\htext(13 9){$2$}\htext(11 9){$2$}
\htext(17.5 6.5){$1$} \htext(16.5 7.5){$0$} \htext(15.5 6.5){$0$}
\htext(14.5 7.5){$1$} \htext(13.5 6.5){$1$} \htext(12.5
7.5){$0$}\htext(11.5 6.5){$0$} \htext(10.5 7.5){$1$}\htext(8.5
7.5){$0$}\htext(9.5 6.5){$1$} \htext(17 5){$2$} \htext(15 5){$2$}
\htext(13 5){$2$}\htext(11 5){$2$}\htext(9 5){$2$} \htext(7
5){$2$}\htext(17 3){$3$} \htext(15 3){$3$} \htext(13
3){$3$}\htext(11 3){$3$}\htext(9 3){$3$} \htext(7 3){$3$}\htext(5
3){$3$} \htext(17 1.5){$4$} \htext(15 1.5){$4$} \htext(13
1.5){$4$}\htext(11 1.5){$4$}\htext(9 1.5){$4$} \htext(7
1.5){$4$}\htext(5 1.5){$4$}\htext(3 1.5){$4$}\htext(17 0.5){$4$}
\htext(15 0.5){$4$} \htext(13 0.5){$4$}\htext(11 0.5){$4$}\htext(9
0.5){$4$}\htext(7 0.5){$4$}\htext(5 0.5){$4$}\htext(3
0.5){$4$}\htext(1 0.5){$4$}
\esegment
\end{texdraw}}
\savebox{\tmpfigb}{\begin{texdraw} \fontsize{7}{7}\selectfont
\drawdim em \setunitscale 0.85 \nc{\dtri}{ \bsegment \lvec(-2 0)
\lvec(-2 2)\lvec(0 2)\lvec(0 0)\ifill f:0.7 \esegment }
\nc{\dtrii}{ \bsegment \lvec(-2 0) \lvec(-2 1)\lvec(0 1)\lvec(0
0)\ifill f:0.7 \esegment } \rlvec(5 0)\rlvec(0 8) \rlvec(-1
0)\rlvec(0 -8)\move(5 6)\rlvec(-3 0)\rlvec(0 -6) \move(5
4)\rlvec(-5 0)\rlvec(0 -4)\move(5 2)\rlvec(-5 0) \htext(4.5
7){$1$}\htext(4.5 5){$\bar{4}$}\htext(3 5){$2$}\htext(4.5
3){$\bar{3}$} \htext(3 3){$0$} \htext(1 3){$3$}\htext(4.5
1){$\bar{2}$} \htext(3 1){$0$} \htext(1 1){$\bar{3}$}
\end{texdraw}}

\vskip 3mm
\begin{example}  If $\frak g=B_4$, $\la=\omega_2+\omega_3+\lambda_4$
and
\begin{center}
$Y=\raisebox{-0.4\height}{\usebox{\tmpfiga}}\in Y(\la)$, then
$T_Y=\raisebox{-0.4\height}{\usebox{\tmpfigb}}$\,.
\end{center}
\end{example}


%
\savebox{\tmpfiga}{\begin{texdraw} \fontsize{9}{9}\selectfont
\drawdim em \setunitscale 0.9 \rlvec(4 0)\rlvec(0 -2) \move(2
0)\rlvec(0 -2)\move(0 -1)\rlvec(4 0)\move(0 0)\rlvec(0 -2)
\htext(1 -0.5){$n$}\htext(3 -0.5){$n$}
\end{texdraw}}

\vskip 3mm Until now, we define a new tableau $T_Y$ of shape
$GRY^{-1}(\la)$ associated with $Y\in Y(\la)$ with a dominant
integral weight $\la$. We will characterize these new tableaux by
the conditions in $Y(\la)$.

At first, as we have seen, if we let $Y\in Y(\la)$, then
${\overline
Y_{\omega_{i_k}}}=(y^1_{\omega_{i_k}},\cdots,y^{i_k}_{\omega_{i_k}})$
corresponds to one column $T_{\omega_{i_k}}$ of $T_Y$. Moreover,
one column $y^a_{\omega_{i_k}}$ ($a=1,\cdots,i_k$) in ${\overline
Y_{\omega_{i_k}}}$ corresponds to one box $T^a_{\omega_{i_k}}$ in
$T_{\omega_{i_k}}$. Then since $Y$ is proper and
$T_{\omega_{i_k}}$ is obtained by stacking from
$T^1_{\omega_{i_k}}$ to $T^{i_k}_{\omega_{i_k}}$, the entries of
the tableau $T_Y$ strictly increase down the columns. Moreover, we
can stack $n$-block as $\raisebox{-0.5\height}{\usebox{\tmpfiga}}$
in ${\overline Y_{\omega_{i_k}}}$ by the pattern of stacking
blocks because the $n$-block is a block of half-unit-height and
$Y$ is proper. So, the 0 element can appear more than once in
tableau $T_{\omega_{i_k}}$.

\vskip 2mm On the other hand, consider the condition that
$Y^{\omega_{i_k}} \subset Y^{\omega_{i_{k+1}}}\,$ and
$Y^{\omega_{i_t}} \subset Y^{\la_n}$. Note that
$y^a_{\omega_{i_k}}$ (resp. $y^a_{\omega_{i_t}}$) and
$y^a_{\omega_{i_{k+1}}}$ (resp. $y^a_{\la_n}$) of ${\overline
Y_{\omega_{i_k}}}$ (resp. ${\overline Y_{\omega_{i_t}}}$) and
${\overline Y_{\omega_{i_{k+1}}}}$ (resp. ${\overline Y_{\la_n}}$)
correspond to the boxes with entries in the same row in
$T_{\omega_{i_k}+\omega_{i_{k+1}}}$ (resp.
$T_{\omega_{i_t}+\la_n}$). Therefore, this condition implies that
the entries are weakly increasing in each rows of $T_Y$.

\savebox{\tmpfiga}{\begin{texdraw} \fontsize{7}{7}\selectfont
\drawdim em \setunitscale 0.8

\rlvec(2 0)\rlvec(0 11) \rlvec(-2 0)\rlvec(0 -11)\move(0
2)\rlvec(2 0)\move(0 9)\rlvec(2 0)

\move(0 9)\clvec (-0.5 9)(-1 7)(-1 6.5)

\htext(-2.5 5.5){$i_k\!\!\!-\!\!2\ge$}

\move(0 2)\clvec (-0.5 2)(-1 4)(-1 4.5)

\htext(1 1){$\overline{1}$}\htext(1 10){$1$}
\end{texdraw}}
\savebox{\tmpfigb}{\begin{texdraw} \fontsize{7}{7}\selectfont
\drawdim em \setunitscale 0.8

\rlvec(2 0)\rlvec(0 9) \rlvec(-2 0)\rlvec(0 -9)\move(0 2)\rlvec(2
0)\move(0 7)\rlvec(2 0)

\move(0 7)\clvec (-0.5 7)(-1 6)(-1 5.5)

\htext(-2.5 4.5){$i_k\!\!\!-\!\!3\ge$}

\move(0 2)\clvec (-0.5 2)(-1 3)(-1 3.5)

 \htext(1
1){$\overline{2}$}\htext(1 8){$2$}
\end{texdraw}}
\savebox{\tmpfigc}{\begin{texdraw} \fontsize{7}{7}\selectfont
\drawdim em \setunitscale 0.8

\rlvec(2 0)\rlvec(0 7) \rlvec(-2 0)\rlvec(0 -7)\move(0 2)\rlvec(2
0)\move(0 5)\rlvec(2 0)

\move(0 5)\clvec (-0.3 5)(-0.7 4.6)(-0.7 4.5)

\htext(-3.5 3.5){$i_k\!\!-\!\!(a\!\!\!+\!\!1)\ge$}

\move(0 2)\clvec (-0.3 2)(-0.7 2.4)(-0.7 2.5)

\htext(1 1){$\overline{a}$}\htext(1 6){$a$}
\end{texdraw}}
\savebox{\tmpfigd}{\begin{texdraw} \fontsize{5}{5}\selectfont
\drawdim em \setunitscale 1

\rlvec(2 0)\rlvec(0 4) \rlvec(-2 0)\rlvec(0 -4)\move(0 2)\rlvec(2
0)

\htext(1 1){$\overline{i_k\!\!\!-\!\!1}$}\htext(1
3){$i_k\!\!\!-\!\!1$}
\end{texdraw}}

\vskip 3mm Consider the condition {\bf (Y1)}
$Y_{\omega_{i_k}}^+\subset |Y_{\omega_{i_k}}^-|$ ($1\leq i_k \leq
n$) in Proposition 2.9. Then this condition implies that the
following can not exist in $T_{\omega_{i_k}}$ as a subtableau :
$$
\raisebox{-0.4\height}{\usebox{\tmpfiga}} \hskip 6mm
\raisebox{-0.4\height}{\usebox{\tmpfigb}} \hskip 3mm \cdots \hskip
10mm \raisebox{-0.4\height}{\usebox{\tmpfigc}}\hskip
3mm\cdots\hskip 3mm \raisebox{-0.4\height}{\usebox{\tmpfigd}}
$$
Therefore, if $T_{\omega_{i_k}}$ has $a$-box and $\bar{a}$-box,
then the number of rows between $a$-box and $\bar{a}$-box is
larger than or equals to $i_k-a$. That is, if $a$-box and
$\bar{a}$-box lie in $p$-th and $q$-th rows of $T_{\omega_{i_k}}$
from bottom to top with $p>q$, we have
\begin{equation}
p-q-1\ge i_k-a.
\end{equation}

\savebox{\tmpfiga}{\begin{texdraw} \fontsize{7}{7}\selectfont
\drawdim em \setunitscale 0.8

\rlvec(2 0)\rlvec(0 10) \rlvec(-2 0)\rlvec(0 -10)\move(0
2)\rlvec(2 0)\move(0 8)\rlvec(2 0)

\move(0 8)\clvec (-0.5 8)(-1 6.5)(-1 6)

\htext(-2.5 5){$i_k\!\!-\!\!(a\!\!\!+\!\!1)$}

\move(0 2)\clvec (-0.5 2)(-1 3.5)(-1 4)

\htext(1 1){$\overline{a}$}\htext(1 9){$a$}
\end{texdraw}}
\savebox{\tmpfigb}{\begin{texdraw} \fontsize{7}{7}\selectfont
\drawdim em \setunitscale 0.8

\rlvec(2 0)\rlvec(0 10) \rlvec(-2 0)\rlvec(0 -10)\move(0
2)\rlvec(2 0)\move(0 8)\rlvec(2 0)

\htext(1 1){$\overline{1}$}\htext(1 9){$a$}
\end{texdraw}}
\savebox{\tmpfigc}{\begin{texdraw} \fontsize{7}{7}\selectfont
\drawdim em \setunitscale 0.8

\rlvec(2 0)\rlvec(0 8) \rlvec(-2 0)\rlvec(0 -8)\move(0 2)\rlvec(2
0)\move(0 6)\rlvec(2 0)

\move(0 6)\clvec (-0.5 6)(-1 5.3)(-1 5)

\htext(-2.5 4){$i_k\!\!-\!\!(a\!\!\!+\!\!2)$}

\move(0 2)\clvec (-0.5 2)(-1 2.7)(-1 3)

\htext(1 1){$\overline{a\!\!+\!\!1}$}\htext(1 7){$a$}
\end{texdraw}}
\savebox{\tmpfigd}{\begin{texdraw} \fontsize{7}{7}\selectfont
\drawdim em \setunitscale 0.8

\rlvec(2 0)\rlvec(0 8) \rlvec(-2 0)\rlvec(0 -8)\move(0 2)\rlvec(2
0)\move(0 6)\rlvec(2 0)

\htext(1 1){$\overline{2}$}\htext(1 7){$a$}
\end{texdraw}}
\savebox{\tmpfige}{\begin{texdraw} \fontsize{5}{5}\selectfont
\drawdim em \setunitscale 1

\rlvec(2 0)\rlvec(0 4) \rlvec(-2 0)\rlvec(0 -4)\move(0 2)\rlvec(2
0)

\htext(1 1){$\overline{i_k\!\!\!-\!\!1}$}\htext(1 3){$a$}
\end{texdraw}}
\savebox{\tmpfigf}{\begin{texdraw} \fontsize{5}{5}\selectfont
\drawdim em \setunitscale 1

\rlvec(2 0)\rlvec(0 4) \rlvec(-2 0)\rlvec(0 -4)\move(0 2)\rlvec(2
0)

\htext(1 1){$\overline{i_k\!\!\!-\!\!a}$}\htext(1 3){$a$}
\end{texdraw}}

\begin{remark} (a) In fact, we know that the following also can not exist
in $T_{\omega_{i_k}}$ as a subtableau.
$$
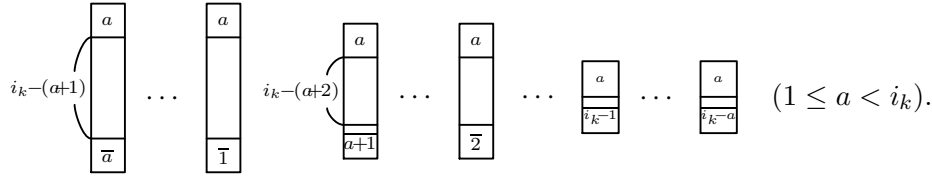
\raisebox{-0.4\height}{\usebox{\tmpfiga}} \hskip 2mm\cdots\hskip
3mm \raisebox{-0.4\height}{\usebox{\tmpfigb}} \hskip 8mm
\raisebox{-0.4\height}{\usebox{\tmpfigc}} \hskip 2mm\cdots\hskip
3mm \raisebox{-0.4\height}{\usebox{\tmpfigd}}\hskip
3mm\cdots\hskip 3mm
\raisebox{-0.4\height}{\usebox{\tmpfige}}\hskip 2mm\cdots\hskip
3mm \raisebox{-0.4\height}{\usebox{\tmpfigf}}\,\quad (1\le a<
i_k).
$$
But, it is evident that this condition is subject to above
condition by the pattern of stacking the blocks in Young walls and
proper condition of $Y$. Anyway, if $T_{\omega_{i_k}}$ has $a$-box
and $\bar{b}$-box, then the number of rows between $a$-box and
$\bar{b}$-box is larger than or equals to $i_k-\text{max}(a,b)$.
That is, if $a$-box and $\bar{b}$-box lie in $p$-th and $q$-th
rows of $T_{\omega_{i_k}}$ from bottom to top with $p>q$, we have
\begin{equation}
p-q-1\ge i_k-\text{max}(a,b).
\end{equation}

(b) It is just the one column condition (1CC) of
Kashiwara-Nakashima tableaux. We will also call it the {\it one
column condition} and  denote by (1CC).
\end{remark}

\savebox{\tmpfiga}{
\begin{texdraw}
\fontsize{7}{7}\selectfont \drawdim em \setunitscale 0.8

\move(0 0)\rlvec(4 0)\rlvec(0 11)\rlvec(-4 0)\rlvec(0 -11)\move(0
2)\rlvec(4 0)\move(0 9)\rlvec(4 0)\move(2 0)\rlvec(0 11)

\move(0 9)\clvec (-0.5 9)(-1 7)(-1 6.5)

\htext(-3 5.5){$n\!\!\!-\!\!a\le$}

\move(0 2)\clvec (-0.5 2)(-1 4)(-1 4.5)

\htext(1 10){$a$}\htext(3 1){$\overline{a}$}
\end{texdraw}}
\savebox{\tmpfigb}{
\begin{texdraw}
\fontsize{7}{7}\selectfont \drawdim em \setunitscale 0.8

\move(0 0)\rlvec(4 0)\rlvec(0 9)\rlvec(-4 0)\rlvec(0 -9)\move(0
2)\rlvec(4 0)\move(0 7)\rlvec(4 0)\move(2 0)\rlvec(0 9)

\move(0 7)\clvec (-0.5 7)(-1 6)(-1 5.5)

\htext(-3.5 4.5){$n\!\!\!-\!\!(a\!\!+\!\!1)\le$}

\move(0 2)\clvec (-0.5 2)(-1 3)(-1 3.5)

\htext(1 8){$a\!\!+\!\!1$}\htext(3 1){$\overline{a\!\!+\!\!1}$}
\end{texdraw}}
\savebox{\tmpfigc}{
\begin{texdraw}
\fontsize{7}{7}\selectfont \drawdim em \setunitscale 0.8

\move(0 0)\rlvec(4 0)\rlvec(0 4)\rlvec(-4 0)\rlvec(0 -4)\move(0
2)\rlvec(4 0)\move(2 0)\rlvec(0 4)

\htext(1 3){$n$}\htext(3 1){$\overline{n}$}
\end{texdraw}}
\savebox{\tmpfigd}{
\begin{texdraw}
\fontsize{7}{7}\selectfont \drawdim em \setunitscale 0.8

\move(0 0)\rlvec(3 0)\rlvec(0 11)\rlvec(-3 0)\rlvec(0 -11)\move(0
2)\rlvec(3 0)\move(0 9)\rlvec(3 0)\move(2 0)\rlvec(0 11)

\move(0 9)\clvec (-0.5 9)(-1 7)(-1 6.5)

\htext(-3 5.5){$n\!\!\!-\!\!a\le$}

\move(0 2)\clvec (-0.5 2)(-1 4)(-1 4.5)

\htext(1 10){$a$}\htext(2.5 1){$\overline{a}$}
\end{texdraw}}
\savebox{\tmpfige}{
\begin{texdraw}
\fontsize{7}{7}\selectfont \drawdim em \setunitscale 0.8

\move(0 0)\rlvec(3 0)\rlvec(0 9)\rlvec(-3 0)\rlvec(0 -9)\move(0
2)\rlvec(3 0)\move(0 7)\rlvec(3 0)\move(2 0)\rlvec(0 9)

\move(0 7)\clvec (-0.5 7)(-1 6)(-1 5.5)

\htext(-3.5 4.5){$n\!\!\!-\!\!(a\!\!+\!\!1)\le$}

\move(0 2)\clvec (-0.5 2)(-1 3)(-1 3.5)

\move(4.5 1)\ravec(-2 0)

\htext(1 8){$a\!\!+\!\!1$}\htext(6 1){$\overline{a\!\!+\!\!1}$}
\end{texdraw}}
\savebox{\tmpfigf}{
\begin{texdraw}
\fontsize{7}{7}\selectfont \drawdim em \setunitscale 0.8

\move(0 0)\rlvec(3 0)\rlvec(0 4)\rlvec(-3 0)\rlvec(0 -4)\move(0
2)\rlvec(3 0)\move(2 0)\rlvec(0 4)

\htext(1 3){$n$}\htext(2.5 1){$\overline{n}$}
\end{texdraw}}
\savebox{\tmpfigg}{
\begin{texdraw}
\fontsize{7}{7}\selectfont \drawdim em \setunitscale 0.8

\move(0 0)\rlvec(4 0)\rlvec(0 2)\rlvec(-4 0)\rlvec(0 -2)\move(2
0)\rlvec(0 2)

\htext(1 1){$0$}\htext(3 1){$0$}
\end{texdraw}}

\vskip 3mm Consider the condition {\bf (Y3)} in Proposition 2.9.
That is,
$$|Y_{\omega_{i_k}}^-| \subset
Y_{\omega_{i_{k+1}}}^+\,\,\,\,\text{and}\,\,\,\,
|Y_{\omega_{i_t}}^-| \subset Y_{\la_n}^+.
$$
This condition says that the following can not exist in
$T_{\omega_{i_k}+\omega_{i_{k+1}}}$ ($k=1,\cdots,t-1$) and
$T_{\omega_{i_t}+\La_{n}}$ as a subtableau :
$$
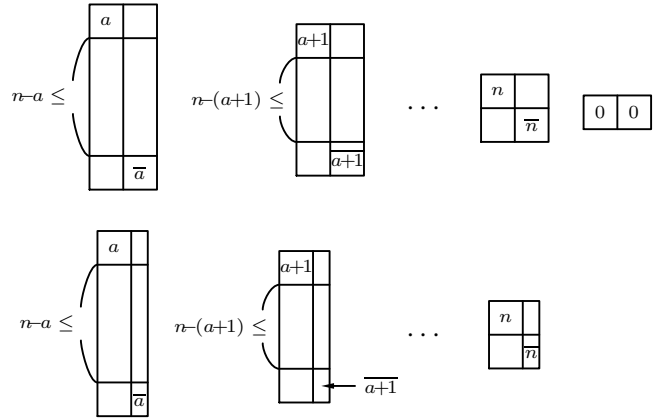
\begin{array}{ccccc}
\raisebox{-0.4\height}{\usebox{\tmpfiga}}\quad
&\raisebox{-0.4\height}{\usebox{\tmpfigb}}&\cdots
&\raisebox{-0.4\height}{\usebox{\tmpfigc}}
&\raisebox{-0.4\height}{\usebox{\tmpfigg}}\\ \\
\raisebox{-0.4\height}{\usebox{\tmpfigd}}\quad
&\raisebox{-0.4\height}{\usebox{\tmpfige}}&\cdots
&\raisebox{-0.4\height}{\usebox{\tmpfigf}}&
\end{array}
$$
Here, $n-i_k+2\le a\le n$ and $n-i_t+2\le a\le n$, respectively.
Therefore, if $T_{\omega_{i_k}}$ (resp. $T_{\omega_{i_t}}$) and
$T_{\omega_{i_{k+1}}}$ (resp. $T_{\la_{n}}$) have $a$ and
$\bar{a}$, respectively, then the number of rows between $a$-box
and $\overline{a}$-box is smaller than $n-a$. That is, if $a$-box
and $\overline{a}$-box lie in $p$-th and $q$-th rows of
$T_{\omega_{i_k}}$ (resp. $T_{\omega_{i_t}}$) and
$T_{\omega_{i_{k+1}}}$  (resp. $T_{\la_{n}}$) from bottom to top
with $p>q$, we have $p-q-1<n-a$ and so
\begin{equation}
p-q\le n-a.
\end{equation}

\savebox{\tmpfiga}{
\begin{texdraw}
\fontsize{7}{7}\selectfont \drawdim em \setunitscale 0.8

\move(0 0)\rlvec(4 0)\rlvec(0 7)\rlvec(-4 0)\rlvec(0 -7)\move(0
2)\rlvec(4 0)\move(0 5)\rlvec(4 0)\move(2 0)\rlvec(0 7)

\move(0 5)\clvec (-0.5 5)(-1 4.3)(-1 4)

\htext(-1.5 3.5){$n\!\!-\!\!a$}

\move(0 2)\clvec (-0.5 2)(-1 2.7)(-1 3)

\htext(1 6){$a$}\htext(3 1){$\overline{a}$}
\end{texdraw}}
\savebox{\tmpfigb}{
\begin{texdraw}
\fontsize{7}{7}\selectfont \drawdim em \setunitscale 0.8

\move(0 0)\rlvec(4 0)\rlvec(0 7)\rlvec(-4 0)\rlvec(0 -7)\move(0
2)\rlvec(4 0)\move(0 5)\rlvec(4 0)\move(2 0)\rlvec(0 7)

\htext(1 6){$a$}\htext(3 1){$0$}
\end{texdraw}}
\savebox{\tmpfigc}{
\begin{texdraw}
\fontsize{7}{7}\selectfont \drawdim em \setunitscale 0.8

\move(0 0)\rlvec(4 0)\rlvec(0 9)\rlvec(-4 0)\rlvec(0 -9)\move(0
2)\rlvec(4 0)\move(0 7)\rlvec(4 0)\move(2 0)\rlvec(0 9)

\move(0 7)\clvec (-0.5 7)(-1 5.8)(-1 5.5)

\htext(-1 4.5){$n\!\!-\!\!(a\!\!\!-\!\!1)$}

\move(0 2)\clvec (-0.5 2)(-1 3.2)(-1 3.5)

\htext(1 8){$a$}\htext(3 1){$\overline{a\!\!-\!\!1}$}
\end{texdraw}}
\savebox{\tmpfigd}{
\begin{texdraw}
\fontsize{7}{7}\selectfont \drawdim em \setunitscale 0.8

\move(0 0)\rlvec(4 0)\rlvec(0 9)\rlvec(-4 0)\rlvec(0 -9)\move(0
2)\rlvec(4 0)\move(0 7)\rlvec(4 0)\move(2 0)\rlvec(0 9)

\htext(1 8){$a$}\htext(3 1){$0$}
\end{texdraw}}
\savebox{\tmpfige}{
\begin{texdraw}
\fontsize{7}{7}\selectfont \drawdim em \setunitscale 0.8

\move(0 0)\rlvec(4 0)\rlvec(0 12)\rlvec(-4 0)\rlvec(0 -12)\move(0
2)\rlvec(4 0)\move(0 10)\rlvec(4 0)\move(2 0)\rlvec(0 12)

\move(0 10)\clvec (-0.5 10)(-1 7.5)(-1 7)

\htext(-1.5 6){$i_k\!\!\!-\!\!2$}

\move(0 2)\clvec (-0.5 2)(-1 4.5)(-1 5)

\htext(1 11){$a$}\htext(3
1){$\overline{n\!\!-\!\!(i_k\!\!-\!\!2)}$}
\end{texdraw}}
\savebox{\tmpfigf}{
\begin{texdraw}
\fontsize{7}{7}\selectfont \drawdim em \setunitscale 0.8

\move(0 0)\rlvec(4 0)\rlvec(0 12)\rlvec(-4 0)\rlvec(0 -12)\move(0
2)\rlvec(4 0)\move(0 10)\rlvec(4 0)\move(2 0)\rlvec(0 12)

\htext(1 11){$a$}\htext(3 1){$0$}
\end{texdraw}}
\savebox{\tmpfigg}{
\begin{texdraw}
\fontsize{7}{7}\selectfont \drawdim em \setunitscale 0.8

\move(0 0)\rlvec(4 0)\rlvec(0 2)\rlvec(-4 0)\rlvec(0 -2)\move(2
0)\rlvec(0 2)

\htext(1 1){$0$}\htext(3 1){$0$}
\end{texdraw}}
\savebox{\tmpfigg}{
\begin{texdraw}
\fontsize{7}{7}\selectfont \drawdim em \setunitscale 0.8

\move(0 0)\rlvec(4 0)\rlvec(0 7)\rlvec(-4 0)\rlvec(0 -7)\move(0
2)\rlvec(4 0)\move(0 5)\rlvec(4 0)\move(2 0)\rlvec(0 7)

\htext(1 6){$a$}\htext(3 1){$\overline{n}$}
\end{texdraw}}
\savebox{\tmpfigh}{
\begin{texdraw}
\fontsize{7}{7}\selectfont \drawdim em \setunitscale 0.8

\move(0 0)\rlvec(4 0)\rlvec(0 9)\rlvec(-4 0)\rlvec(0 -9)\move(0
2)\rlvec(4 0)\move(0 7)\rlvec(4 0)\move(2 0)\rlvec(0 9)

\htext(1 8){$a$}\htext(3 1){$\overline{n}$}
\end{texdraw}}
\savebox{\tmpfigi}{
\begin{texdraw}
\fontsize{7}{7}\selectfont \drawdim em \setunitscale 0.8

\move(0 0)\rlvec(4 0)\rlvec(0 12)\rlvec(-4 0)\rlvec(0 -12)\move(0
2)\rlvec(4 0)\move(0 10)\rlvec(4 0)\move(2 0)\rlvec(0 12)

\htext(1 11){$a$}\htext(3 1){$\overline{n}$}
\end{texdraw}}
\savebox{\tmpfigj}{
\begin{texdraw}
\fontsize{7}{7}\selectfont \drawdim em \setunitscale 0.8

\move(0 0)\rlvec(4 0)\rlvec(0 2)\rlvec(-4 0)\rlvec(0 -2)\move(2
0)\rlvec(0 2)

\htext(1 1){$0$}\htext(3 1){$1$}
\end{texdraw}}
\savebox{\tmpfigk}{
\begin{texdraw}
\fontsize{7}{7}\selectfont \drawdim em \setunitscale 0.8

\move(0 0)\rlvec(4 0)\rlvec(0 4)\rlvec(-4 0)\rlvec(0 -4)\move(2
0)\rlvec(0 4)\move(0 2)\rlvec(4 0)

\move(5 1)\ravec(-2 0)

\htext(1 3){$0$}\htext(7.5 1){$0$ or $\overline{n}$}
\end{texdraw}}
\savebox{\tmpfigl}{
\begin{texdraw}
\fontsize{7}{7}\selectfont \drawdim em \setunitscale 0.8

\move(0 0)\rlvec(4 0)\rlvec(0 6)\rlvec(-4 0)\rlvec(0 -6)\move(2
0)\rlvec(0 6)\move(0 2)\rlvec(4 0)\move(0 4)\rlvec(4 0)

\move(5 1)\ravec(-2 0)

\htext(1 5){$0$}\htext(7.5 1){$0$ or $\overline{n}$}
\end{texdraw}}
\savebox{\tmpfigm}{
\begin{texdraw}
\fontsize{7}{7}\selectfont \drawdim em \setunitscale 0.8

\move(0 0)\rlvec(4 0)\rlvec(0 12)\rlvec(-4 0)\rlvec(0 -12)\move(2
0)\rlvec(0 12)\move(0 2)\rlvec(4 0)\move(0 10)\rlvec(4 0)

\move(5 1)\ravec(-2 0)

\move(0 10)\clvec (-0.5 10)(-1 7.5)(-1 7)

\htext(-1.5 6){$i_k\!\!\!-\!\!2$}

\move(0 2)\clvec (-0.5 2)(-1 4.5)(-1 5)

\htext(1 11){$0$}\htext(7.5 1){$0$ or $\overline{n}$}
\end{texdraw}}

\begin{remark} (a) In fact, we know that the following also can not exist
in $T_{\omega_{i_k}+\omega_{i_{k+1}}}$ as a subtableau.
$$
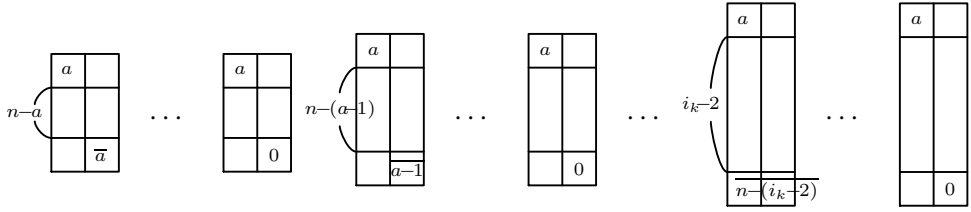
\begin{array}{cccccccccc}
\raisebox{-0.4\height}{\usebox{\tmpfiga}}\,&\cdots
&\raisebox{-0.4\height}{\usebox{\tmpfigb}}\,\,\,
&\raisebox{-0.4\height}{\usebox{\tmpfigc}}\,&\cdots
&\raisebox{-0.4\height}{\usebox{\tmpfigd}}\,&\cdots
&\raisebox{-0.4\height}{\usebox{\tmpfige}}\,&\cdots
&\raisebox{-0.4\height}{\usebox{\tmpfigf}}\\ \\
\end{array}
$$
where $n-i_k+2\le a\le n$. Moreover, the following can not exist
in $T_{\omega_{i_k}+\omega_{i_{k+1}}}$ as a subtableau.
$$
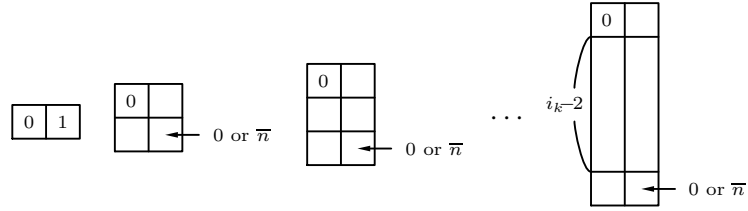
\begin{array}{ccccc}
\raisebox{-0.4\height}{\usebox{\tmpfigj}}
&\raisebox{-0.4\height}{\usebox{\tmpfigk}}\quad
&\raisebox{-0.4\height}{\usebox{\tmpfigl}}\quad &\cdots
&\raisebox{-0.4\height}{\usebox{\tmpfigm}}
\end{array}
$$
But, it is clear that this condition is subject to above condition
(1.7). Anyway, if $T_{\omega_{i_k}}$ and $T_{\omega_{i_{k+1}}}$
have $a$ and $\bar{b}$, respectively, the number of rows between
$a$-box and $\overline{b}$-box is smaller than
$n-\text{min}(a,b)$. That is, if $a$-box and $\overline{b}$-box
lie in $p$-th and $q$-th rows of $T_{\omega_{i_k}}$ and
$T_{\omega_{i_{k+1}}}$ from bottom to top with $p>q$, we have
$p-q-1< n-\text{min}(a,b)$ and so
\begin{equation}
p- q\le n-\text{min}(a,b).
\end{equation}
Similarly, it also holds for $T_{\omega_{i_t}+\lambda_n}$.

(b) This condition contains the condition very similar to the two
column condition (2CC) for the tableau of Kashiwara and Nakashima
in the $(a,n)$, $(n,n)$-configuration. We will call it the {\it
first two column condition} and  denote by (2CC-1).
\end{remark}

\savebox{\tmpfiga}{\begin{texdraw} \fontsize{7}{7}\selectfont
\drawdim em \setunitscale 0.9

\move(0 -0.5)\lvec(0 3)\rlvec(2 0)\rlvec(0 -3.5)\move(0
1.5)\rlvec(2 0) \htext(1 0.75){$a\!\!-\!\!2$}\htext(1
2.25){$a\!\!-\!\!1$}
\end{texdraw}}
\savebox{\tmpfigb}{\begin{texdraw} \fontsize{7}{7}\selectfont
\drawdim em \setunitscale 0.9 \move(0 -0.5)\lvec(0 3)\rlvec(2
0)\rlvec(0 -3.5)\move(0 1.5)\rlvec(2 0) \htext(1
0.75){$a\!\!+\!\!1$}\htext(1 2.25){$a$}
\end{texdraw}}
\vskip 3mm Consider the condition {\bf (Y4)} in Proposition 2.9.
That is,
$$\aligned
&Y_{\omega_{i_k}}^+(p,q,a)\subset |Y_{\omega_{i_k}}^-(p,q,a)|,\,\,
Y_{\omega_{i_{k+1}}}^+(p,q,a)\subset
|Y_{\omega_{i_{k+1}}}^-(p,q,a)|\quad (1\le k\le t-1),\\
&Y_{\omega_{i_t}}^+(p,q,a)\subset |Y_{\omega_{i_t}}^-(p,q,a)|.
\endaligned
$$

At first, we consider the assumptions of the condition {\bf (Y4)}.
The top of ${\overline Y_{\omega_{i_k}}}$ (resp. ${\overline
Y_{\omega_{i_t}}}$) at $p$-th column from the right is
$\raisebox{-0.4\height}{\usebox{\tmpfiga}}$ and the top of
${\overline Y_{\omega_{i_{k+1}}}}$ (resp. ${\overline Y_{\la_n}}$)
at $q$-th column from the right is
$\raisebox{-0.4\height}{\usebox{\tmpfigb}}$ with $p>q$. Then these
assumptions can be represented in the tableau
$T_{\omega_{i_k}+\omega_{i_{k+1}}}$ and $T_{\omega_{i_t}+\la_n}$
as follows.

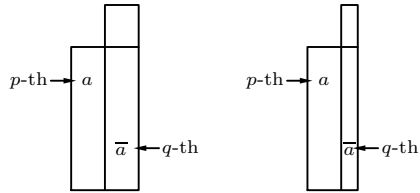
\begin{center}
\begin{texdraw}
\fontsize{7}{7}\selectfont \drawdim em \setunitscale 0.8

\bsegment \rlvec(4 0)\rlvec(0 11)\rlvec(-2 0)\rlvec(0
-2.5)\rlvec(-2 0)\rlvec(0 -8.5) \move(2 8.5)\rlvec(2 0) \move(2
8.5)\rlvec(0 -8.5)

\move(-1.2 6.5)\ravec(1.5 0)\htext(-2.5 6.5){$p$-th}\htext(1
6.5){$a$}

\move(5.2 2.5)\ravec(-1.5 0)\htext(6.5 2.5){$q$-th}\htext(3
2.5){$\overline{a}$} \esegment

\move(14 0)\bsegment \rlvec(3 0)\rlvec(0 11)\rlvec(-1 0)\rlvec(0
-2.5)\rlvec(-2 0)\rlvec(0 -8.5) \move(2 8.5)\rlvec(1 0)\move(2
8.5)\rlvec(0 -8.5)

\move(-1.2 6.5)\ravec(1.5 0)\htext(-2.5 6.5){$p$-th} \htext(1
6.5){$a$}

\move(4.2 2.5)\ravec(-1.5 0)\htext(5.5 2.5){$q$-th}\htext(2.5
2.5){$\overline{a}$} \esegment
\end{texdraw}
\end{center}

\savebox{\tmpfiga}{\begin{texdraw} \fontsize{7}{7}\selectfont
\drawdim em \setunitscale 0.8

\rlvec(2 0)\rlvec(0 11) \rlvec(-2 0)\rlvec(0 -11)\move(0
2)\rlvec(2 0)\move(0 9)\rlvec(2 0)

\move(0 9)\clvec (-0.5 9)(-1 7)(-1 6.5)

\htext(-2.5 5.5){$p\!\!-\!\!q\!\!-\!\!1\ge$}

\move(0 2)\clvec (-0.5 2)(-1 4)(-1 4.5)

\htext(1 1){$\overline{a}$}\htext(1 10){$a$}
\end{texdraw}}
\savebox{\tmpfigb}{\begin{texdraw} \fontsize{7}{7}\selectfont
\drawdim em \setunitscale 0.8

\rlvec(2 0)\rlvec(0 9) \rlvec(-2 0)\rlvec(0 -9)\move(0 2)\rlvec(2
0)\move(0 7)\rlvec(2 0)

\move(0 7)\clvec (-0.5 7)(-1 6)(-1 5.5)

\htext(-2.5 4.5){$p\!\!-\!\!q\!\!-\!\!2\ge$}

\move(0 2)\clvec (-0.5 2)(-1 3)(-1 3.5)

 \htext(1
1){$\overline{a\!\!+\!\!1}$}\htext(1 8){$a\!\!+\!\!1$}
\end{texdraw}}
\savebox{\tmpfigc}{\begin{texdraw} \fontsize{7}{7}\selectfont
\drawdim em \setunitscale 0.8

\rlvec(2 0)\rlvec(0 7) \rlvec(-2 0)\rlvec(0 -7)\move(0 2)\rlvec(2
0)\move(0 5)\rlvec(2 0)

\move(0 5)\clvec (-0.3 5)(-0.7 4.6)(-0.7 4.5)

\htext(-4 3.5){$p\!\!-\!\!q\!\!-\!\!(k\!\!+\!\!1)\ge$}

\move(0 2)\clvec (-0.3 2)(-0.7 2.4)(-0.7 2.5)

\htext(1 1){$\overline{a\!\!+\!\!k}$}\htext(1 6){$a\!\!+\!\!k$}
\end{texdraw}}
\savebox{\tmpfigd}{\begin{texdraw} \fontsize{7}{7}\selectfont
\drawdim em \setunitscale 0.8

\rlvec(2 0)\rlvec(0 4) \rlvec(-2 0)\rlvec(0 -4)\move(0 2)\rlvec(2
0)

\move(3 1)\ravec(-2 0)\move(3 3)\ravec(-2 0)

\htext(7
1){$\overline{a\!\!+\!\!(p\!\!-\!\!q\!\!-\!\!1)}$}\htext(7
3){$a\!\!+\!\!(p\!\!-\!\!q\!\!-\!\!1)$}
\end{texdraw}}

\noindent Now, consider the subtableau
$T^{(p,q)}_{\omega_{i_k}+\omega_{i_{k+1}}}$ and
$T^{(p,q)}_{\omega_{i_t}+\la_{n}}$ consisting of boxes lying
between $q$-th row and $p$-th row in
$T_{\omega_{i_k}+\omega_{i_{k+1}}}$ and $T_{\omega_{i_t}+\la_n}$.
Then the condition $Y_{\omega_{i_k}}^+(p,q,a)\subset
|Y_{\omega_{i_k}}^-(p,q,a)|$ (resp.
$Y_{\omega_{i_t}}^+(p,q,a)\subset |Y_{\omega_{i_t}}^-(p,q,a)|$),
$Y_{\omega_{i_{k+1}}}^+(p,q,a)\subset
|Y_{\omega_{i_{k+1}}}^-(p,q,a)|$ imply that the following can not
exist in $T_{\omega_{i_k}}\cap
T^{(p,q)}_{\omega_{i_k}+\omega_{i_{k+1}}}$,
$T_{\omega_{i_{k+1}}}\cap
T^{(p,q)}_{\omega_{i_k}+\omega_{i_{k+1}}}$ and
$T_{\omega_{i_t}}\cap T^{(p,q)}_{\omega_{i_t}+\la_{n}}$ as a
subtableau.
$$
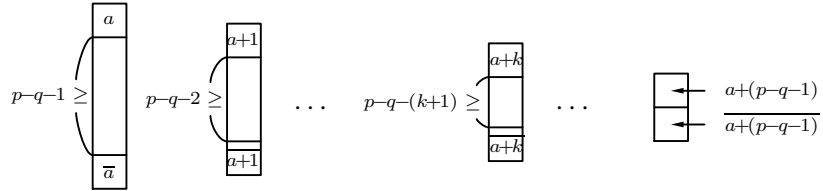
\raisebox{-0.4\height}{\usebox{\tmpfiga}}\, \hskip 7mm
\raisebox{-0.4\height}{\usebox{\tmpfigb}}\, \hskip 3mm \cdots
\hskip 12mm \raisebox{-0.4\height}{\usebox{\tmpfigc}}\,\hskip
3mm\cdots\hskip 8mm \raisebox{-0.4\height}{\usebox{\tmpfigd}}
$$
Therefore, if $T_{\omega_{i_k}}$ and $T_{\omega_{i_{k+1}}}$ has
$b$-box and $\bar{b}$-box between $p$-th column and $q$-th column,
then the number of rows of between $b$-box and $\bar{b}$-box in
$T_{\omega_{i_k}}$ or $T_{\omega_{i_{k+1}}}$ is larger than to
$p-q-(b-a+1)$. That is, if $b$-box and $\bar{b}$-box lie in $r$-th
and $s$-th rows of $T_{\omega_{i_k}}$ or $T_{\omega_{i_{k+1}}}$
from bottom to top with $r>s$, we have $r-s-1> p-q-(b-a+1)$ and so
\begin{equation}
(p-r)+(s-q)< b-a.
\end{equation}

\savebox{\tmpfiga}{\begin{texdraw} \fontsize{7}{7}\selectfont
\drawdim em \setunitscale 0.8

\rlvec(2 0)\rlvec(0 10) \rlvec(-2 0)\rlvec(0 -10)\move(0
2)\rlvec(2 0)\move(0 8)\rlvec(2 0)

\move(0 8)\clvec (-0.5 8)(-1 6.5)(-1 6)

\htext(-3 5){$p\!\!-\!\!q\!\!-\!\!(k\!\!+\!\!a)$}

\move(0 2)\clvec (-0.5 2)(-1 3.5)(-1 4)

\htext(1 1){$\overline{a\!\!+\!\!k}$}\htext(1 9){$a\!\!+\!\!k$}
\end{texdraw}}
\savebox{\tmpfigb}{\begin{texdraw} \fontsize{7}{7}\selectfont
\drawdim em \setunitscale 0.8

\rlvec(2 0)\rlvec(0 10) \rlvec(-2 0)\rlvec(0 -10)\move(0
2)\rlvec(2 0)\move(0 8)\rlvec(2 0)

\htext(1 1){$\overline{a}$}\htext(1 9){$a\!\!+\!\!k$}
\end{texdraw}}
\savebox{\tmpfigc}{\begin{texdraw} \fontsize{7}{7}\selectfont
\drawdim em \setunitscale 0.8

\rlvec(2 0)\rlvec(0 8) \rlvec(-2 0)\rlvec(0 -8)\move(0 2)\rlvec(2
0)\move(0 6)\rlvec(2 0)

\move(0 6)\clvec (-0.5 6)(-1 5.3)(-1 5)

\htext(-3 4){$p\!\!-\!\!q\!\!-\!\!(k\!\!+\!\!a\!\!+\!\!1)$}

\move(0 2)\clvec (-0.5 2)(-1 2.7)(-1 3)

\htext(1 1){$\overline{a\!\!+\!\!k\!\!+\!\!1}$}\htext(1
7){$a\!\!+\!\!k$}
\end{texdraw}}
\savebox{\tmpfigd}{\begin{texdraw} \fontsize{7}{7}\selectfont
\drawdim em \setunitscale 0.8

\rlvec(2 0)\rlvec(0 8) \rlvec(-2 0)\rlvec(0 -8)\move(0 2)\rlvec(2
0)\move(0 6)\rlvec(2 0)

\htext(1 1){$\overline{a\!\!+\!\!1}$}\htext(1 7){$a\!\!+\!\!k$}
\end{texdraw}}
\savebox{\tmpfige}{\begin{texdraw} \fontsize{7}{7}\selectfont
\drawdim em \setunitscale 0.8

\rlvec(2 0)\rlvec(0 4) \rlvec(-2 0)\rlvec(0 -4)\move(0 2)\rlvec(2
0)

\htext(1 1){$\overline{p\!\!-\!\!q}$}\htext(1 3){$a\!\!+\!\!k$}
\end{texdraw}}
\savebox{\tmpfigf}{\begin{texdraw} \fontsize{7}{7}\selectfont
\drawdim em \setunitscale 0.8

\rlvec(2 0)\rlvec(0 4) \rlvec(-2 0)\rlvec(0 -4)\move(0 2)\rlvec(2
0)

\move(3 1)\ravec(-2 0)

\htext(7
1){$\overline{p\!\!-\!\!q\!\!+\!\!1\!\!-\!\!(a\!\!+\!\!k)}$,}\htext(1
3){$a\!\!+\!\!k$}
\end{texdraw}}

\begin{remark} (a) In fact, we know that the following also can not exist in
$T_{\omega_{i_k}}\cap T^{(p,q)}_{\omega_{i_k}+\omega_{i_{k+1}}}$,
$T_{\omega_{i_{k+1}}}\cap
T^{(p,q)}_{\omega_{i_k}+\omega_{i_{k+1}}}$ and
$T_{\omega_{i_t}}\cap T^{(p,q)}_{\omega_{i_t}+\la_{n}}$ as a
subtableau.
$$
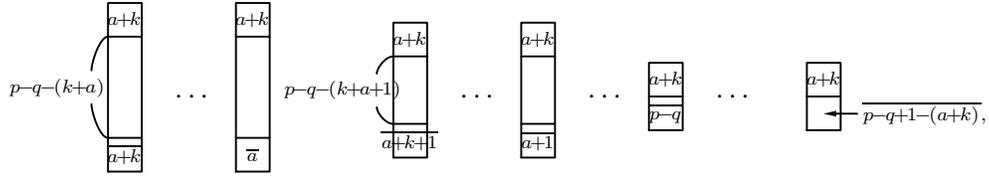
\raisebox{-0.4\height}{\usebox{\tmpfiga}}\,\hskip 3mm\cdots\hskip
3mm \raisebox{-0.4\height}{\usebox{\tmpfigb}}\, \hskip 9mm
\raisebox{-0.4\height}{\usebox{\tmpfigc}}\,\hskip 3mm\cdots\hskip
3mm \raisebox{-0.4\height}{\usebox{\tmpfigd}}\,\hskip
3mm\cdots\hskip 3mm
\raisebox{-0.4\height}{\usebox{\tmpfige}}\,\hskip 3mm\cdots\hskip
7mm \raisebox{-0.4\height}{\usebox{\tmpfigf}}
$$
where $0\le k\le p-q-a+1$. Therefore, if $T_{\omega_{i_k}}$ or
$T_{\omega_{i_{k+1}}}$ has $b$-box and $\bar{c}$-box between
$p$-th column and $q$-th column, then the number of rows of
between $b$-box and $\bar{c}$-box is larger than
$p-q-(\text{max}\,(b,c)-a+1)$. That is, if $b$-box and
$\bar{c}$-box lie in $r$-th and $s$-th rows of $T_{\omega_{i_k}}$
from bottom to top with $r>s$, we have $r-s-1>
p-q-(\text{max}\,(b,c)-a+1)$ and so
\begin{equation}
(p-r)+(s-q)< \text{max}\,(b,c)-a.
\end{equation}

(b) This condition is very similar to the two column condition
(2CC) for the tableau of Kashiwara and Nakashima in the
$(a,b)$-configuration. We will call it the second two column
condition and denote by (2CC-2).

(c) The first two column condition (2CC-1) covers a lot of parts
of this condition (2CC-2).
\end{remark}

\savebox{\tmpfiga}{\begin{texdraw} \fontsize{7}{7}\selectfont
\drawdim em \setunitscale 0.8

\rlvec(0 2)\rlvec(2 0)\rlvec(0 -2)\rlvec(-2 0)

\htext(1 1){$x_1$}
\end{texdraw}}
\savebox{\tmpfigb}{\begin{texdraw} \fontsize{7}{7}\selectfont
\drawdim em \setunitscale 0.8

\rlvec(0 2)\rlvec(2 0)\rlvec(0 -2)\rlvec(-2 0)

\htext(1 1){$x_t$}
\end{texdraw}}
\savebox{\tmpfigc}{\begin{texdraw} \fontsize{7}{7}\selectfont
\drawdim em \setunitscale 0.8

\rlvec(0 2)\rlvec(2 0)\rlvec(0 -2)\rlvec(-2 0)

\htext(1 1){$x_k$}
\end{texdraw}}

\vskip 3mm Now, we are ready to give an another explicit
realization of crystal graph $B(\la)$ over $B_n$.

\vskip 3mm For a dominant integral weight $\la$, we define
$T(\la)$ for $\frak g=B_n$ by the set of tableaux of shape
$GRY^{-1}(\la)$ with entries $\{i,\overline{i}\,|1\le i\le
n\}\cup\{0\}$ such that
\begin{enumerate}
\item the entries of $T$ weakly increase along the rows, but the
element $0$ cannot appear more than once,
\item the entries of $T$ strictly increase down the columns, but
the element $0$ can appear more than once,
\item for a half column $C$ of $T$, $i$ and $\overline{i}$ can not
appear at the same time,
\item for each column $C$ of $T$, (1CC) holds,
\item for each pair of adjacent columns $C$, $C'$ of $T$, (2CC-1)
and (2CC-2) holds,
\end{enumerate}
and we identify a tableau $T$ of shape $GRY^{-1}(\la)$ with the
vector
$$\raisebox{-0.3\height}{\usebox{\tmpfiga}}\otimes\cdots\otimes
\raisebox{-0.3\height}{\usebox{\tmpfigb}}\in {\bf B}^{\otimes
t}\,\,\,\,\text{or}\,\,\,\,
v_{sp}\otimes\raisebox{-0.3\height}{\usebox{\tmpfiga}}\otimes\cdots\otimes
\raisebox{-0.3\height}{\usebox{\tmpfigb}} \in B_{sp}\otimes{\bf
B}^{\otimes t},$$ where $x_i$ is an entry of $T$ by reading from
top to bottom and from right to left. Then we have:

\savebox{\tmpfiga}{\begin{texdraw} \fontsize{7}{7}\selectfont
\drawdim em \setunitscale 0.85

\lvec(2 0) \rlvec(0 2)\rlvec(-2 0)\rlvec(0 -2)

\htext(1 1){$i$}
\end{texdraw}}
\savebox{\tmpfigb}{\begin{texdraw} \fontsize{7}{7}\selectfont
\drawdim em \setunitscale 0.85

\lvec(2 0) \rlvec(0 2)\rlvec(-2 0)\rlvec(0 -2)

\htext(1 1){$i\!\!+\!\!1$}
\end{texdraw}}
\savebox{\tmpfigc}{\begin{texdraw} \fontsize{7}{7}\selectfont
\drawdim em \setunitscale 0.85

\lvec(2 0) \rlvec(0 2)\rlvec(-2 0)\rlvec(0 -2)

\htext(1 1){$\overline{i\!\!+\!\!1}$}
\end{texdraw}}
\savebox{\tmpfigd}{\begin{texdraw} \fontsize{7}{7}\selectfont
\drawdim em \setunitscale 0.85

\lvec(2 0) \rlvec(0 2)\rlvec(-2 0)\rlvec(0 -2)

\htext(1 1){$\overline{i}$}
\end{texdraw}}

\vskip 3mm
\begin{thm} \label{thm:B_n}
For a dominant integral weight $\lambda$, there exists a crystal
isomorphism for $U_q(B_n)$-modules $\varphi:Y(\la)\rightarrow
T(Y)$sending $Y$ to $T_Y$.
\end{thm}

\begin{proof} It is clear that $\varphi$ is a bijection.
So, it suffices to show that it is a crystal morphism. Suppose
that $\fit Y$ ($1\le i<n$) is obtained by stacking the block
$\raisebox{-0.3\height}{\usebox{\tmpfiga}}$ on top of some column
$y_a$ of $Y$. Then these columns of $Y$ and $\fit Y$ correspond to
the boxes  $\raisebox{-0.3\height}{\usebox{\tmpfiga}}$ and
$\raisebox{-0.3\height}{\usebox{\tmpfigb}}$, or
$\raisebox{-0.3\height}{\usebox{\tmpfigc}}$ and
$\raisebox{-0.3\height}{\usebox{\tmpfigd}}$ of $T_Y$ and $T_{\fit
Y}$. Moreover, it is easy to see that the $i$-signatures of each
column $y_k$ of $Y$ and a corresponded box $\varphi(y_k)$ of $T_Y$
are the same. Therefore, by the tensor product rule of Kashiwara
operators, it  is easy to see that $\fit T_Y$ is obtained by
acting on $\varphi(y_a)$ of $T_Y$ and so $\varphi(\fit Y)= \fit
\varphi(Y).$ For $i=n$, it is proved by the same method. Moreover,
for $\eit$, it is also proved by the similar argument.
\end{proof}


\vskip 3mm By the same algorithm, we can obtain new tableau
realization of crystal bases of irreducible highest weight modules
over other classical Lie algebras. The basic data in {\bf Step 1}
and {\bf Step 2} over other classical Lie algebras are as follows:

\savebox{\tmpfiga}{\begin{texdraw} \fontsize{7}{7}\selectfont
\drawdim em \setunitscale 1.7 \move(-1 0)\lvec(-1 1)\lvec(0
1)\lvec(-1 0)\lvec(0 0)\lvec(0 1) \htext(-0.7 0.75){$0$}
\end{texdraw}}
\savebox{\tmpfigb}{\begin{texdraw} \fontsize{7}{7}\selectfont
\drawdim em \setunitscale 1.7 \move(-1 0)\lvec(-1 1)\lvec(0
1)\lvec(-1 0)\lvec(0 0)\lvec(0 1) \htext(-0.3 0.3){$0$}
\end{texdraw}}
\savebox{\tmpfigc}{\begin{texdraw} \fontsize{7}{7}\selectfont
\drawdim em \setunitscale 1.7 \move(0 0)\lvec(-1 0)\lvec(-1
1)\lvec(0 1)\lvec(0 0) \htext(-0.5 0.5){$1$}
\end{texdraw}}
\savebox{\tmpfigd}{\begin{texdraw} \fontsize{7}{7}\selectfont
\drawdim em \setunitscale 1.7 \move(-1 0)\lvec(-1 1)\lvec(0
1)\lvec(-1 0)\lvec(0 0)\lvec(0 1) \htext(-0.7 0.75){$1$}
\end{texdraw}}
\savebox{\tmpfige}{\begin{texdraw} \fontsize{7}{7}\selectfont
\drawdim em \setunitscale 1.7 \move(-1 0)\lvec(-1 1)\lvec(0
1)\lvec(-1 0)\lvec(0 0)\lvec(0 1) \htext(-0.3 0.3){$1$}
\end{texdraw}}
\savebox{\tmpfigf}{\begin{texdraw} \fontsize{7}{7}\selectfont
\drawdim em \setunitscale 1.7 \move(0 0)\lvec(-1 0)\lvec(-1
1)\lvec(0 1)\lvec(0 0) \htext(-0.5 0.5){$\bar{1}$}
\end{texdraw}}
\savebox{\tmpfigg}{\begin{texdraw} \fontsize{7}{7}\selectfont
\drawdim em \setunitscale 1.7 \lpatt(0.1 0.15) \move(0 1)\lvec(0
2)\lvec(-1 2)\lvec(-1 1) \move(-1 1)\lvec(0 2) \lpatt(1 0) \move(0
0)\lvec(-1 0)\lvec(-1 1)\lvec(0 1)\lvec(0 0) \htext(-0.5 0.5){$2$}
\end{texdraw}}
\savebox{\tmpfigh}{\begin{texdraw} \fontsize{7}{7}\selectfont
\drawdim em \setunitscale 1.7 \move(-1 0)\lvec(-1 2)\lvec(0
2)\lvec(0 0)\lvec(-1 0) \move(-1 0)\lvec(0 1)\lvec(-1 1)
\htext(-0.5 1.5){$2$}
\end{texdraw}}
\savebox{\tmpfigi}{\begin{texdraw} \fontsize{7}{7}\selectfont
\drawdim em \setunitscale 1.7 \move(0 0)\lvec(-1 0)\lvec(-1
1)\lvec(0 1)\lvec(0 0) \htext(-0.5 0.5){$2$}
\end{texdraw}}
\savebox{\tmpfigj}{\begin{texdraw} \fontsize{7}{7}\selectfont
\drawdim em \setunitscale 1.7 \lpatt(0.1 0.15) \move(-1 1)\lvec(-1
2)\lvec(0 2)\lvec(0 1) \lpatt(1 0) \move(0 0)\lvec(-1 0)\lvec(-1
1)\lvec(0 1)\lvec(0 0) \htext(-0.5 0.5){$3$} \htext(-0.5 1.5){$2$}
\end{texdraw}}
\savebox{\tmpfigk}{\begin{texdraw} \fontsize{7}{7}\selectfont
\drawdim em \setunitscale 1.7 \move(0 0)\lvec(-1 0)\lvec(-1
1)\lvec(0 1)\lvec(0 0) \htext(-0.5 0.5){$\bar{2}$}
\end{texdraw}}
\savebox{\tmpfigl}{\begin{texdraw} \fontsize{7}{7}\selectfont
\drawdim em \setunitscale 1.7 \move(0 0)\lvec(-1 0) \move(0
1)\lvec(-1 1)\move(0 0)\lvec(0 2)\lvec(-1 2)\lvec(-1 0)
\htext(-0.5 0.5){$j\!\!-\!\!2$} \htext(-0.5 1.5){$j\!\!-\!\!1$}
\end{texdraw}}
\savebox{\tmpfigm}{\begin{texdraw} \fontsize{7}{7}\selectfont
\drawdim em \setunitscale 1.7 \move(0 0)\lvec(-1 0)\lvec(-1
1)\lvec(0 1)\lvec(0 0) \htext(-0.5 0.5){$j$}
\end{texdraw}}
\savebox{\tmpfign}{\begin{texdraw} \fontsize{7}{7}\selectfont
\drawdim em \setunitscale 1.7 \move(0 0)\lvec(-1 0) \move(0
1)\lvec(-1 1)\move(0 0)\lvec(0 2)\lvec(-1 2)\lvec(-1 0)
\htext(-0.5 1.5){$j$} \htext(-0.5 0.5){$j\!\!+\!\!1$}
\end{texdraw}}
\savebox{\tmpfigo}{\begin{texdraw} \fontsize{7}{7}\selectfont
\drawdim em \setunitscale 1.7 \move(0 0)\lvec(-1 0)\lvec(-1
1)\lvec(0 1)\lvec(0 0) \htext(-0.5 0.5){$\bar{j}$}
\end{texdraw}}
\savebox{\tmpfigp}{\begin{texdraw} \fontsize{7}{7}\selectfont
\drawdim em \setunitscale 1.7 \move(0 0)\lvec(-1 0)\move(0
1)\lvec(-1 1) \move(0 0)\lvec(0 2)\lvec(-1 2)\lvec(-1 0)
\htext(-0.5 0.5){$n\!\!-\!\!1$} \htext(-0.5 1.5){$n$}
\end{texdraw}}
\savebox{\tmpfigq}{\begin{texdraw} \fontsize{7}{7}\selectfont
\drawdim em \setunitscale 1.7 \move(0 0)\lvec(-1 0)\lvec(-1
1)\lvec(0 1)\lvec(0 0) \htext(-0.5 0.5){$0$}
\end{texdraw}}
\savebox{\tmpfigr}{\begin{texdraw} \fontsize{7}{7}\selectfont
\drawdim em \setunitscale 1.7 \move(0 0)\lvec(-1 0)\lvec(-1
1)\lvec(0 1)\lvec(0 0) \move(-1 0)\lvec(-0.5 0.5)\move(-0.05
0.95)\lvec(0 1) \htext(-0.5 0.75){$n\!\!-\!\!1$}
\end{texdraw}}
\savebox{\tmpfigs}{\begin{texdraw} \fontsize{7}{7}\selectfont
\drawdim em \setunitscale 1.7 \move(0 0)\lvec(-1 0)\lvec(-1
1)\lvec(0 1)\lvec(0 0) \move(-1 0)\lvec(-0.95 0.05)\move(-0.5
0.5)\lvec(0 1) \htext(-0.5 0.25){$n\!\!-\!\!1$}
\end{texdraw}}
\savebox{\tmpfigt}{\begin{texdraw} \fontsize{7}{7}\selectfont
\drawdim em \setunitscale 1.7 \move(0 0)\lvec(-1 0)\lvec(-1
1)\lvec(0 1)\lvec(0 0) \htext(-0.5 0.5){$n$}
\end{texdraw}}
\savebox{\tmpfigu}{\begin{texdraw} \fontsize{7}{7}\selectfont
\drawdim em \setunitscale 1.7 \move(-1 0)\lvec(-1 1)\lvec(0
1)\lvec(-1 0)\lvec(0 0)\lvec(0 1) \htext(-0.7 0.75){$n$}
\end{texdraw}}
\savebox{\tmpfigv}{\begin{texdraw} \fontsize{7}{7}\selectfont
\drawdim em \setunitscale 1.7 \move(-1 0)\lvec(-1 1)\lvec(0
1)\lvec(-1 0)\lvec(0 0)\lvec(0 1) \htext(-0.3 0.3){$n$}
\end{texdraw}}
\savebox{\tmpfigw}{\begin{texdraw} \fontsize{7}{7}\selectfont
\drawdim em \setunitscale 1.7 \move(0 0)\lvec(-1 0)\lvec(-1
1)\lvec(0 1)\lvec(0 0) \htext(-0.5 0.5){$\bar{n}$}
\end{texdraw}}
\savebox{\tmpfigx}{\begin{texdraw} \fontsize{7}{7}\selectfont
\drawdim em \setunitscale 1.7 \move(0 0)\lvec(0 1)\lvec(-1
1)\lvec(-1 0) \move(0 0)\lvec(-1 0)\lvec(0 1) \htext(-0.3
0.3){$1$} \htext(-0.7 0.75){$0$}
\end{texdraw}}
\savebox{\tmpfigy}{\begin{texdraw} \fontsize{7}{7}\selectfont
\drawdim em \setunitscale 1.7 \move(0 0)\lvec(0 1)\lvec(-1
1)\lvec(-1 0) \move(0 0)\lvec(-1 0)\lvec(0 1) \htext(-0.3
0.3){$0$} \htext(-0.7 0.75){$1$}
\end{texdraw}}
\savebox{\tmpfigz}{\begin{texdraw} \fontsize{7}{7}\selectfont
\drawdim em \setunitscale 1.7 \lpatt(0.1 0.15) \move(-1 1)\lvec(-1
2)\lvec(0 2)\lvec(0 1) \lpatt(1 0) \move(0 0)\lvec(-1 0)\lvec(-1
1)\lvec(0 1)\lvec(0 0) \htext(-0.5 0.5){$n\!\!-\!\!2$} \htext(-0.5
1.5){$n\!\!-\!\!1$}
\end{texdraw}}
\savebox{\tmpfigaa}{\begin{texdraw} \fontsize{7}{7}\selectfont
\drawdim em \setunitscale 1.7 \move(0 0)\lvec(-1 0)\lvec(-1
1)\lvec(0 1)\lvec(0 0) \htext(-0.5 0.5){$n\!\!-\!\!1$}
\end{texdraw}}
\savebox{\tmpfigad}{\begin{texdraw} \fontsize{7}{7}\selectfont
\drawdim em \setunitscale 1.7 \move(0 0)\lvec(-1 0)\lvec(-1
1)\lvec(0 1)\lvec(0 0) \htext(-0.5 0.5){$\overline{n\!\!-\!\!1}$}
\end{texdraw}}
\savebox{\tmpfigab}{\begin{texdraw} \fontsize{7}{7}\selectfont
\drawdim em \setunitscale 1.7 \move(0 0)\lvec(-1 0) \move(0
1)\lvec(-1 1)\lvec(-0.5 1.5) \move(-0.05 1.95)\lvec(0 2) \move(0
0)\lvec(0 2)\lvec(-1 2)\lvec(-1 0) \htext(-0.5 0.5){$n\!\!-\!\!2$}
\htext(-0.3 1.3){$n$} \htext(-0.5 1.75){$n\!\!-\!\!1$}
\end{texdraw}}
\savebox{\tmpfigac}{\begin{texdraw} \fontsize{7}{7}\selectfont
\drawdim em \setunitscale 1.7 \move(0 0)\lvec(-1 0) \move(0
1)\lvec(-1 1)\lvec(-0.95 1.05) \move(-0.5 1.5)\lvec(0 2) \move(0
0)\lvec(0 2)\lvec(-1 2)\lvec(-1 0) \htext(-0.5 0.5){$n\!\!-\!\!2$}
\htext(-0.7 1.75){$n$} \htext(-0.5 1.25){$n\!\!-\!\!1$}
\end{texdraw}}
\savebox{\tmpfigag}{\begin{texdraw} \fontsize{7}{7}\selectfont
\drawdim em \setunitscale 1.7 \move(0 0)\lvec(-1 0)\lvec(-1
1)\lvec(0 1)\lvec(0 0) \htext(-0.5 0.5){$\overline{n\!\!-\!\!2}$}
\end{texdraw}}
\savebox{\tmpfigae}{\begin{texdraw} \fontsize{7}{7}\selectfont
\drawdim em \setunitscale 1.7 \move(0 0)\lvec(-1 0)\lvec(-0.5 0.5)
\move(-0.05 0.95)\lvec(0 1) \move(0 1)\lvec(-1 1) \move(0
0)\lvec(0 2)\lvec(-1 2)\lvec(-1 0) \htext(-0.5 1.5){$n\!\!-\!\!2$}
\htext(-0.5 0.75){$n\!\!-\!\!1$} \htext(-0.3 0.3){$n$}
\end{texdraw}}
\savebox{\tmpfigaf}{\begin{texdraw} \fontsize{7}{7}\selectfont
\drawdim em \setunitscale 1.7 \move(0 0)\lvec(-1 0)\lvec(-0.95
0.05) \move(-0.5 0.5)\lvec(0 1) \move(0 1)\lvec(-1 1) \move(0
0)\lvec(0 2)\lvec(-1 2)\lvec(-1 0) \htext(-0.5 1.5){$n\!\!-\!\!2$}
\htext(-0.5 0.25){$n\!\!-\!\!1$} \htext(-0.7 0.75){$n$}
\end{texdraw}}
\savebox{\tmpfigah}{\begin{texdraw} \fontsize{7}{7}\selectfont
\drawdim em \setunitscale 1.7 \move(0 0)\lvec(0 0.5)\lvec(-1
0.5)\lvec(-1 0)\lvec(0 0) \htext(-0.5 0.25){$0$}
\end{texdraw}}
\savebox{\tmpfigai}{\begin{texdraw} \fontsize{7}{7}\selectfont
\drawdim em \setunitscale 1.7 \move(0 0)\lvec(0 1)\lvec(-1
1)\lvec(-1 0)\lvec(0 0) \move(0 0.5)\lvec(-1 0.5) \htext(-0.5
0.25){$0$} \htext(-0.5 0.75){$0$}
\end{texdraw}}
\savebox{\tmpfigaj}{\begin{texdraw} \fontsize{7}{7}\selectfont
\drawdim em \setunitscale 1.7 \move(0 0)\lvec(0 2)\lvec(-1
2)\lvec(-1 0) \move(0 0)\lvec(-1 0)\move(0 1)\lvec(-1 1)
\htext(-0.5 0.5){$3$} \htext(-0.5 1.5){$2$}
\end{texdraw}}
\savebox{\tmpfigak}{\begin{texdraw} \fontsize{7}{7}\selectfont
\drawdim em \setunitscale 1.7 \move(0 0)\lvec(0 2)\lvec(-1
2)\lvec(-1 0) \move(0 0)\lvec(-1 0)\move(0 1)\lvec(-1 1)
\htext(-0.5 0.5){$n\!\!-\!\!2$} \htext(-0.5 1.5){$n\!\!-\!\!1$}
\end{texdraw}}
\savebox{\tmpfigal}{\begin{texdraw} \fontsize{7}{7}\selectfont
\drawdim em \setunitscale 1.7 \move(0 0)\lvec(0 0.5)\lvec(-1
0.5)\lvec(-1 0)\lvec(0 0) \htext(-0.5 0.26){$n$}
\end{texdraw}}
\savebox{\tmpfigam}{\begin{texdraw} \fontsize{7}{7}\selectfont
\drawdim em \setunitscale 1.7 \move(0 0)\lvec(0 1)\lvec(-1
1)\lvec(-1 0)\lvec(0 0) \move(0 0.5)\lvec(-1 0.5) \htext(-0.5
0.26){$n$} \htext(-0.5 0.76){$n$}
\end{texdraw}}
\savebox{\tmpfigan}{\begin{texdraw} \fontsize{7}{7}\selectfont
\drawdim em \setunitscale 1.7 \move(0 0)\lvec(-1 0)\lvec(-1
1)\lvec(0 1)\lvec(0 0) \htext(-0.5 0.5){$\bar{3}$}
\end{texdraw}}
\savebox{\tmpfigao}{\begin{texdraw} \fontsize{7}{7}\selectfont
\drawdim em \setunitscale 1.7 \move(0 0)\lvec(-1 0)\lvec(-1
1)\lvec(0 1)\lvec(0 0) \htext(-0.5 0.5){$j\!\!+\!\!1$}
\end{texdraw}}
\savebox{\tmpfigap}{\begin{texdraw} \fontsize{7}{7}\selectfont
\drawdim em \setunitscale 1.7 \move(0 0)\lvec(0 1)\lvec(-1
1)\lvec(-1 0) \move(0 0)\lvec(-1 0)\lvec(0 1) \move(-1 1)\lvec(-1
2)\lvec(0 2)\lvec(0 1) \htext(-0.3 0.3){$1$} \htext(-0.7
0.75){$0$} \htext(-0.5 1.5){$2$}
\end{texdraw}}
\savebox{\tmpfigaq}{\begin{texdraw} \fontsize{7}{7}\selectfont
\drawdim em \setunitscale 1.7 \move(0 0)\lvec(0 1)\lvec(-1
1)\lvec(-1 0) \move(0 0)\lvec(-1 0)\lvec(0 1) \move(-1 1)\lvec(-1
2)\lvec(0 2)\lvec(0 1) \htext(-0.3 0.3){$0$} \htext(-0.7
0.75){$1$} \htext(-0.5 1.5){$2$}
\end{texdraw}}
\savebox{\tmpfigar}{\begin{texdraw} \fontsize{7}{7}\selectfont
\drawdim em \setunitscale 1.7 \move(0 0)\lvec(-1 0)\lvec(-1
1)\lvec(0 1)\lvec(0 0) \htext(-0.5 0.5){$3$}
\end{texdraw}}
\savebox{\tmpfigas}{\begin{texdraw} \fontsize{7}{7}\selectfont
\drawdim em \setunitscale 1.7 \move(0 0)\lvec(-1 0)\lvec(-1
1)\lvec(0 1)\lvec(0 0) \htext(-0.5 0.5){$\bar{n}$}
\end{texdraw}}
\savebox{\tmpfigat}{\begin{texdraw} \fontsize{7}{7}\selectfont
\drawdim em \setunitscale 1.7 \move(0 0)\lvec(-1 0)\lvec(-1
1)\lvec(0 1)\lvec(0 0) \move(-1 0)\lvec(-0.5 0.5)\move(-0.05
0.95)\lvec(0 1) \htext(-0.5 0.75){$n\!\!-\!\!1$} \htext(-0.3
0.3){$n$}
\end{texdraw}}
\savebox{\tmpfigau}{\begin{texdraw} \fontsize{7}{7}\selectfont
\drawdim em \setunitscale 1.7 \move(0 0)\lvec(-1 0)\lvec(-1
1)\lvec(0 1)\lvec(0 0) \move(-1 0)\lvec(-0.95 0.05)\move(-0.5
0.5)\lvec(0 1) \htext(-0.5 0.25){$n\!\!-\!\!1$} \htext(-0.7
0.75){$n$}
\end{texdraw}}
\savebox{\tmpfigav}{\begin{texdraw} \fontsize{7}{7}\selectfont
\drawdim em \setunitscale 1.7 \move(0 0)\lvec(-1 0)\lvec(-1
2)\lvec(0 2)\lvec(0 0) \move(-1 0.5)\lvec(0 0.5) \move(-1
1)\lvec(0 1) \htext(-0.5 0.26){$0$} \htext(-0.5 0.76){$0$}
\htext(-0.5 1.5){$1$}
\end{texdraw}}
\savebox{\tmpfigaw}{\begin{texdraw} \fontsize{7}{7}\selectfont
\drawdim em \setunitscale 1.7 \move(0 0)\lvec(-1 0)\lvec(-1
2)\lvec(0 2)\lvec(0 0) \move(-1 0.5)\lvec(0 0.5) \move(-1
1)\lvec(0 1) \htext(-0.5 0.26){$n$} \htext(-0.5 0.76){$n$}
\htext(-0.5 1.5){$n\!\!-\!\!1$}
\end{texdraw}}
\savebox{\tmpfigax}{\begin{texdraw} \fontsize{7}{7}\selectfont
\drawdim em \setunitscale 1.7 \move(0 0)\lvec(-1 0)\lvec(-1
1)\lvec(0 1)\lvec(0 0) \htext(-0.5 0.5){$i$}
\end{texdraw}}
\savebox{\tmpfigay}{\begin{texdraw} \fontsize{7}{7}\selectfont
\drawdim em \setunitscale 1.7 \move(0 0)\lvec(-1 0)\lvec(-1
1)\lvec(0 1)\lvec(0 0) \htext(-0.5 0.5){$i\!\!+\!\!1$}
\end{texdraw}}
\savebox{\tmpfigaz}{\begin{texdraw} \fontsize{7}{7}\selectfont
\drawdim em \setunitscale 1.7 \move(0 0)\lvec(-1 0)\lvec(-1
1)\lvec(0 1)\lvec(0 0) \htext(-0.5 0.5){$\overline{i\!\!+\!\!1}$}
\end{texdraw}}
\savebox{\tmpfigba}{\begin{texdraw} \fontsize{7}{7}\selectfont
\drawdim em \setunitscale 1.7 \move(0 0)\lvec(-1 0)\lvec(-1
1)\lvec(0 1)\lvec(0 0) \htext(-0.5 0.5){$\bar{i}$}
\end{texdraw}}
\savebox{\tmpfigbb}{\begin{texdraw} \fontsize{7}{7}\selectfont
\drawdim em \setunitscale 1.7 \lpatt(0.1 0.15) \move(-1 1)\lvec(-1
2)\lvec(0 2)\lvec(0 1) \lpatt(1 0) \move(0 0)\lvec(-1 0)\lvec(-1
1)\lvec(0 1)\lvec(0 0) \htext(-0.5 1.5){$i$} \htext(-0.5
0.5){$i\!\!-\!\!1$}
\end{texdraw}}
\savebox{\tmpfigbc}{\begin{texdraw} \fontsize{7}{7}\selectfont
\drawdim em \setunitscale 1.7 \move(0 0)\lvec(0 2)\lvec(-1
2)\lvec(-1 0) \move(0 0)\lvec(-1 0)\move(0 1)\lvec(-1 1)
\htext(-0.5 1.5){$i$} \htext(-0.5 0.5){$i\!\!-\!\!1$}
\end{texdraw}}
\savebox{\tmpfigbd}{\begin{texdraw} \fontsize{7}{7}\selectfont
\drawdim em \setunitscale 1.7 \move(0 0)\lvec(0 2)\lvec(-1
2)\lvec(-1 0) \move(0 0)\lvec(-1 0)\move(0 1)\lvec(-1 1)
\htext(-0.5 1.5){$i$} \htext(-0.5 0.5){$i\!\!+\!\!1$}
\end{texdraw}}
\savebox{\tmpfigbe}{\begin{texdraw} \fontsize{7}{7}\selectfont
\drawdim em \setunitscale 1.7 \lpatt(0.1 0.15) \move(-1 1)\lvec(-1
2)\lvec(0 2)\lvec(0 1) \lpatt(1 0) \move(0 0)\lvec(-1 0)\lvec(-1
1)\lvec(0 1)\lvec(0 0) \htext(-0.5 1.5){$i$} \htext(-0.5
0.5){$i\!\!+\!\!1$}
\end{texdraw}}
\savebox{\tmpfigbf}{\begin{texdraw} \fontsize{7}{7}\selectfont
\drawdim em \setunitscale 1.7 \move(0 0)\lvec(0 2)\lvec(-1
2)\lvec(-1 0)\lvec(0 0) \move(0 1)\lvec(-1 1)\lvec(0 2)
\htext(-0.3 1.3){$0$} \htext(-0.5 0.5){$2$}
\end{texdraw}}
\savebox{\tmpfigbg}{\begin{texdraw} \fontsize{7}{7}\selectfont
\drawdim em \setunitscale 1.7 \move(0 0)\lvec(0 2)\lvec(-1
2)\lvec(-1 0)\lvec(0 0) \move(0 1)\lvec(-1 1)\lvec(0 2)
\htext(-0.7 1.75){$0$} \htext(-0.5 0.5){$2$}
\end{texdraw}}
%
%
\vskip 3mm

\noindent{\bf Step 1.}

\begin{multicols}{2}[\text{1)  $A_n$ ($n\geq1$)}]
\hskip 5mm {\usebox{\tmpfigm}} \araise{$\longmapsto$}
{\usebox{\tmpfigao}}
\quad \araise{($0\leq j\leq n$), }\\[4pt]
\end{multicols}
\vskip 3mm
\begin{multicols}{2}[\text{2)  $C_n$ ($n\geq3$)}]
\hskip 5mm
\usebox{\tmpfigb} \araise{ or } \usebox{\tmpfiga}
\araise{$\longmapsto$} \usebox{\tmpfigc}\\[4pt]
\vskip 2mm \hskip 5mm
\usebox{\tmpfigx} \araise{ or } \usebox{\tmpfigy}
\araise{$\longmapsto$} \usebox{\tmpfigi}\\[4pt]
\vskip 2mm \hskip 5mm
\usebox{\tmpfigap} \braise{ or } \usebox{\tmpfigaq}
\braise{$\longmapsto$} \craise{\usebox{\tmpfigar}}\\[4pt]
\vskip 2mm \hskip 5mm
\usebox{\tmpfigl}
\braise{$\longmapsto$} \craise{\usebox{\tmpfigm}}
\quad \draise{($4\leq j \leq n$)}\\[4pt]
\vskip 2mm
\usebox{\tmpfigp}
\braise{$\longmapsto$} \craise{\usebox{\tmpfigw}}\\[4pt]
\vskip 4mm
\usebox{\tmpfign}
\braise{$\longmapsto$} \craise{\usebox{\tmpfigo}}
\quad \draise{($2\leq j \leq n-1$)}\\[4pt]
\vskip 4mm
\usebox{\tmpfigd} \araise{ or } \usebox{\tmpfige}
\araise{$\longmapsto$} \usebox{\tmpfigf}
\end{multicols}

\vskip 3mm
\begin{multicols}{2}[\text{3)  $D_n$ ($n\geq4$)}]
\hskip 5mm
\usebox{\tmpfigb} \araise{ or } \usebox{\tmpfiga}
\araise{$\longmapsto$} \usebox{\tmpfigc}\\[4pt]

\vskip 2mm \hskip 5mm
\usebox{\tmpfigx} \araise{ or } \usebox{\tmpfigy}
\araise{$\longmapsto$} \usebox{\tmpfigi}\\[4pt]

\vskip 2mm \hskip 5mm
\usebox{\tmpfigap} \braise{ or } \usebox{\tmpfigaq}
\braise{$\longmapsto$} \craise{\usebox{\tmpfigar}}\\[4pt]

\vskip 2mm \hskip 5mm
\usebox{\tmpfigl}
\braise{$\longmapsto$} \craise{\usebox{\tmpfigm}}
\quad \draise{($4\leq j \leq n-1$)}\\[4pt]

\vskip 2mm \hskip 5mm
\usebox{\tmpfigs} \araise{ or } \usebox{\tmpfigr}
\araise{$\longmapsto$} \usebox{\tmpfigt}\\[4pt]

\vskip 2mm
\usebox{\tmpfigv} \araise{ or } \usebox{\tmpfigu}
\araise{$\longmapsto$} \usebox{\tmpfigas}\\[4pt]

\vskip 2mm
\usebox{\tmpfigat} \araise{ or } \usebox{\tmpfigau}
\araise{$\longmapsto$} \usebox{\tmpfigad}\\[4pt]

\vskip 2mm
\usebox{\tmpfigae} \braise{ or } \usebox{\tmpfigaf}
\braise{$\longmapsto$} \craise{\usebox{\tmpfigag}}\\[4pt]

\vskip 2mm
\usebox{\tmpfign}
\braise{$\longmapsto$} \craise{\usebox{\tmpfigo}}
\quad \draise{($2\leq j \leq n-3$)}\\[4pt]

\vskip 2mm
\usebox{\tmpfigd} \araise{ or } \usebox{\tmpfige}
\araise{$\longmapsto$} \usebox{\tmpfigf}
\end{multicols}


%
\savebox{\tmpfiga}{\begin{texdraw} \fontsize{7}{7}\selectfont
\drawdim em \setunitscale 0.85

\rlvec(8 0)\rlvec(0 4)\move(2 0)\rlvec(2 2)\rlvec(0 -2) \move(6
2)\rlvec(-1.3 -1.3)\move(6 2)\rlvec(0 -2)\rlvec(2 2)\move(8
2)\rlvec(-6 0)\rlvec(0 -2)

\htext(3.5 0.5){$n$} \htext(5 0.5){$n\!\!-\!\!1$}\htext(7.5
0.5){$n$}\htext(1 0.5){$\cdots$}
\end{texdraw}}
\savebox{\tmpfigb}{\begin{texdraw} \fontsize{7}{7}\selectfont
\drawdim em \setunitscale 0.85

\rlvec(8 0)\rlvec(0 4)\move(4 0)\rlvec(0 2)\rlvec(-1.3
-1.3)\move(4 0)\rlvec(2 2)\rlvec(0 -2) \move(8 2) \rlvec(-1.3
-1.3) \move(8 2)\rlvec(-6 0)\rlvec(0 -2)

\htext(3 0.5){$n\!\!-\!\!1$} \htext(5.5 0.5){$n$}\htext(7
0.5){$n\!\!-\!\!1$}\htext(1 0.5){$\cdots$}
\end{texdraw}}
\savebox{\tmpfigc}{\begin{texdraw} \fontsize{7}{7}\selectfont
\drawdim em \setunitscale 0.8

\rlvec(1 0)\rlvec(0 9)\rlvec(-1 0)\rlvec(0 -9)\move(0 2)\rlvec(1
0)\move(0 7)\rlvec(1 0)\move(0 4)\rlvec(1 0)

\htext(0.5 1){$\overline{n}$} \htext(0.5
3){$n\!\!-\!\!1$}\htext(0.5 8){$1$}\htext(0.5 5.5){$\vdots$}
\end{texdraw}}
\savebox{\tmpfigd}{\begin{texdraw} \fontsize{7}{7}\selectfont
\drawdim em \setunitscale 0.8

\rlvec(1 0)\rlvec(0 9)\rlvec(-1 0)\rlvec(0 -9)\move(0 2)\rlvec(1
0)\move(0 7)\rlvec(1 0)\move(0 4)\rlvec(1 0)

\htext(0.5 1){$n$} \htext(0.5 3){$n\!\!-\!\!1$}\htext(0.5
8){$1$}\htext(0.5 5.5){$\vdots$}
\end{texdraw}}
\savebox{\tmpfige}{\begin{texdraw} \fontsize{7}{7}\selectfont
\drawdim em \setunitscale 0.85

\rlvec(2 0)\rlvec(0 2)\rlvec(-2 0)\rlvec(0 -2)\rlvec(2 2)

\htext(0.5 1.5){$n\!\!-\!\!1$} \htext(1.5 0.5){$n$}
\end{texdraw}}
\savebox{\tmpfigf}{\begin{texdraw} \fontsize{7}{7}\selectfont
\drawdim em \setunitscale 0.85

\rlvec(2 0)\rlvec(0 2)\rlvec(-2 0)\rlvec(0 -2)\rlvec(2 2)

\htext(1.5 0.5){$n\!\!-\!\!1$} \htext(0.5 1.5){$n$}
\end{texdraw}}
\savebox{\tmpfigg}{\begin{texdraw}\fontsize{7}{7}\selectfont
\drawdim em \setunitscale 0.85

\rlvec(1 0)\rlvec(0 2)\rlvec(-1 0)\rlvec(0 -2)

\htext(0.5 1){$\overline{n\!\!-\!\!1}$}
\end{texdraw}}
\savebox{\tmpfigh}{\begin{texdraw} \fontsize{7}{7}\selectfont
\drawdim em \setunitscale 0.85

\rlvec(2 0)\rlvec(0 4)\rlvec(-2 0)\rlvec(0 -4)\rlvec(2 2)\rlvec(-2
0)

\htext(0.5 1.5){$n\!\!-\!\!1$}\htext(1.5 0.5){$n$}\htext(1
3){$n\!\!-\!\!2$}
\end{texdraw}}
\savebox{\tmpfigi}{\begin{texdraw}\fontsize{7}{7}\selectfont
\drawdim em \setunitscale 0.85

\rlvec(2 0)\rlvec(0 4)\rlvec(-2 0)\rlvec(0 -4)\rlvec(2 2)\rlvec(-2
0)

\htext(1.5 0.5){$n\!\!-\!\!1$}\htext(0.5 1.5){$n$}\htext(1
3){$n\!\!-\!\!2$}
\end{texdraw}}
\savebox{\tmpfigj}{\begin{texdraw} \fontsize{7}{7}\selectfont
\drawdim em \setunitscale 0.85

\rlvec(1 0)\rlvec(0 2)\rlvec(-1 0)\rlvec(0 -2)

\htext(0.5 1){$\overline{n\!\!-\!\!2}$}
\end{texdraw}}
\savebox{\tmpfigk}{\begin{texdraw} \fontsize{7}{7}\selectfont
\drawdim em \setunitscale 0.85

\rlvec(2 0)\rlvec(0 4)\rlvec(-2 0)\rlvec(0 -4)\move(0 2)\rlvec(2
0)

\htext(1 1){$j\!\!+\!\!1$}\htext(1 3){$j$}
\end{texdraw}}
\savebox{\tmpfigl}{\begin{texdraw} \fontsize{7}{7}\selectfont
\drawdim em \setunitscale 0.85

\rlvec(1 0)\rlvec(0 2)\rlvec(-1 0)\rlvec(0 -2)

\htext(0.5 1){$\overline{j}$}
\end{texdraw}}
\savebox{\tmpfigm}{\begin{texdraw} \fontsize{7}{7}\selectfont
\drawdim em \setunitscale 0.85

\rlvec(2 0)\rlvec(0 2)\rlvec(-2 0)\rlvec(0 -2)\rlvec(2 2)

\htext(0.5 1.5){$1$}
\end{texdraw}}
\savebox{\tmpfign}{\begin{texdraw} \fontsize{7}{7}\selectfont
\drawdim em \setunitscale 0.85

\rlvec(2 0)\rlvec(0 2)\rlvec(-2 0)\rlvec(0 -2)\rlvec(2 2)

\htext(1.5 0.5){$1$}
\end{texdraw}}
\savebox{\tmpfigo}{\begin{texdraw}\fontsize{7}{7}\selectfont
\drawdim em \setunitscale 0.85

\rlvec(1 0)\rlvec(0 2)\rlvec(-1 0)\rlvec(0 -2)

\htext(0.5 1){$\overline{1}$}
\end{texdraw}}
\savebox{\tmpfigp}{\begin{texdraw}\fontsize{7}{7}\selectfont
\drawdim em \setunitscale 0.85

\rlvec(2 0)\rlvec(0 2)\rlvec(-2 -2)\move(0 0)\rlvec(0 2)\rlvec(2
0)

\htext(1.5 0.5){$n$}
\end{texdraw}}
\savebox{\tmpfigq}{\begin{texdraw} \fontsize{7}{7}\selectfont
\drawdim em \setunitscale 0.85

\rlvec(1 0)\rlvec(0 2)\rlvec(-1 0)\rlvec(0 -2)

\htext(0.5 1){$\overline{n}$}
\end{texdraw}}

\vskip 3mm {\bf Step 2.}\, Now, consider the ${\overline
Y_{\la_{n-1}}}$ and ${\overline Y_{\la_n}}$ for $\frak g=D_n$.  If
$\overline{Y}_{\la_{n-1}}=H_{\lambda_{n-1}}=\raisebox{-0.25\height}{\usebox{\tmpfiga}}$
and
$\overline{Y}_{\la_{n}}=H_{\lambda_n}=\raisebox{-0.25\height}{\usebox{\tmpfigb}}$\,,
they correspond to the tableaux
$$
\raisebox{-0.25\height}{\usebox{\tmpfigc}}\quad\text{and}\quad
\raisebox{-0.25\height}{\usebox{\tmpfigd}}\,,\quad\text{respectively}.
$$

Suppose that there are some added blocks above $H_{\lambda_{n-1}}$
and $H_{\lambda_n}$ in ${\overline Y_{\la_{n-1}}}$ and ${\overline
Y_{\la_n}}$. Then the columns having added blocks above
$H_{\lambda_{n-1}}$ or $H_{\lambda_n}$ correspond to boxes with
entries as follows.
$$\begin{array}{ll}
\raisebox{-0.25\height}{\usebox{\tmpfige}}\,\,\,\text{or}\,\,\,
\raisebox{-0.25\height}{\usebox{\tmpfigf}}\,\,\,\longmapsto
\,\,\,\raisebox{-0.25\height}{\usebox{\tmpfigg}} &\qquad\quad
\raisebox{-0.25\height}{\usebox{\tmpfigm}}\,\,\,\text{or}\,\,\,
\raisebox{-0.25\height}{\usebox{\tmpfign}}\,\,\,\longmapsto
\,\,\,\raisebox{-0.25\height}{\usebox{\tmpfigo}} \\ \\
\raisebox{-0.25\height}{\usebox{\tmpfigh}}\,\,\,\text{or}\,\,\,
\raisebox{-0.25\height}{\usebox{\tmpfigi}}\,\,\,\longmapsto
\,\,\,\raisebox{-0.25\height}{\usebox{\tmpfigj}} &\qquad\quad
\raisebox{-0.25\height}{\usebox{\tmpfigk}}\,\,\,\longmapsto
\,\,\,\raisebox{-0.25\height}{\usebox{\tmpfigl}}\,\,\,(2\le j\le
n-3)
\end{array}$$
Moreover, if a column whose right column is a full column is
$\raisebox{-0.25\height}{\usebox{\tmpfigp}}$\,, a part of
ground-state wall, then it corresponds to
$\raisebox{-0.25\height}{\usebox{\tmpfigq}}$\,. Finally, if
$k$-many columns with blocks above $H_{\lambda_{n-1}}$ or
$H_{\lambda_n}$ in ${\overline Y_{\lambda_{n-1}}}$ and ${\overline
Y_{\lambda_n}}$ correspond to $k$-many boxes with entries
$i_1,\cdots,i_k\in\{\overline{n},\cdots,\overline{1}\}$, then
($n-k$)-many columns with no block above $H_{\lambda_{n-1}}$ or
$H_{\lambda_n}$ correspond to boxes with entries on
$\{1,\cdots,n\}\setminus\{\overline{i_1},\cdots,\overline{i_k}\}$.
Now, ${\overline Y}_{\la_{n-1}}$ and ${\overline Y}_{\la_{n}}$
correspond to tableaux of half column of length $n$ obtained by
stacking $T^a_{\la_{n-1}}$ and $T^a_{\la_n}$ ($a=1,\cdots,n$) for
the entries of $T^a_{\la_{n-1}}$ and $T^a_{\la_n}$ to strictly
increase down the columns, respectively.


%
\savebox{\spinaf}{\begin{texdraw}\fontsize{7}{7}\selectfont
\drawdim em \setunitscale 0.85 \move(-1 0)\lvec(8 0)\rlvec(0
7)\move(8 2)\rlvec(-2 -2)\rlvec(0 2)\rlvec(-2 -2)\rlvec(0
2)\rlvec(-2 -2)\rlvec(0 2)\rlvec(-2 -2)\rlvec(0 2)\move (8
6)\rlvec(-2 -2)\rlvec(2 0)\move (8 2)\rlvec(-2 0)\rlvec(0
2)\move(6 2)\rlvec(-6 0)\move(8 6)\rlvec(-2 0)\rlvec(0 -2)
\htext(7.5 4.5){$1$}\htext(7 3){$2$} \htext(1.5
0.5){$3$}\htext(3.5 0.5){$4$}\htext(4.5 1.5){$4$} \htext(5.5
0.5){$3$} \htext(7.5 0.5){$4$}\htext(6.5 1.5){$3$}
\end{texdraw}}
\savebox{\spinag}{\begin{texdraw}\fontsize{7}{7}\selectfont
\drawdim em \setunitscale 0.85

\lvec(1 0)\rlvec(0 2)\rlvec(-1 0)\rlvec(0 -2)

\htext(0.5 1){$\overline{1}$}
\end{texdraw}}
\savebox{\spinah}{\begin{texdraw}\fontsize{7}{7}\selectfont
\drawdim em \setunitscale 0.85

\lvec(1 0)\rlvec(0 2)\rlvec(-1 0)\rlvec(0 -2)

\htext(0.5 1){$\overline{3}$}
\end{texdraw}}
\savebox{\spinai}{\begin{texdraw}\fontsize{7}{7}\selectfont
\drawdim em \setunitscale 0.85

\lvec(1 0)\rlvec(0 2)\rlvec(-1 0)\rlvec(0 -2)

\htext(0.5 1){$\overline{4}$}
\end{texdraw}}
\savebox{\spinaj}{\begin{texdraw}\fontsize{7}{7}\selectfont
\drawdim em \setunitscale 0.85

\lvec(2 0)\rlvec(0 2)\rlvec(-2 0)\rlvec(0 -2)\rlvec(2 2)

\htext(1 0.5){$4$}
\end{texdraw}}
\savebox{\spinbf}{\begin{texdraw}\fontsize{7}{7}\selectfont
\drawdim em \setunitscale 0.8 \lvec(0 8)\rlvec(1 0)\rlvec(0
-8)\rlvec(-1 0)\move(0 2)\rlvec(1 0)\move(0 4)\rlvec(1 0)\move(0
6)\rlvec(1 0) \htext(0.5 1){$\bar{1}$}\htext(0.5 3){$\bar{3}$}
\htext(0.5 5){$\bar{4}$}\htext(0.5 7){$2$}
\end{texdraw}}

\vskip 3mm
\begin{example} Let $\frak g=D_4$, $\la=\la_3$ and
$$Y=(y^1_{\la_3},y^2_{\la_3},y^3_{\la_3},y^4_{\la_3})=\raisebox{-0.4\height}{\usebox{\spinaf}}\,.$$
Then $y^1_{\la_3}$ and $y^2_{\la_3}$ correspond to
$\raisebox{-0.3\height}{\usebox{\spinag}}$ and
$\raisebox{-0.3\height}{\usebox{\spinah}}$\,, respectively.
Moreover, since $y^3_{\la_3}$ is
$\raisebox{-0.3\height}{\usebox{\spinaj}}$ and the its right
column $y^2_{\la_3}$ is a full column, it corresponds to
$\raisebox{-0.3\height}{\usebox{\spinai}}$\,. Therefore, we have
$$T_Y=\raisebox{-0.3\height}{\usebox{\spinbf}}\,.$$
\end{example}

\vskip 3mm  Then we have the following realization of crystal
bases for the classical Lie algebras $\frak g=A_n,C_n$ and $D_n$.

\vskip 3mm At first, we define $T(\la)$ for $\frak g=A_n$ and
$C_n$ as follows:

\begin{itemize}
\item [1)] $\frak g=A_n$ : the set of tableaux of shape $GRY^{-1}(\la)$
with entries $\{1,\cdots, n\}$ such that
\begin{itemize}
\item [(a)] the entries of $T$ weakly increase along the rows,
\item [(b)] the entries of $T$ strictly increase down the columns,
\end{itemize}

\item [2)] $\frak g=C_n$ : the set of tableaux of shape
$GRY^{-1}(\la)$ with entries $\{i,\overline{i}\,|i=1,\cdots,n\}$
such that
\begin{itemize}
\item [(a)] the entries of $T$ weakly increase along the rows,
\item [(b)] the entries of $T$ strictly increase down the columns,
\item [(c)] for each column $C$ of $T$, (1CC) holds,
\item [(d)] for each pair of adjacent columns $C$, $C'$ of $T$, (2CC-1)
and (2CC-2) holds.
\end{itemize}
\end{itemize}

\vskip 3mm
\begin{thm}
For a dominant integral weight $\lambda$, there exists a crystal
isomorphism for $U_q(\frak g)$-modules $\varphi:Y(\la)\rightarrow
T(\la)$ sending $Y$ to $T_Y$, where $\frak g=A_n$ and $C_n$.
\end{thm}

\begin{proof} It is similar to the proof of Theorem Theorem
\ref{thm:B_n}, we omit it.
\end{proof}

\vskip 3mm For the type of $\frak g=D_n$, we have another
conditions. Let $C$ and $C'$ be adjacent columns of length $N$ and
$M$ ($1\le M\le N\le n$) consisting of the entries (reading from
bottom to top) $i_1$,$i_2$,$\cdots$,$i_M$ and
$j_1$,$j_2$,$\cdots$,$j_N$, respectively. (Note that $C$ can be a
half-column.)

\vskip 2mm
\begin{defi}
For $1\le a<n$, we say that $C$ and $C'$ is in the
$a$-odd-configuration (resp. $a$-even-configuration) if there
exist $M\ge p\ge q>r\ge s\ge 1$ such that $q-r+1$ is odd (resp.
even) and
$$\aligned &(i_p,j_q,i_r,j_s)=(a,n,\overline{n},\overline{a})\,\,\,\text{or}
\,\,\,(a,\overline{n},n,\overline{a})\\
(\text{resp.}\,\,\,\,
&(i_p,j_q,i_r,j_s)=(a,n,n,\overline{a})\,\,\,\text{or}
\,\,\,(a,\overline{n},\overline{n},\overline{a})).
\endaligned$$
\end{defi}

\vskip 3mm Now,  for $C$ and $C'$ in the $a$-odd-configuration or
$a$-even-configuration, we define
$$q(a;C,C')=p-s.$$

\vskip 3mm {\bf (2CC-3)} Given adjacent two columns $C$ and $C'$,
we say that $C$ and $C'$ satisfy the {\it third two-column
condition} (2CC-3) if for every $a$-odd-configuration or
$a$-even-configuration, $q(a;C,C')<n-a$.

\vskip 3mm Now, we can describe the crystal graph $B(\la)$ for
$\frak g=D_n$ as follows.

\vskip 2mm \noindent We define $T(\la)$ for $\frak g=D_n$ by the
set of tableaux of shape $GRY^{-1}(\la)$ with entries
$\{i,\overline{i}\,|i=1,\cdots,n\}$ such that
\begin{enumerate}
\item the entries of $T$ weakly increase along the rows,
\item the entries of $T$ strictly increase down the columns, but
the element $\overline{n}$ (resp. $n$) can appear below the
element $n$ (resp. $\overline{n}$),
\item for a half column $C$ of $T$, $i$ and $\overline{i}$ can not
appear at the same time,
\item for each column $C$ of $T$, (1CC) holds,
\item for each pair of adjacent columns $C$, $C'$ of $T$,
(2CC-1), (2CC-2) and (2CC-3) holds.
\end{enumerate}

\vskip 3mm
\begin{thm}
For a dominant integral weight $\la$, there is a crystal
isomorphism for $U_q(D_n)$-modules $\varphi: Y(\la)\rightarrow
T(\la)$ sending $Y$ to $T_Y$.
\end{thm}

\begin{proof} It is similar to the argument of the proof of
Theorem \ref{thm:B_n}, we will omit it.
\end{proof}

%
\savebox{\tmpfiga}{
\begin{texdraw}
\fontsize{7}{7}\selectfont \drawdim em \setunitscale 0.8
\nc{\dtri}{ \bsegment \lvec(-2 0) \lvec(-2 2)\lvec(0 2)\lvec(0
0)\ifill f:0.7 \esegment } \nc{\dtrii}{ \bsegment \lvec(-2 0)
\lvec(-2 1)\lvec(0 1)\lvec(0 0)\ifill f:0.7 \esegment }

\bsegment \setgray 0.6 \move(-1 0) \lvec(16 0)\lvec(16 15)
\move(16 12)\rlvec(-2 0)\rlvec(0 -12) \move(16 10)\rlvec(-4
0)\rlvec(0 -10) \move(16 8)\rlvec(-6 0)\rlvec(0 -8)\move(16
6)\rlvec(-10 0)\rlvec(0 -6)\move(8 6)\rlvec(0 -6)\move(16
4)\rlvec(-10 0)\move(16 2)\rlvec(-16 0)\rlvec(0 -2) \move(16
10)\rlvec(-2 -2)\move(14 10)\rlvec(-2 -2) \move(16 6)\rlvec(-2
-2)\move(14 6)\rlvec(-2 -2) \move(12 6)\rlvec(-2 -2)\move(10
6)\rlvec(-2 -2)\move(8 6)\rlvec(-2 -2) \move(16 2)\rlvec(-2
-2)\move(14 2)\rlvec(-2 -2)\move(12 2)\rlvec(-2 -2) \move(10
2)\rlvec(-2 -2)\move(8 2)\rlvec(-2 -2)\move(6 2)\rlvec(-2
-2)\move(4 2)\rlvec(-2 -2)\move(2 2)\rlvec(-2 -2)\move(6
0)\rlvec(0 2)\move(4 0)\rlvec(0 2)\move(2 0)\rlvec(0 2)

\move(16 10)\linewd 0.15 \setgray 0.4 \rlvec(-2 -2)\rlvec(-2 0)
\rlvec(0 -2) \rlvec(-2 0)\rlvec(-2 -2)\rlvec(0 -2) \rlvec(-2 -2)
\rlvec(0 2)\rlvec(-2 -2) \rlvec(0 2)\rlvec(-2 -2) \rlvec(0
2)\rlvec(-2 -2)

\htext(15 11){$2$} \htext(15.5 8.5){$3$}\htext(14.5
9.5){$4$}\htext(12.5 9.5){$3$} \htext(15 7){$2$} \htext(13
7){$2$}\htext(11 7){$2$} \htext(15.5 4.5){$1$}\htext(14.5
5.5){$0$} \htext(13.5 4.5){$0$}\htext(12.5 5.5){$1$}\htext(11.5
4.5){$1$}\htext(10.5 5.5){$0$}\htext(9.5 4.5){$0$}\htext(7.5
4.5){$1$} \htext(15 3){$2$}\htext(13 3){$2$}\htext(11
3){$2$}\htext(9 3){$2$}\htext(7 3){$2$} \htext(15.5
0.5){$3$}\htext(14.5 1.5){$4$}\htext(13.5 0.5){$4$} \htext(12.5
1.5){$3$}\htext(11.5 0.5){$3$} \htext(10.5 1.5){$4$} \htext(9.5
0.5){$4$}\htext(8.5 1.5){$3$} \htext(7.5 0.5){$3$}\htext(6.5
1.5){$4$}  \htext(5.5 0.5){$4$}\htext(4.5 1.5){$3$}  \htext(3.5
0.5){$3$}\htext(1.5 0.5){$4$} \esegment
\end{texdraw}}
\savebox{\tmpfigb}{\begin{texdraw} \fontsize{7}{7}\selectfont
\drawdim em \setunitscale 0.8 \nc{\dtri}{ \bsegment \lvec(-2 0)
\lvec(-2 2)\lvec(0 2)\lvec(0 0)\ifill f:0.7 \esegment }
\nc{\dtrii}{ \bsegment \lvec(-2 0) \lvec(-2 1)\lvec(0 1)\lvec(0
0)\ifill f:0.7 \esegment } \rlvec(3 0)\rlvec(0 8) \rlvec(-3
0)\rlvec(0 -8)\move(2 8)\rlvec(0 -8) \move(3 6)\rlvec(-3 0)
\move(3 4)\rlvec(-3 0)\move(3 2)\rlvec(-3 0) \htext(2.5
7){$2$}\htext(2.5 5){$4$}\htext(2.5 3){$\bar{3}$}\htext(2.5
1){$\bar{1}$} \htext(1 7){$1$} \htext(1 5){$3$}\htext(1 3){$4$}
\htext(1 1){$\bar{2}$}
\end{texdraw}}

\vskip 3mm
\begin{example}
If $\frak g=D_4$, $\lambda=\omega_4+\la_4=3\la_4$ and
\begin{center}
$Y=\raisebox{-0.4\height}{\usebox{\tmpfiga}}\in Y(\la)$, then
$\varphi(Y)=\raisebox{-0.4\height}{\usebox{\tmpfigb}}$\,.
\end{center}
\end{example}

\vskip 3mm
\begin{remark} We know that the crystal graph $B(\la)$ for the quantum classical
algebra $U_q(\frak g)$ appears as a connected component in the
crystal graph  $B(\La)$ of basic representation for the quantum
affine algebra $U_q(\widehat{\frak g})$ without $0$-arrows.
Moreover, there are many connected components of $B(\La)$ without
$0$-arrows that are isomorphic to $B(\la)$. Therefore, we can
obtain a lot of Young wall realizations of the crystal graph
$B(\la)$. In particular, if we choose the Young wall realization
of $B(\la)$ corresponding to the connected components having the
largest number of blocks, then applying a similar map to the map
$\varphi$, we can obtain the tableau realization of $B(\la)$ given
by Kashiwara-Nakashima (We shall not describe this Young wall
realization and the map similar to $\varphi$ in this paper)
\end{remark}

\vskip 3mm  Finally, we will close this section to give a 1-1
correspondence between our new realization and the realization of
Kashiwara and Nakashima for $\frak g=A_n,C_n,B_n$ and $D_n$.
Moreover, we will prove that it is just a crystal isomorphism for
$U_q(\frak g)$.

\begin{thm} Let $U_q(\frak g)$ $(\frak g=A_n,C_n,B_n,D_n)$
be a quantum classical Lie algebra and let $\la$ be a dominant
integral weight. Then there is a crystal isomorphism $\psi:T(\la)
\rightarrow B(\la)$ for $U_q(\frak g)$-modules.
\end{thm}

\begin{proof} Let $T$ be a tableau of $T(\la)$ consisting of
columns $T_1$, $T_2$, $\cdots$, $T_p$ reading from right to left.
We define $\psi(T)$ by
\begin{center}
$T_1\leftarrow T_2\leftarrow \cdots \leftarrow T_p.$
\end{center}
Here, the notation $T\leftarrow U$ for the tableaux $T, U$
consisting of one column is referred to \cite{KS}. Then by Theorem
4.14 of \cite{KS}, it is clear that the map $\psi$ is a crystal
isomorphism.
\end{proof}

\savebox{\tmpfiga}{\begin{texdraw}
\fontsize{7}{7}\selectfont\drawdim em \setunitscale 0.8

\rlvec(6 0)\rlvec(0 6)\rlvec(-2 0)\rlvec(0 -6)\move(0 0)\rlvec(0
2) \rlvec(6 0) \move(2 0)\rlvec(0 4)\rlvec(4 0)

\htext(1 1){$2$}\htext(3 1){$3$}\htext(3 3){$1$}\htext(5
1){$\overline{3}$}\htext(5 3){$3$}\htext(5 5){$2$}
\end{texdraw}}
\savebox{\tmpfigb}{\begin{texdraw}
\fontsize{7}{7}\selectfont\drawdim em \setunitscale 0.8

\rlvec(0 6)\rlvec(2 0)\rlvec(0 -6)\rlvec(-2 0)\move(0 2)\rlvec(2
0)\move(0 4)\rlvec(2 0)

\htext(1 1){$\overline{3}$}\htext(1 3){$3$}\htext(1 5){$2$}
\end{texdraw}}
\savebox{\tmpfigc}{\begin{texdraw}
\fontsize{7}{7}\selectfont\drawdim em \setunitscale 0.8

\rlvec(0 4)\rlvec(2 0)\rlvec(0 -4)\rlvec(-2 0)\move(0 2)\rlvec(2
0)

\htext(1 1){$3$}\htext(1 3){$1$}
\end{texdraw}}
\savebox{\tmpfigd}{\begin{texdraw}
\fontsize{7}{7}\selectfont\drawdim em \setunitscale 0.8

\rlvec(0 2)\rlvec(2 0)\rlvec(0 -2)\rlvec(-2 0)

\htext(1 1){$2$}
\end{texdraw}}
\savebox{\tmpfige}{\begin{texdraw}
\fontsize{7}{7}\selectfont\drawdim em \setunitscale 0.8

\rlvec(0 6)\rlvec(6 0)\rlvec(0 -2)\rlvec(-6 0)\move(0 0)\rlvec(2
0)\rlvec(0 6)\move(0 2)\rlvec(4 0)\rlvec(0 4)

\htext(1 3){$2$}\htext(1 1){$3$}\htext(1 5){$1$}\htext(3
5){$1$}\htext(5 5){$\overline{1}$}\htext(3 3){$2$}
\end{texdraw}}
\savebox{\tmpfigf}{\begin{texdraw}
\fontsize{7}{7}\selectfont\drawdim em \setunitscale 0.8

\rlvec(5 0)\rlvec(0 8)\rlvec(-1 0)\rlvec(0 -8)\move(0 0)\rlvec(0
2) \rlvec(5 0) \move(2 0)\rlvec(0 4)\rlvec(3 0)\move(4 6)\rlvec(1
0)

\htext(1 1){$1$}\htext(3 1){$\overline{3}$}\htext(3
3){$\overline{4}$}\htext(4.5 1){$\overline{1}$}\htext(4.5
3){$\overline{4}$}\htext(4.5 5){$3$}\htext(4.5 7){$2$}
\end{texdraw}}
\savebox{\tmpfigg}{\begin{texdraw}
\fontsize{7}{7}\selectfont\drawdim em \setunitscale 0.8

\rlvec(0 8)\rlvec(1 0)\rlvec(0 -8)\rlvec(-1 0)\move(0 2)\rlvec(1
0)\move(0 4)\rlvec(1 0)\move(0 6)\rlvec(1 0)

\htext(0.5 1){$\overline{1}$}\htext(0.5
3){$\overline{4}$}\htext(0.5 5){$3$}\htext(0.5 7){$2$}
\end{texdraw}}
\savebox{\tmpfigh}{\begin{texdraw}
\fontsize{7}{7}\selectfont\drawdim em \setunitscale 0.8

\rlvec(0 4)\rlvec(2 0)\rlvec(0 -4)\rlvec(-2 0)\move(0 2)\rlvec(2
0)

\htext(1 1){$\overline{3}$}\htext(1 3){$\overline{4}$}
\end{texdraw}}
\savebox{\tmpfigi}{\begin{texdraw}
\fontsize{7}{7}\selectfont\drawdim em \setunitscale 0.8

\rlvec(0 2)\rlvec(2 0)\rlvec(0 -2)\rlvec(-2 0)

\htext(1 1){$1$}
\end{texdraw}}
\savebox{\tmpfigj}{\begin{texdraw}
\fontsize{7}{7}\selectfont\drawdim em \setunitscale 0.8

\rlvec(0 8)\rlvec(5 0)\rlvec(0 -2)\rlvec(-5 0)\move(0 0)\rlvec(1
0)\rlvec(0 8)\move(0 2)\rlvec(1 0)\move(0 4)\rlvec(3 0)\rlvec(0 4)

\htext(0.5 1){$\overline{3}$}\htext(0.5
3){$\overline{4}$}\htext(0.5 5){$2$}\htext(0.5 7){$1$} \htext(2
5){$\overline{2}$}\htext(2 7){$2$} \htext(4 7){$\overline{4}$}
\end{texdraw}}

\vskip 3mm
\begin{example}
(a) Let $\frak g=B_4$, $\la=\omega_1+\omega_2+\omega_3$ and
$$T=\raisebox{-0.4\height}{\usebox{\tmpfiga}}\,.$$
Then
$$\psi(T)=\raisebox{-0.4\height}{\usebox{\tmpfigb}} \,\,\,\leftarrow\,\,\,
\raisebox{-0.4\height}{\usebox{\tmpfigc}}\,\,\,\leftarrow\,\,\,
\raisebox{-0.3\height}{\usebox{\tmpfigd}}=\raisebox{-0.4\height}{\usebox{\tmpfige}}\,.$$

(b) Let $\frak g=D_4$, $\la=\omega_1+\omega_2+\la_4$ and
$$T=\raisebox{-0.4\height}{\usebox{\tmpfigf}}\,.$$
Then
$$\psi(T)=\raisebox{-0.4\height}{\usebox{\tmpfigg}} \,\,\,\leftarrow\,\,\,
\raisebox{-0.4\height}{\usebox{\tmpfigh}}\,\,\,\leftarrow\,\,\,
\raisebox{-0.3\height}{\usebox{\tmpfigi}}=\raisebox{-0.4\height}{\usebox{\tmpfigj}}\,.$$
\end{example}

\vskip 10mm \noindent {\bf Acknowledgments.}\, The authors are
grateful to professor S.-J. Kang for valuable discussions. This
research was supported by KOSEF Grant \# 98-0701-01-5-L and BK21
Mathematical Sciences Division, Seoul National University.

\vskip 15mm


\begin{thebibliography}{1}

\bibitem{Con}
C. De Concini, \emph{Symplectic  standard tableaux}, Adv. Math.,
\textbf{34} (1979), 1--27.

\bibitem{Drin}
V. G. Drinfeld, \emph{Hopf algebras and the quantum Yang-Baxter
equation}, Soviet Math. Dokl., \textbf{32} (1985), 254--258.

\bibitem{Fulton}
W.~Fulton, \emph{Young Tableaux {\rm :} with applications to
representation theory and geometry}, Cambridge University Press,
1997.

\bibitem{HK}
J.~Hong, S.-J. Kang, \emph{Introduction to Quantum Groups and
Crystal Bases}, Graduate Studies in Mathematics \textbf{42}, Amer.
Math. Soc.

\bibitem{Jim}
M.~Jimbo, \emph{A $q$-difference analogue of $U(\frak g)$ and the
Yang-Baxter equation}, Lett. Math. Phys. \textbf{10} (1985),
63--69.

\bibitem{JMMO}
M.~Jimbo, K.~C. Misra, T.~Miwa, M.~Okado, \emph{Combinatorics of
representations of ${U}_q(\widehat{\mathfrak{sl}}(n))$ at $q=0$},
Comm. Math. Phys. \textbf{136} (1991), 543--566.

\bibitem{Kac90}
V.~Kac, \emph{Infinite Dimensional Lie Algebras}, Cambridge
University Press, 3rd ed., 1990.

\bibitem{Kang}
S.-J. Kang, \emph{Crystal bases for quantum affine algebras and
combinatorics of Young walls}, to appear in Proc. London Math.
Soc.

\bibitem{KKLS}
S.-J. Kang, J.-A. Kim, H. Lee, D.-U. Shin \emph{Young wall
realization of crystal bases for classical Lie algebras},
submitted

\bibitem{KM}
S.-J. Kang, K.~C. Misra, \emph{Crystal bases and tensor product
decompositions of $U_q(G_2)$-modules}, J. Algebra 163 (1994),
675--691.

\bibitem{Kas90}
M.~Kashiwara, \emph{Crystalizing the $q$-analogue of universal
enveloping
  algebras}, Comm. Math. Phys. \textbf{133} (1990), 249--260.

\bibitem{Kas91}
\bysame, \emph{On crystal bases of the $q$-analogue of universal
enveloping algebras}, Duke Math. J. \textbf{63} (1991), 465--516.

\bibitem{Kas96}
\bysame, \emph{Similarity of crystal bases}, in Lie algebras and
the representations, Seoul(1995), Contemp. Math. \textbf{194}
(1996), American Math. Soc., 177-186.

\bibitem{KasNak}
M.~Kashiwara, T.~Nakashima, \emph{Crystal graphs for
representations of the $q$-analogue of classical Lie algebras}, J.
Algebra \textbf{165} (1994), 295-345.

\bibitem{KS}
J.-A. Kim, D.-U. Shin, \emph{Insertion scheme for the crystal of
the classical Lie algebras}, submitted

\bibitem{Li1}
P.~Littelmann, \emph{Crystal graphs and Young tableaux}, J.
Algebra \textbf{175} (1995), 65--87.

\bibitem{Li2}
\bysame, \emph{ Paths and root operators in representation
theory}, Ann. of Math. \textbf{142} (1995), 499--525.

\bibitem{Li3}
\bysame, \emph{ A Littlewood-Richardson rule for symmetrizable
Kac-Moody algebras}, Invent. Math. 116 (1994), 329--346.

\bibitem{Mac}
I. G. Macdonald, \emph{Symmetric Functions and Hall Polynomials},
Oxford University Press, Oxford, 2nd ed., 1995.

\bibitem{Nak}
T.~Nakashima, \emph{Crystal base and a generalization of the
Littlewood-Richardson rule for classical Lie algebras}, Comm.
Math. Phys. \textbf{154} (1993), 215--243.

\bibitem{Sheats}
J. T. Sheats, \emph{A symplectic jeu de taquin bijection beyween
the tableaux of King and of De Concini}, Trans. Amer. Math. Soc.
\textbf{351} (1999), 3569--3607.

\end{thebibliography}

\providecommand{\bysame}{\leavevmode\hbox
to5em{\hrulefill}\thinspace}

\end{document}